\documentclass[twoside,11pt]{article}

\usepackage{blindtext}

\usepackage[abbrvbib, preprint]{jmlr2e}

\usepackage{lastpage}

\usepackage{graphicx}
\usepackage{caption}
\usepackage{comment}
\usepackage{amsmath} 
\usepackage{amssymb}  
\usepackage{bm}
\usepackage{lipsum}
\usepackage[linesnumbered,ruled,vlined]{algorithm2e}
\usepackage{algorithmic}
\usepackage{color}
\usepackage{array}
\usepackage{hyperref}
\usepackage{upgreek}
\usepackage{todonotes}
\usepackage{cite}
\usepackage[capitalize,noabbrev]{cleveref}

\usepackage{multirow}

\usepackage{thmtools}
\usepackage{thm-restate}
\usepackage{float}

\usepackage{color}
\usepackage{colortbl}
\definecolor{lightgray}{gray}{0.9}

\usepackage[capitalize,noabbrev]{cleveref}

\usepackage{microtype}
\usepackage{graphicx}
\usepackage{subfigure}
\usepackage{multirow}
\usepackage{makecell}

\usepackage{tikz}
\usetikzlibrary{shapes.geometric, arrows, positioning, fit, calc}
\tikzset{
	box/.style  = {draw, rectangle, minimum width=3.0cm, minimum height=1.7cm, text centered, text width=3.5cm,   
		font=\normalsize},
	smallbox/.style  = {draw, rectangle, minimum width=2.0cm, minimum height=1.0cm, text centered, text width=1.8cm,   
		font=\normalsize},
}
\tikzstyle{arrow} = [->, >=stealth, -triangle 60]
\tikzstyle{botharrow} = [<->, >=stealth, {triangle 60}-{triangle 60}]

\newcommand{\sech}{\operatorname{sech}}
\newcommand{\csch}{\operatorname{csch}}
\newcommand{\sinhc}{\operatorname{sinhc}}
\newcommand{\tanhc}{\operatorname{tanhc}}
\newcommand{\cothc}{\operatorname{cothc}}
\newcommand{\cschc}{\operatorname{cschc}}

\newcommand{\nagc}{NAG-\texttt{C}}
\newcommand{\nagsc}{NAG-\texttt{SC}}
\newcommand{\nagcsp}{NAG-\texttt{C} }
\newcommand{\nagscsp}{NAG-\texttt{SC} }
\newcommand{\bft}{\mathbf{t}}
\newcommand{\bfk}{\mathbf{k}}
\newcommand{\tst}{0}
\newcommand{\bfT}{\mathbf{T}}

\allowdisplaybreaks

\ShortHeadings{Unifying NAG for Convex and Strongly Convex Objective Functions}{Kim and Yang}
\firstpageno{1}

\begin{document}

\title{Unifying Nesterov’s Accelerated Gradient Methods for Convex and Strongly Convex Objective Functions: From Continuous-Time Dynamics to Discrete-Time Algorithms}

\author{\name Jungbin Kim \email kjb2952@snu.ac.kr \\
       \addr Department of Electrical and Computer Engineering\\
       Seoul National University\\
       Seoul 08826, Korea
       \AND
       \name Insoon Yang \email insoonyang@snu.ac.kr \\
       \addr Department of Electrical and Computer Engineering\\
       Seoul National University\\
       Seoul 08826, Korea}

\editor{My editor}

\maketitle

\begin{abstract}
Although Nesterov's accelerated gradient (NAG) methods have been studied from various perspectives, it remains unclear why the most popular forms of NAG must handle convex and strongly convex objective functions separately. Motivated by this inconsistency, we propose an NAG method that unifies the existing ones for the convex and strongly convex cases. We first design a Lagrangian function that continuously extends the \emph{first} Bregman Lagrangian to the strongly convex setting. As a specific case of the Euler--Lagrange equation for this Lagrangian, we derive an ordinary differential equation (ODE) model, which we call the \emph{unified NAG ODE}, that bridges the gap between the ODEs that model NAG for convex and strongly convex objective functions. We then design the \emph{unified NAG}, a novel momentum method whereby the continuous-time limit corresponds to the unified ODE. The coefficients and the convergence rates of the unified NAG and unified ODE are continuous in the strong convexity parameter $\mu$ on $[0, +\infty)$. 
Unlike the existing popular algorithm and ODE for strongly convex objective functions, the unified NAG and the unified NAG ODE always have superior convergence guarantees compared to the known algorithms and ODEs for non-strongly convex objective functions. This property is beneficial in practical perspective when considering strongly convex objective functions with small $\mu$. 
Furthermore, we extend our unified dynamics and algorithms to the higher-order setting. {Last but not least, we propose the \emph{unified NAG-G ODE}, a novel ODE model for minimizing the gradient norm of strongly convex objective functions. Our unified Lagrangian framework is crucial in the process of constructing this ODE. Fascinatingly, using our novel tool, called the \emph{differential kernel}, we observe that the unified NAG ODE and the unified NAG-G ODE have an anti-transpose relationship.}
\end{abstract}

\begin{keywords}
 Convex optimization, first-order methods, Nesterov acceleration
\end{keywords}

\section{Introduction}
\label{sec:introduction}

We consider the optimization problem
\begin{equation}
	\label{eq:problem}
	\min_{x \in \mathbb{R}^n} \; f(x),
\end{equation}
where $f:\mathbb{R}^n\rightarrow\mathbb{R}$ is a continuously differentiable function whose gradient is $L$-Lipschitz continuous. We assume that the objective function $f$ has a minimizer $x^{*}$. One of the most popular first-order method for solving this problem is gradient descent (GD):
\begin{equation}
	\label{eq:gd}
	x_{k+1}=x_k-s\nabla f\left(x_k\right)
\end{equation}
with the algorithmic stepsize $s>0$. When $f$ is convex, GD with $s\leq1/L$ achieves an $O(\|x_{0}-x^{*}\|^{2}/k)$ convergence rate \citep[see][Section~4.2]{daspremont2021acceleration}. When $f$ is $\mu$-strongly convex, GD with $s\leq1/L$ achieves an $O((1-\mu s)^k\|x_{0}-x^{*}\|^{2})$ convergence rate \citep[see][Section~4.5]{daspremont2021acceleration}.

\paragraph{Nesterov acceleration.}
A natural and important question is whether there are other first-order methods that outperform gradient descent. \citet{nesterov1983} proposed an accelerated gradient method that achieves a faster convergence rate compared to gradient descent. Given the initial point $x_0=z_0$, a general \emph{three-sequence scheme} for Nesterov's accelerated gradient (NAG) methods can be written as
\begin{subequations}
	\label{eq:nag}
	\begin{align}
		\label{eq:nag1}
		y_{k} & =x_{k}+\tau_{k}\left(z_{k}-x_{k}\right)\\
		\label{eq:nag2}
		x_{k+1} & =y_{k}-s\nabla f\left(y_{k}\right)\\
		\label{eq:nag3}
		z_{k+1} & =z_{k}+\delta_{k}\left(\mu y_{k}-\mu z_{k}-\nabla f\left(y_{k}\right)\right)
	\end{align}
\end{subequations}
with $s>0$, where the parameters $\tau_k$ and $\delta_k$ usually satisfy the \emph{collinearity condition}\footnote{This condition ensures that the points $x_{k}$, $x_{k+1}$, $z_{k+1}$ are collinear (see Section~\ref{sec:limiting_two}). Thus, one can write the updating rule for $y_k$ as $y_{k+1}=x_{k+1}+\beta_k(x_{k+1}-x_{k})$ for some $\beta_k \in \mathbb{R}$. This property provides a clear momentum effect: The point $y_{k+1}$ is defined by adding a momentum term $\beta_k\left(x_{k+1}-x_{k}\right)$ to the previous point $x_{k+1}$. This property is useful when generalizing NAG methods to handle non-smooth terms \citep[see][Algorithm~20]{daspremont2021acceleration}.}
\begin{equation}
	\label{eq:collinear}
	1-\mu\delta_{k}-(1/s-\mu)\tau_{k}\delta_{k}=0.
\end{equation}
In particular, for $\mu$-strongly (possibly with $\mu=0$) convex objective functions, Nesterov considered the following algorithm: Given an initial point $x_0=z_0\in\mathbb{R}^n$ and $\gamma_0>0$, the \emph{constant step scheme I} \citep[Equation~2.2.19]{nesterov2018lectures} (we will refer to this algorithm as the \emph{original NAG}) updates the iterates as
\begin{equation}
	\label{eq:original_nag}
	\begin{aligned}
		\gamma_{k+1} & =\left(1-\alpha_{k}\right)\gamma_{k}+\mu\alpha_{k}\\
		y_{k} & =\frac{1}{\gamma_{k}+\mu\alpha_{k}}\left(\alpha_{k}\gamma_{k}z_{k}+\gamma_{k+1}x_{k}\right)\\
		x_{k+1} & =y_{k}-s\nabla f\left(y_{k}\right)\\
		z_{k+1} & =\frac{1}{\gamma_{k+1}}\left(\left(1-\alpha_{k}\right)\gamma_{k}z_{k}+\mu\alpha_{k}y_{k}-\alpha_{k}\nabla f\left(y_{k}\right)\right),
	\end{aligned}
\end{equation}
where the sequence $(\alpha_k)_{k=0}^{\infty}$ in $(0,1)$ is inductively defined by the equation
\begin{equation}
	\label{eq:alpha_k_update}
	\frac{1}{s}\alpha_k^2 = \left(1-\alpha_k\right)\gamma_k + \mu \alpha_k.
\end{equation}
Using the estimate sequence technique, \citet[Theorem~2.2.1]{nesterov2018lectures} showed that the iterates of the original NAG \eqref{eq:original_nag} satisfy the inequality 
\begin{equation}
	\label{eq:original_nag_convergence}
	f\left(x_{k}\right)-f\left(x^{*}\right)\leq\left(\prod_{i=0}^{k-1}\left(1-\alpha_{i}\right)\right)\left(f\left(x_{0}\right)-f\left(x^{*}\right)+\frac{\gamma_{0}}{2}\left\Vert x_{0}-x^{*}\right\Vert ^{2}\right)
\end{equation} 
when $s\leq1/L$. Although the original NAG achieves a faster convergence rate than gradient descent, it is difficult to analyze this algorithm because it involves auxiliary sequences $\alpha_k$ and $\gamma_k$ which are defined inductively. However, when $\gamma_0=\mu$ (here we need $\mu>0$ because $\gamma_0>0$ is assumed), we simply have $\alpha_{k}=\sqrt{\mu s}$ and $\gamma_{k}=\mu$ for all $k\geq0$. In this case, the original NAG \eqref{eq:original_nag} can be expressed as the three-sequence scheme \eqref{eq:nag} with $\tau_{k}=\frac{\sqrt{\mu s}}{1+\sqrt{\mu s}}$ and $\delta_{k}=\sqrt{\frac{s}{\mu}}$: 
\begin{equation}
	\label{eq:nag-sc}
	\begin{aligned}
		y_{k} & =x_{k}+\frac{\sqrt{\mu s}}{1+\sqrt{\mu s}}\left(z_{k}-x_{k}\right)\\
		x_{k+1} & =y_{k}-s\nabla f\left(y_{k}\right)\\
		z_{k+1} & =z_{k}+\sqrt{\frac{s}{\mu}}\left(\mu y_{k}-\mu z_{k}-\nabla f\left(y_{k}\right)\right).
	\end{aligned}
\end{equation}
We refer to this algorithm as \nagsc. Letting $\alpha_i=\sqrt{\mu s}$ in \eqref{eq:original_nag_convergence}, we can see that this algorithm achieves an $O((1-\sqrt{\mu s})^{k}(f(x_{0})-f(x^{*})+\frac{\mu}{2}\|x_{0}-x^{*}\|^{2}))$ convergence rate when $s\leq 1/L$. A major drawback of \nagscsp is that we cannot apply it to non-strongly convex objective functions ($\mu=0$). For non-strongly convex objective functions, \citet{tseng2008accelerated} proposed a simple alternative algorithm to the original NAG \eqref{eq:original_nag}. They set the algorithmic parameters as $\tau_k=\frac{2}{k+1}$ and $\delta_k=\frac{s(k+1)}{2}$ to obtain the following simple algorithm, which we call \nagc:
\begin{equation}
	\label{eq:nag-c}
	\begin{aligned}
		y_{k} & =x_{k}+\frac{2}{k+1}\left(z_{k}-x_{k}\right)\\
		x_{k+1} & =y_{k}-s\nabla f\left(y_{k}\right)\\
		z_{k+1} & =z_{k}-\frac{s(k+1)}{2}\nabla f\left(y_{k}\right).
	\end{aligned}
\end{equation}
When $s\leq 1/L$, this algorithm achieves an $O\left(\|x_{0}-x^{*}\|^{2}/k^2\right)$ convergence rate (see Section~\ref{sec:lyapunov}).

Although there are many variants of NAG, most recent studies on acceleration \citep{diakonikolas2019approximate,shi2019acceleration,siegel2019accelerated,alimisis2020continuous,shi2021understanding,wilson2021lyapunov,kim2022accelerated} focus on these two particular algorithms because of their simplicity. Unfortunately, these two algorithms should be handled separately because \nagscsp \eqref{eq:nag-sc} does not recover \nagcsp \eqref{eq:nag-c} as $\mu\rightarrow0$. 
\begin{center}
	\emph{{\sf Inconsistency I.} \nagscsp does not recover \nagcsp as $\mu\to0$.}
\end{center}
Moreover, \nagscsp has the following drawbacks:
\begin{itemize}
	\item It cannot be applied to non-strongly convex objective functions.
	\item When $\mu$ is very small, the convergence guarantee for \nagscsp is worse than that for \nagcsp in early stages because $(1-\sqrt{\mu s})^k$ converges to $0$ very slowly. 
	\item The convergence rate of \nagscsp depends on both the initial squared distance $\|x_0-x^*\|^2$ and the initial function value accuracy $f(x_0)-f(x^*)$, while the convergence rate of \nagcsp depends only on the squared initial distance $\|x_0-x^*\|^2$.
\end{itemize}
As most of recent works on Nesterov acceleration are based on these two specific algorithms, similar inconsistencies can be found in the literature. We discuss more inconsistencies below.

\subsection{Inconsistencies between convex and strongly convex cases}
\label{sec:inconsistency}

\subsubsection{Continuous-time models.}\label{sec:ctm}

In this subsection, we first informally derive the limiting ODE of the three-sequence scheme \eqref{eq:nag}. To identify a discrete-time sequence $(x_{k})_{k=0}^{\infty}$ with a continuous-time curve $X:[\tst,\infty)\rightarrow\mathbb{R}^{n}$, given the algorithmic stepsize $s$, we introduce a strictly increasing sequence $(\bft_{k})_{k=0}^{\infty}$ (depending on $s$) in $[0,\infty)$ and make the identification $X(\bft_k)=x_k$. We denote the inverse of the sequence $\bft:\{0,1,2,\ldots\}\rightarrow\mathbb{R}$ as ${\bfk}$, that is, ${\bfk}(\bft_k)=k$ for all $k\geq0$. For convenience, we extend the function ${\bfk}$ to a piecewise linear function defined on $[0,\infty)$.

We assume that
\begin{equation}
	\label{eq:tk_condition0}
	\lim_{s\rightarrow0}\bft_{0}=\tst
\end{equation}
and that the timesteps  are asymptotically equivalent to $\sqrt{s}$ as $s\rightarrow0$ in the sense that
\begin{equation}
	\label{eq:tk_condition}
	\lim_{s\rightarrow0}\frac{\bft_{\bfk(t)+1}-t}{\sqrt{s}}=1\text{ for all }t\in\left(\tst,\infty\right).
\end{equation}
Note that the popular choice $\bft_k=t_k:=k\sqrt{s}$ (we will use the notation $t_k$ for this specific sequence throughout the paper) used in \citep{su2014,wibisono2016,shi2021understanding} satisfies these conditions. 

For the iterates of three-sequence scheme \eqref{eq:nag}, we have
\begin{align*}
	\frac{x_{k+1}-x_{k}}{\sqrt{s}} & =\frac{\tau_{k}}{\sqrt{s}}\left(z_{k}-x_{k}\right)-\sqrt{s}\nabla f\left(y_{k}\right)\\
	\frac{z_{k+1}-z_{k}}{\sqrt{s}} & =\frac{\delta_{k}}{\sqrt{s}}\left(\mu y_{k}-\mu z_{k}-\nabla f\left(y_{k}\right)\right).
\end{align*}
We introduce two sufficiently smooth curves $X,Z:[\tst,\infty)\rightarrow\mathbb{R}^{n}$ (possibly depending on $s$ now) such that $X(t)= x_{\bfk(t)}$ and $Z(t)= z_{\bfk(t)}$. Since $\Vert x_{k+1}-y_{k}\Vert=o(\sqrt{s})$ and $\nabla f$ is Lipschitz continuous, we have
\begin{align*}
	\dot{X}(t) & =\lim_{s\rightarrow0}\frac{x_{\bfk(t)+1}-x_{\bfk(t)}}{t_{\bfk(t)+1}-t}=\lim_{s\rightarrow0}\frac{x_{\bfk(t)+1}-x_{\bfk(t)}}{\sqrt{s}}=\lim_{s\rightarrow0}\left\{ \frac{\tau_{\bfk(t)}}{\sqrt{s}}\right\} (Z(t)-X(t))\\
	\dot{Z}(t) & =\lim_{s\rightarrow0}\frac{z_{\bfk(t)+1}-z_{\bfk(t)}}{t_{\bfk(t)+1}-t}=\lim_{s\rightarrow0}\frac{z_{\bfk(t)+1}-z_{\bfk(t)}}{\sqrt{s}}=\lim_{s\rightarrow0}\left\{ \frac{\delta_{\bfk(t)}}{\sqrt{s}}\right\} \left(\mu X(t)-\mu Z(t)-\nabla f(X(t))\right)
\end{align*}
for all $t>0$. Thus, if the limits 
\begin{equation}
	\label{eq:taut_deltat}
	\begin{aligned}
		\tau(t) & =\lim_{s\rightarrow0}\frac{\tau_{\bfk(t)}}{\sqrt{s}}\\
		\delta(t) & =\lim_{s\rightarrow0}\frac{\delta_{\bfk(t)}}{\sqrt{s}}
	\end{aligned}
\end{equation}
exist for all $t\in(\tst,\infty)$, then as $s \to 0$, the iterates generated by the three-sequence scheme \eqref{eq:nag} converge to a solution to the following system of ODEs:
\begin{equation}
	\label{eq:limiting_ode}
	\begin{aligned}
		\dot{X}(t) & =\tau(t)(Z(t)-X(t))\\
		\dot{Z}(t) & =\delta(t)(\mu X(t)-\mu Z(t)-\nabla f(X(t)))
	\end{aligned}
\end{equation}
with the initial conditions $X(\tst)=Z(\tst)=x_0$. We can equivalently write this as the following second-order ODE:
\begin{equation}
	\label{eq:limiting_ode_oneline}
	\ddot{X}+\left(\tau(t)-\frac{\dot{\tau}(t)}{\tau(t)}+\mu\delta(t)\right)\dot{X}+\tau(t)\delta(t)\nabla f(X)=0.
\end{equation}
{Furthermore, when the collinearity condition \eqref{eq:collinear} holds, we have
\begin{equation}
	\label{eq:reciprocal}
	\delta(t)=\lim_{s\rightarrow0}\frac{\delta_{k}}{\sqrt{s}}=\lim_{s\rightarrow0}\frac{1}{\sqrt{s}\left(\mu+(1/s-\mu)\tau_{k}\right)}=\lim_{s\rightarrow0}\frac{\sqrt{s}}{\mu s+(1-\mu s)\tau_{k}}=\frac{1}{\tau(t)}.
\end{equation}}

\paragraph{Limiting ODE of \nagc.}

Recall that \nagcsp \eqref{eq:nag-c} is the three-sequence scheme \eqref{eq:nag} with $\tau_{k}=\frac{2}{k+1}$ and $\delta_{k}=\frac{s(k+1)}{2}$. With the sequence $\bft_k=k\sqrt{s}$, we have
\begin{align*}
	\tau(t) & =\lim_{s\rightarrow0}\frac{\tau_{\bfk(t)}}{\sqrt{s}}=\lim_{s\rightarrow0}\frac{2}{\sqrt{s}\left(t/\sqrt{s}+1\right)}=\frac{2}{t}\\
	\delta(t) & =\lim_{s\rightarrow0}\frac{\delta_{\bfk(t)}}{\sqrt{s}}=\lim_{s\rightarrow0}\frac{\sqrt{s}\left(t/\sqrt{s}+1\right)}{2}=\frac{t}{2}.
\end{align*}
Thus, as $s \to 0$, \nagcsp converges to the following ODE system, which we call \nagcsp system:
\begin{equation}
	\label{eq:c_ode_twoline}
	\begin{aligned}
		\dot{X} & =\frac{2}{t}(Z-X)\\
		\dot{Z} & =-\frac{t}{2}\nabla f(X)
	\end{aligned}
\end{equation}
with $X(0)=Z(0)=x_0$. This system can be written in the following second-order ODE, which we call \nagcsp ODE:
\begin{equation}
	\label{eq:c_ode}
	\ddot{X}+\frac{3}{t}\dot{X}+\nabla f(X)=0
\end{equation}
with $X(0)=x_0$ and $\dot{X}(0)=0$. \citet{su2014} first derived this ODE and showed that the solution to \eqref{eq:c_ode} satisfies an $O(\|x_0-x^*\|^2/t^2)$ convergence rate. 

\paragraph{Limiting ODE of \nagsc.}
Recall that \nagscsp \eqref{eq:nag-sc} is the three-sequence scheme \eqref{eq:nag} with $\tau_{k}=\frac{\sqrt{\mu s}}{1+\sqrt{\mu s}}$ and $\delta_{k}=\sqrt{\frac{s}{\mu}}$. With the sequence $\bft_{k}=-k\frac{\log(1-\sqrt{\mu s})}{\sqrt{\mu}}$,\footnote{Although the sequence $\bft_k=k\sqrt{s}$ leads to the same limiting dynamics, this particular sequence makes a clear connection between the convergence analysis of \nagscsp and that of \nagscsp ODE (see Section~\ref{sec:lyapunov}).} we have
\begin{align*}
	\tau(t) & =\lim_{s\rightarrow0}\frac{\tau_{\bfk(t)}}{\sqrt{s}}=\lim_{s\rightarrow0}\frac{\sqrt{\mu}}{1+\sqrt{\mu s}}=\sqrt{\mu}\\
	\delta(t) & =\lim_{s\rightarrow0}\frac{\delta_{\bfk(t)}}{\sqrt{s}}=\lim_{s\rightarrow0}\frac{1}{\sqrt{\mu}}=\frac{1}{\sqrt{\mu}}.
\end{align*}
Thus, as $s \to 0$, \nagscsp converges to the following ODE system, which we call \nagscsp system: 
\begin{equation}
	\label{eq:sc_ode_twoline}
	\begin{aligned}
		\dot{X} & =\sqrt{\mu}(Z-X)\\
		\dot{Z} & =\frac{1}{\sqrt{\mu}}\left(\mu X-\mu Z-\nabla f(X)\right)
	\end{aligned}
\end{equation}
with $X(0)=Z(0)=x_0$, or equivalently, the following \nagscsp ODE: 
\begin{equation}
	\label{eq:sc_ode}
	\ddot{X}+2\sqrt{\mu}\dot{X}+\nabla f(X)=0\\
\end{equation}
with $X(0)=x_0$ and $\dot{X}(0)=0$. \citet{wilson2021lyapunov} showed that the solution to this ODE satisfies an $O(e^{-\sqrt{\mu}t}(f(x_0)-f(x^*)+\frac{\mu}{2}\|x_0-x^*\|^2))$ convergence rate. Just like in the discrete-time case, \nagcsp ODE \eqref{eq:c_ode} and \nagscsp ODE \eqref{eq:sc_ode} should be handled as separate cases because \nagscsp ODE does not recover \nagcsp ODE as $\mu\to0$.
\begin{center}
	\emph{{\sf Inconsistency II.} \nagscsp ODE does not recover \nagcsp ODE as $\mu\to0$.}
\end{center}
Moreover, \nagscsp ODE has the following drawbacks:
\begin{itemize}
	\item The solution to \nagscsp ODE with $\mu=0$ may not converge to the minimizer of $f$: For the objective function $f(x)=\frac{1}{2}x^2$ on $\mathbb{R}$, the solution to \nagscsp ODE with $x_0=1$ is $X(t)=\cos(t)$, which does not converge to the minimizer $x^{*}=0$.
	\item When $\mu$ is very small, the convergence guarantee for \nagscsp ODE is worse than that for \nagcsp ODE in early stages because $e^{-\sqrt{\mu}t}$ converges to $0$ very slowly. 
	\item The convergence rate of \nagscsp ODE depends on both the initial squared distance $\|x_0-x^*\|^2$ and the initial function value accuracy $f(x_0)-f(x^*)$, while the convergence rate of \nagcsp ODE depends only on the squared initial distance $\|x_0-x^*\|^2$.
\end{itemize}

\subsubsection{Bregman Lagrangians}

To systematically study the acceleration phenomenon of momentum methods, \citet{wibisono2016} introduced the following \emph{first} Bregman Lagrangian:
\begin{equation}
	\label{eq:first_lagrangian}
	\mathcal{L}_{\text{1st}}\left(X,\dot{X},t\right)=e^{\alpha+\gamma}\left(D_{h}\left(X+e^{-\alpha}\dot{X},X\right)-e^{\beta}f(X)\right),
\end{equation}
where $\alpha,\beta,\gamma:[\tst,\infty)\rightarrow\mathbb{R}$ are continuously differentiable functions, $h$ is a continuously differentiable  strictly convex function, and $D_h$ is the Bregman divergence (see Section~\ref{sec:convexity} for its definition). In order to obtain accelerated convergence rates,  the following \emph{ideal scaling conditions} are introduced:
\begin{subequations}
	\label{eq:ideal_scaling}
	\begin{align}
		\label{eq:ideal_scaling_a}
		\dot{\gamma} & =e^{\alpha}\\
		\label{eq:ideal_scaling_b}
		\dot{\beta} & \leq e^{\alpha}.
	\end{align}
\end{subequations}
Under the ideal scaling condition \eqref{eq:ideal_scaling_a}, the Euler--Lagrange equation
\begin{equation}
	\label{eq:e-l-eq}
	\frac{d}{dt}\left\{ \frac{\partial\mathcal{L}}{\partial \dot{X}}\left(X,\dot{X},t\right)\right\} =\frac{\partial\mathcal{L}}{\partial X}\left(X,\dot{X},t\right)
\end{equation}
for the first Bregman Lagrangian \eqref{eq:first_lagrangian} reduces to the following system of first-order equations:
\begin{subequations}
	\label{eq:first_family}
	\begin{align}
		\dot{X} & =e^{\alpha}(Z-X)\\
		\frac{d}{dt}\nabla h(Z) & =-e^{\alpha+\beta}\nabla f(X).
	\end{align}
\end{subequations}
When $f$ is convex, any solution to the system of ODEs \eqref{eq:first_family} reduces the objective function value accuracy at an $O(e^{-\beta(t)})$ convergence rate (see Section~\ref{sec:lyapunov}). In particular, setting $\alpha(t)=\log\frac{2}{t}$ and $\beta(t)=\log\frac{t^{2}}{4}$, we recover \nagcsp system \eqref{eq:c_ode_twoline} and its convergence rate.

Although the first Bregman Lagrangian \eqref{eq:first_lagrangian} generates a large family of momentum dynamics, it does not include \nagscsp system \eqref{eq:sc_ode_twoline}. To handle strongly convex cases, \citet{wilson2021lyapunov} introduced the \emph{second} Bregman Lagrangian, defined as
\begin{equation}
	\label{eq:second_lagrangian}
	\mathcal{L}_{\text{2nd}}\left(X,\dot{X},t\right)=e^{\alpha+\beta+\gamma}\left(\mu D_{h}\left(X+e^{-\alpha}\dot{X},X\right)-f(X)\right).
\end{equation}
Under the ideal scaling condition \eqref{eq:ideal_scaling_a}, the Euler--Lagrange equation \eqref{eq:e-l-eq} for the second Bregman Lagrangian \eqref{eq:second_lagrangian} reduces to the following system of first-order equations:
\begin{subequations}
	\label{eq:second_family}
	\begin{align}
		\dot{X} & =e^{\alpha}(Z-X)\\
		\frac{d}{dt}\nabla h(Z) & =\dot{\beta}\left(\nabla h(X)-\nabla h(Z)\right)-\frac{e^{\alpha}}{\mu}\nabla f(X).
	\end{align}
\end{subequations}
When $f$ is $\mu$-uniformly convex with respect to $h$ (see Section~\ref{sec:convexity}), any solution to the system of ODEs \eqref{eq:second_family} satisfies an $O(e^{-\beta(t)})$ convergence rate (see Section~\ref{sec:lyapunov}). In particular, letting $\alpha(t)=\log\sqrt{\mu}$ and $\beta(t)=\sqrt{\mu}t$, we recover \nagscsp system \eqref{eq:sc_ode_twoline} and its convergence rate. Here, we observe an inconsistency between the two Bregman Lagrangians.

\begin{center}
	\emph{{\sf Inconsistency III.} The second Bregman Lagrangian does not recover the first Bregman Lagrangian as $\mu\to0$.}
\end{center}

\subsection{Contributions}
\label{sec:motivation}

\begin{figure}[!ht] 
	\begin{center}
		\begin{tikzpicture}[node distance=2.3cm]
			\node (ode1) [box, draw=black] {\nagcsp ODE \citep{su2014}};
			\node (lag1) [box, draw=black, above of=ode1] {First Bregman Lagrangian \citep{wibisono2016}};
			\node (lagu) [box, draw=black, above of=ode1, xshift=+5.5cm] {\textbf{Unified Bregman Lagrangian (Section~\ref{sec:unified_lagrangian})}};
			\node (lag2) [box, draw=black, above of=ode1, xshift=+11cm] {Second Bregman Lagrangian \citep{wilson2021lyapunov}};
			\node (odeu) [box, draw=black, below of=lagu] {\textbf{Unified NAG ODE (Section~\ref{sec:unified_ode})}};
			\node (ode2) [box, draw=black, below of=lag2] {\nagscsp ODE \citep{wilson2021lyapunov}};
			\node (nag1) [box, draw=black, below of=ode1] {\nagcsp \citep{tseng2008accelerated}};
			\node (nagu) [box, draw=black, below of=odeu] {\textbf{Unified NAG (Section~\ref{sec:unified_nag})}};
			\node (nag2) [box, draw=black, below of=ode2] {\nagscsp \citep{nesterov2018lectures}};
			\node (bigbox) [dashed, inner sep=0.3cm, box, fit=(lag1)(nag2)] {};
			\node (higher) [box, draw=black, below of=bigbox, xshift=-2.75cm, yshift=-2.6cm] {\textbf{Higher-order optimization (Section~\ref{sec:higher-order})}};
			\node (gradient) [box, draw=black, below of=bigbox, xshift=+2.75cm, yshift=-2.6cm] {\textbf{Gradient norm minimization (Section~\ref{sec:minimizing_gradient})}};
			
			\draw [arrow] (bigbox.south -| higher.north) -- (higher.90);
			\draw [arrow] (bigbox.south -| gradient.north) -- (gradient.90);
			
			\draw [arrow] (lag1) --node [left] {special case} (ode1);
			\draw [arrow] (lag2) --node [right] {special case} (ode2);
			\draw [arrow] (ode1.250) --node [left] {discretize} (nag1.110);
			\draw [arrow] (nag1.70) --node [right] {limit} (ode1.290);
			\draw [arrow] (ode2.250) --node [left] {discritize} (nag2.110);
			\draw [arrow] (nag2.70) --node [right] {limit} (ode2.290);
			\draw [arrow] (lag1) --node [above] {unify} (lagu);
			\draw [arrow] (lag2) --node [above] {unify} (lagu);
			\draw [arrow] (lagu) -- (odeu);
			\draw [arrow] (odeu.250) -- (nagu.110);
			\draw [arrow] (nagu.70) -- (odeu.290);
			\draw [arrow] [dashed] (odeu) --node [above] {recover} (ode1);
			\draw [arrow] [dashed] (odeu) --node [below] {$\mu=0$} (ode1);
			\draw [arrow] [dashed] (nagu) --node [above] {recover} (nag1);
			\draw [arrow] [dashed] (nagu) --node [below] {$\mu=0$} (nag1);
			\draw [arrow] [dashed] (odeu) --node [above] {recover} (ode2);
			\draw [arrow] [dashed] (odeu) --node [below] {$t\rightarrow\infty$} (ode2);
			\draw [arrow] [dashed] (nagu) --node [above] {recover} (nag2);
			\draw [arrow] [dashed] (nagu) --node [below] {$k\rightarrow\infty$} (nag2);
		\end{tikzpicture}
	\end{center}
	\caption{An illustration of our framework and contributions.}
	\label{fig:chart}
\end{figure}
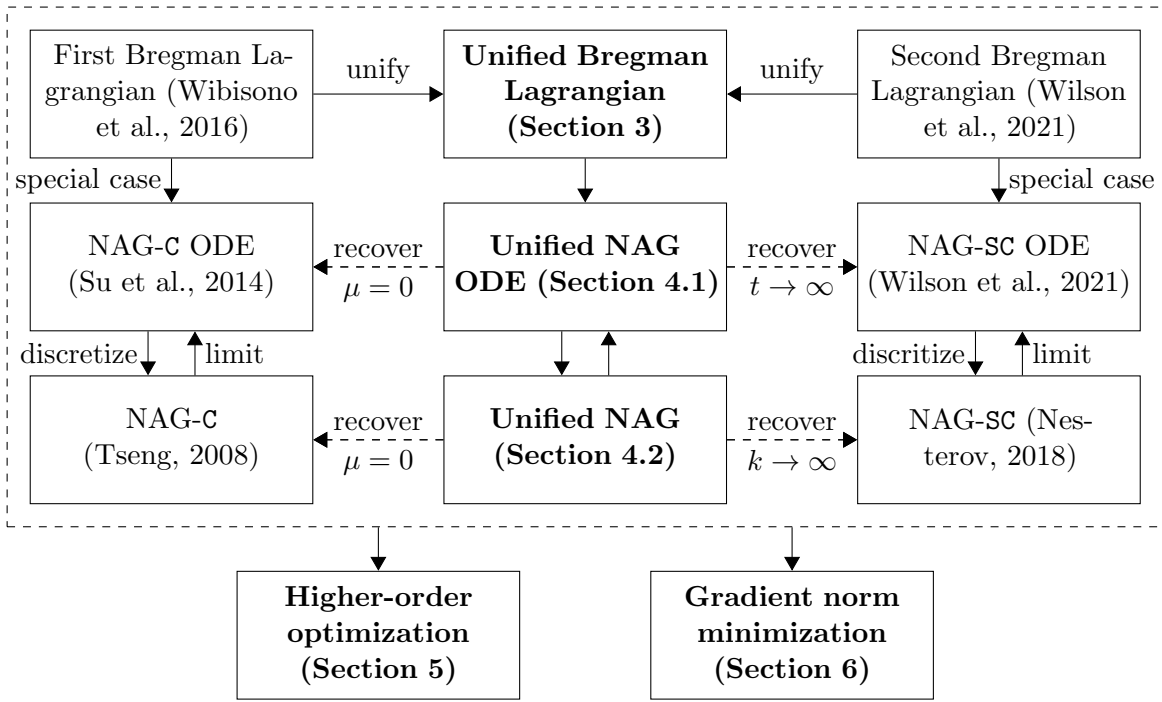

In this paper, we propose a novel unified framework for Lagrangians, ODE models and algorithms 
to address the inconsistencies between the convex case and the strongly convex case mentioned above.
The proposed framework  seamlessly bridges the gap between the two cases as illustrated in Figure~\ref{fig:chart}. The main contributions of this work can be summarized as follows:
\begin{itemize}
	\item We propose the \emph{unified} Bregman Lagrangian (Section~\ref{sec:unified_lagrangian}). Unlike the second Bregman Lagrangian, the unified Bregman Lagrangian recovers the first Bregman Lagrangian  when $\mu=0$. As the Euler--Lagrange equation for the unified Bregman Lagrangian, we obtain a family of continuous-time dynamics (Proposition~\ref{prop:e-l-eq}). Using a Lyapunov function, we analyze the convergence rate for these flows (Theorem~\ref{thm:mainthm_lagrangian}).
	\item We derive the \emph{unified} NAG ODE \eqref{eq:u_ode_oneline} as a special case of the unified Bregman Lagrangian flows (Section~\ref{sec:unified_ode}). Unlike \nagscsp ODE \eqref{eq:sc_ode}, for non-strongly convex objective functions ($\mu=0$), the unified NAG ODE and its convergence rate (Theorem~\ref{thm:mainthm_discrete}) recover \nagcsp ODE \eqref{eq:c_ode} and its convergence rate. Furthermore, for any $\mu > 0$, the unified NAG ODE and its convergence rate (Corollary~\ref{cor:maincor_continuous}) recover \nagscsp ODE \eqref{eq:sc_ode} and its convergence rate as $t\rightarrow\infty$. 
	\item We devise the \emph{unified} NAG family \eqref{eq:nag-u}, a family of momentum algorithms that converge to the unified NAG ODE as $s \to 0$ (Section~\ref{sec:unified_nag}). As a special case, we have the \emph{unified} NAG \eqref{eq:unified_nag_specific}, a simple algorithm which unifies \nagcsp \eqref{eq:nag-c} and \nagscsp \eqref{eq:nag-sc}. Moreover, using an adaptive timestep in the unified NAG family, we constructively recover the original NAG \eqref{eq:original_nag} with $\gamma_0>\mu$ and its convergence rate \eqref{eq:original_nag_convergence}.
	\item We extend the unified NAG ODE and the unified NAG family to the higher-order non-Euclidean setting (mirror descent setup) (Section~\ref{sec:higher-order}). Our novel dynamics and algorithms can be viewed as continuous extensions of the \emph{accelerated tensor method (convex case)} and its limiting ODE in \citep{wibisono2016} to the strongly convex setting.
\end{itemize}
We also made the following contributions that are not closely related to our major goal but may deserve independent attention:
\begin{itemize}
	\item We compute the general limiting ODEs of the three-sequence scheme \eqref{eq:nag}, the two-sequence scheme \eqref{eq:two_sequence_scheme}, and the fixed-step first-order scheme \eqref{eq:discrete_FSFO}. In particular, we introduce a novel tool, called the \emph{differential kernel} $H(t,\tau)$, to derive the limiting ODE of the fixed-step first-order scheme. We show that an anti-transpose relationship \eqref{eq:symmetry_b_c} between OGM and OGM-G can be naturally shifted to a continuous-time setting by this tool.
	\item We propose the \emph{unified NAG-G ODE}, an ODE model for minimizing the gradient norm of strongly convex objective functions (Section~\ref{sec:minimizing_gradient}). Surprisingly, the differential kernels corresponding to the unified NAG ODE and the unified NAG-G ODE have an anti-transpose relationship, just like it does between OGM ODE and OGM-G ODE.
\end{itemize}

\renewcommand{\arraystretch}{1.5}
\begin{table*}[ht]
	\begin{center}
		\scalebox{0.85}{\begin{tabular}{|c|c|}
			\hline
			\rowcolor{lightgray}
			Dynamics & Convergence rate \\ \hline
			{Unified NAG ODE} & $f\big(X\big(t\big)\big)-f\big(x^{*}\big)\leq O\big(\min\big\{1/t^{2},e^{-\sqrt{\mu}t}\big\}\big\Vert x_{0}-x^{*}\big\Vert^{2}\big)$ \\ \hline
			Unified accelerated tensor flow & $f\big(X\big(t\big)\big)-f\big(x^{*}\big)\leq O\big(\min\big\{1/t^{p},e^{-pC^{1/p}\mu^{1/p}t}\big\} D_{h}\big(x^{*},x_{0}\big)\big)$ \\ \hline
			Unified NAG-G ODE & $\big\Vert\nabla f(X(T))\big\Vert^{2}\leq O\big(\min\big\{1/T^{2},e^{-\sqrt{\mu}T}\big\}\big(f\big(x_0\big)-f\big(x^{*}\big)\big)\big)$ \\ \hline
			\rowcolor{lightgray}
			Algorithm & Convergence rate\\ \hline
			{Unified NAG} & $f\big(x_{k}\big)-f\big(x^{*}\big)\leq O\big(\min\big\{1/k^{2},\big(1-\sqrt{\mu s}\big)^{k}\big\}\big\Vert x_{0}-x^{*}\big\Vert^{2}\big)$ \\ \hline
			Unified accelerated tensor method & $f\big(x_{k}\big)-f\big(x^{*}\big)\leq O\big(\min\big\{1/k^{p},\big(1+C^{1/p}p\mu^{1/p}s^{1/p}\big)^{-k}\big\}D_{h}\big(x^{*},x_{0}\big)\big)\big)$ \\ \hline
		\end{tabular}}
	\end{center}
	\caption{Convergence rates of the momentum dynamics and algorithms proposed in this paper.}
	\label{table:convergence}
\end{table*}
\renewcommand{\arraystretch}{1}

We summarize the convergence rates for our continuous-time dynamics and discrete-time algorithms in Table~\ref{table:convergence}. In addition to theoretical and algorithmic perspectives, 
we discuss the need for unified acceleration methods from a practical perspective.

\paragraph{Practical perspective.}
Many optimization problems in machine learning can be formulated as
\begin{equation}
	\label{eq:ml_problem}
	\min_{x\in\mathbb{R}^{n}}\;f(x)=\frac{1}{m}\left(\sum_{i=1}^{m}f_{i}(x)+\lambda R(x)\right),
\end{equation}
where $f_i$ is the loss function corresponding to the $i$-th sample, $\lambda>0$ is the regularization parameter, and $R(x)$ is the regularization term \citep[Equation~1.1]{bubeck2015convex}. Consider the problem \eqref{eq:ml_problem} where the functions $f_i$ are convex and $L$-smooth, and $R(x)=\Vert x\Vert ^{2}$. Then, $f$ is $\mu$-strongly convex and $L$-smooth, where $\mu={2\lambda}/{m}$. As the sample size $m$ grows or the regularization parameter $\lambda$ decreases, the strong convexity parameter $\mu$ decreases. Thus, improving the convergence rate for ill-conditioned strongly convex objective functions (where $\mu$ is small) is quite significant, as emphasized in \citep[Section~3.6]{bubeck2015convex}.

As mentioned above, the convergence guarantee of \nagscsp \eqref{eq:nag-sc} is no better than that of \nagcsp \eqref{eq:nag-c} when $\mu$ is small. In our numerical experiments (see Section~\ref{sec:experiments}), it is observed that the performance of \nagscsp is worse than that of \nagcsp when $\mu$ is very small. Thus, it is desirable to design a strongly convex optimization algorithm whose convergence guarantee is not worse than that of \nagcsp even when $\mu$ is very small. In the experiments, we observe that for a logistic regression problem, when $\mu$ is small, our algorithm is comparable to \nagc, while \nagscsp underperforms \nagc.

{
\paragraph{Existing unified methods and dynamics.}
To clarify what is our novel contribution and what is not, we review existing algorithms and dynamics that can handle the non-strongly convex case and the strongly convex case in a unified way. The original NAG \eqref{eq:original_nag} is an accelerated algorithm that can handle both convex objective functions and strongly convex objective functions. In Section~\ref{sec:variable}, we show that the original NAG can be constructively recovered by our unified Lagrangian formulation. \citet{luo2021differential} designed the following ODE model for the original NAG, which we call the \emph{original NAG system}:
\begin{equation}
	\label{eq:original_nag_flow}
	\begin{aligned}
		\dot{\gamma} & =\mu-\gamma\\
		\dot{X} & =Z-X\\
		\dot{Z} & =\frac{1}{\gamma}(\mu X-\mu Z-\nabla f(X))
	\end{aligned}
\end{equation}
with $X(0)=Z(0)=x_0$ and $\gamma(0)=\gamma_0>0$. \citet[Section~6.2]{luo2021differential} showed that the original NAG can be viewed as a discretization scheme with the timestep $\alpha_i$, which is inductively defined in \eqref{eq:alpha_k_update}. 

Using time rescaling technique, \citet{luo2021differential} also proposed the following system of ODEs (although most of their results directly deal with Equation~\ref{eq:original_nag_flow}):
\begin{equation}
	\label{eq:original_nag_flow_scaled}
	\begin{aligned}
		\dot{X}(t) & =a(t)(Z(t)-X(t))\\
		b(t)\dot{Z}(t) & =a(t)(\mu X(t)-\mu Z(t)-\nabla f(X(t))),
	\end{aligned}
\end{equation}
where $a:[0,\infty)\to[0,\infty)$ is an arbitrary function and $$b(t)=\gamma\left(\int_0^t a(s)\,ds\right).$$ 
This ODE system is closely related to the \emph{unified Bregman Lagrangian flow} \eqref{eq:unified_family} and the \emph{unified NAG system} \eqref{eq:u_ode} proposed in this paper. In Appendix~\ref{app:existing_unified1}, we show that the rescaled original NAG flow \eqref{eq:original_nag_flow_scaled} can be expressed as the unified Bregman Lagrangian flow \eqref{eq:unified_family}. Conversely, the unified Bregman Lagrangian flow can be expressed as the rescaled original NAG flow if the ideal scaling condition \eqref{eq:ideal_scaling_b} holds with equality and the distance-generating function $h$ is Euclidean ($h(x)=\frac{1}{2}\|x\|^2$). Therefore, our unified Bregman Lagrangian generates a strictly larger family compared to \eqref{eq:original_nag_flow_scaled}. To emphasize, only our family can deal with the non-Euclidean setup (mirror descent setup). In addition, the derivation of our unified family \eqref{eq:unified_family} is more constructive because it comes from a Lagrangian formulation, whereas \citet{luo2021differential} designed the family \eqref{eq:original_nag_flow_scaled} through heuristic speculation. 

\subsection{Related work}
\citet{nesterov1983} first proposed the original NAG \eqref{eq:original_nag} with $\mu=0$. The original NAG with $\mu>0$ was first analyzed using the estimate sequence technique \citep{nesterov2018lectures}. \citet{tseng2008accelerated} proposed \nagcsp \eqref{eq:nag-c} and its generalization to composite optimization problems. \citet{su2014} derived \nagcsp ODE \eqref{eq:c_ode} by taking the limit $s\rightarrow0$ in NAG-C. This ODE has further been generalized and investigated in \citep{krichene2015accelerated,attouch2018fast}. \citet{wibisono2016} proposed the first Bregman Lagrangian \eqref{eq:first_lagrangian} that systematically generates a family of ODEs \eqref{eq:first_family} including \nagcsp ODE and its higher-order extensions. \citet{wilson2021lyapunov} extended this framework to the strongly convex case. They proposed the second Bregman Lagrangian \eqref{eq:second_lagrangian}, which generates a family of continuous-time flows \eqref{eq:second_family} including \nagscsp ODE \eqref{eq:sc_ode_twoline}, and strengthened the connection between continuous-time dynamics and discrete-time algorithms via Lyapunov function arguments. However, as mentioned in Section~\ref{sec:inconsistency}, their work is not consistent with \citep{wibisono2016} because the second Bregman Lagrangian does not recover the first Bregman Lagrangian as $\mu \to 0$. Based on Lagrangian formulations, \citet{betancourt2018symplectic} studied a symplectic integrator to obtain discrete-time algorithms from continuous-time dynamics. \citet{shi2019acceleration,shi2021understanding} derived high-resolution ODEs for \nagcsp and \nagsc, and then obtained algorithms with accelerated convergence rates by applying the symplectic Euler method to the high-resolution ODEs. \citet{luo2021differential} understood acceleration using the $\mathcal{A}$-stability theory and designed an ODE model for the original NAG method.
\citet{zhang2021revisiting} obtained an accelerated algorithm by applying the explicit Euler method to a variant of high-resolution ODEs. \citet{diakonikolas2019approximate} proposed the approximate duality gap technique to construct and analyze accelerated algorithms. Using conservation laws in dilated coordinate systems, \citet{suh2022continuous} recovered \nagcsp ODE and \nagscsp ODE and showed that a semi-second-order symplectic Euler discretization in the dilated coordinate yields accelerated methods.

\section{Preliminaries}
In this section, we review the basic notions that we will use throughout the paper. While Sections~\ref{sec:convexity} and \ref{sec:lyapunov} review the standard concepts in the literature, Sections~\ref{sec:hyperbolic} and \ref{sec:limiting} contain novel ideas and results.

\subsection{Convex analysis}
\label{sec:convexity}

\paragraph{Convexity and smoothness.}
Let $f:\mathbb{R}^n\rightarrow\mathbb{R}$ be a $C^{\infty}$ function. Then for $\mu\geq0$, the function $f$ is called \emph{$\mu$-strongly convex} if the inequality
\[
f(y)\geq f(x)+\left\langle \nabla f(x),y-x\right\rangle +\frac{\mu}{2}\left\Vert y-x\right\Vert ^{2}
\]
holds for all $x,y\in\mathbb{R}^n$. In particular, the function $f$ is called \emph{convex} if it is strongly convex with the strong convexity parameter $\mu=0$. For $L>0$, the function $f$ is called \emph{$L$-smooth} if its gradient is $L$-Lipschitz continuous, that is, the inequality
\[
\left\Vert \nabla f(x)-\nabla f(y)\right\Vert \leq L\left\Vert x-y\right\Vert 
\]
holds for all $x,y\in\mathbb{R}^n$. It is known that when $f$ is $L$-smooth, the inequality 
\[
f(y)\leq f(x)+\left\langle \nabla f(x),y-x\right\rangle +\frac{L}{2}\left\Vert y-x\right\Vert ^{2}
\]
holds for all $x,y\in\mathbb{R}^n$. For most of the remaining sections of this paper (Sections~\ref{sec:unified_sections} and \ref{sec:minimizing_gradient}), we make the following assumptions, which we call the \emph{standard smooth strongly convex setting}:
\begin{itemize}
	\item The objective function $f$ is $(1/s)$-smooth, where $s>0$ is the algorithmic stepsize. 
	\item The objective function $f$ is $\mu$-strongly (possibly with $\mu=0$) convex.
\end{itemize}

\paragraph{Higher-order convexity and smoothness.}
The notions of convexity and smoothness can be generalized to the higher-order setting. The function $f$ is called \emph{$\mu$-uniformly convex of order $p\geq2$} if the inequality
\begin{equation}
	\label{eq:higher_order_uniform_convexity}
	f(y)\geq f(x)+\left\langle \nabla f(x),y-x\right\rangle +\frac{\mu}{p}\left\Vert y-x\right\Vert ^{p}
\end{equation}
holds for all $x,y\in\mathbb{R}^n$. The function $f$ is called \emph{$L$-smooth of order $p-1$} if the inequality
\begin{equation}
	\label{eq:higher_order_smoothness}
	\left\Vert \nabla^{p-1}f(y)-\nabla^{p-1}f(x)\right\Vert \leq L\left\Vert y-x\right\Vert 
\end{equation}
holds for all $x,y\in\mathbb{R}^n$. Note that these definitions recover the standard notions of convexity and smoothness when $p=2$.

\paragraph{Bregman divergences.}
In the optimization literature, a common way to consider a non-Euclidean setting is by using the Bregman divergence, instead of the Euclidean distance. For a continuously differentiable function $h:\mathbb{R}^n\rightarrow\mathbb{R}$ which is convex and essentially smooth ($\Vert \nabla h(x)\Vert \rightarrow \infty$ as $\Vert x \Vert \rightarrow\infty$), the \emph{Bregman divergence} $D_h:\mathbb{R}^n\times\mathbb{R}^n\rightarrow[0,\infty)$ of $h$ is defined as
\begin{equation}
	\label{eq:bregman_divergence}
	D_{h}(y,x)=h(y)-h(x)-\left\langle \nabla h(x),y-x\right\rangle.
\end{equation}
Note that when $h(x)=\frac{1}{2}\Vert x\Vert ^{2}$, the Bregman divergence of $h$ is the squared Euclidean distance $\frac{1}{2}\Vert y-x\Vert ^{2}$. For all $x,y,z\in\mathbb{R}^n$, the \emph{three-point identity} \citep[see][Proposition~5]{wilson2021lyapunov}
\begin{equation}
	\label{eq:three-point}
	D_{h}(x,y)-D_{h}(x,z)=-\left\langle \nabla h(y)-\nabla h(z),x-y\right\rangle -D_{h}(y,z)
\end{equation}
holds. For $\mu\geq0$, the function $f$ is called \emph{$\mu$-uniformly convex with respect to $h$} if the inequality
\begin{equation}
	\label{eq:uniform_convexity}
	D_f(x,y)\geq\mu D_h(x,y)
\end{equation}
holds for all $x,y\in\mathbb{R}^n$. Note that this condition is equivalent to the $\mu$-storng convexity of $f$ when $h(x)=\frac{1}{2}\left\Vert x\right\Vert ^{2}$.

\subsection{Lyapunov arguments for convergence analyses}
\label{sec:lyapunov}

A popular method for proving the convergence rates of momentum dynamics and algorithms is constructing an energy function decreasing over time, called the Lyapunov function \citep{lyapunov1992general}. The particular analyses presented in this section handle discrete-time algorithms and the corresponding continuous-time dynamics using a single Lyapunov function, as in \citep{krichene2015accelerated}. To prove the convergence rates of the given algorithm and associated dynamics, we take the following steps:
\begin{enumerate}
	\item Define a time-dependent Lyapunov function $V:\mathbb{R}^n\times\mathbb{R}^n\times[0,\infty)\rightarrow [0, \infty)$.
	\item Show that the continuous-time energy functional $\mathcal{E}(t)=V(X(t),Z(t),t)$ is monotonically decreasing along the solution trajectory $(X,Z):[0,\infty)\rightarrow\mathbb{R}^n\times\mathbb{R}^n$ of the ODE system.
	\item Show that the discrete-time energy functional $\mathcal{E}_k=V(x_k,z_k,\bft_k)$ is monotonically decreasing along the iterates $(x_k,z_k):\{0,1,2,\ldots\}\rightarrow\mathbb{R}^n\times\mathbb{R}^n$ of the algorithm.
\end{enumerate}
The remainder of this subsection shows how we can apply this strategy to known algorithms. We assume the standard smooth (strongly) convex setting (see Section~\ref{sec:convexity}).

\paragraph{\nagcsp and \nagcsp ODE.}
We define a time-dependent Lyapunov function as 
\begin{equation}
	\label{eq:c_ode_lyapunov}
	V(X,Z,t):=\frac{1}{2}\left\Vert Z-x^{*}\right\Vert ^{2}+\frac{t^{2}}{4}\left(f(X)-f\left(x^{*}\right)\right).
\end{equation}
Then, the continuous-time energy functional
\[
\mathcal{E}(t)=V(X(t),Z(t),t) = \frac{1}{2}\left\Vert Z(t)-x^{*}\right\Vert ^{2}+\frac{t^{2}}{4}\left(f(X(t))-f\left(x^{*}\right)\right)
\]
is monotonically decreasing along the solution trajectory of \nagcsp ODE \eqref{eq:c_ode_twoline} \citep[see][]{su2016differential}. Writing $\mathcal{E}(t)\leq \mathcal{E}(0)$ explicitly, we obtain an $O(1/t^2)$ convergence rate as
\[
f(X(t))-f(x^*) \leq \frac{4}{t^2}\mathcal{E}(t) \leq \frac{4}{t^2}\mathcal{E}(0) = \frac{2}{t^2}\left\Vert x_0 - x^* \right\Vert ^2.
\]
For the iterates of \nagcsp \eqref{eq:nag-c}, the discrete-time energy function 
\begin{equation}
	\label{eq:energy_nagc}
	\mathcal{E}_k=V(x_k,z_k,\bft_k)= \frac{1}{2}\left\Vert z_k-x^{*}\right\Vert ^{2}+\frac{sk^2}{4}\left(f(x_k)-f\left(x^{*}\right)\right),
\end{equation}
where $\bft_k=k\sqrt{s}$, is monotonically decreasing \citep[see][Chapter~12]{ryu2022large}. Hence, we obtain an $O(1/k^2)$ convergence rate.

\paragraph{\nagscsp and \nagscsp ODE.}
We define a time-dependent Lyapunov function as
\begin{equation}
	\label{eq:sc_ode_lyapunov}
	V(X,Z,t):=e^{\sqrt{\mu}t}\left(\frac{\mu}{2}\left\Vert Z-x^{*}\right\Vert ^{2}+f(X)-f\left(x^{*}\right)\right).
\end{equation}
Then we can show that \nagscsp ODE \eqref{eq:sc_ode_twoline} achieves an $O(e^{-\sqrt{\mu}t})$ convergence rate by showing that the energy functional
\[
\mathcal{E}(t)=V(X(t),Z(t),t)=e^{\sqrt{\mu}t}\left(\frac{\mu}{2}\left\Vert Z(t)-x^{*}\right\Vert ^{2}+f(X(t))-f\left(x^{*}\right)\right)
\]
is monotonically decreasing along the solution trajectory of \nagscsp ODE \citep[see][]{wilson2021lyapunov}. Similarly, we can show that \nagscsp \eqref{eq:nag-sc} achieves an $O((1-\sqrt{\mu s})^k)$ convergence rate by showing that the energy functional
\[
\mathcal{E}_{k}=V(x_{k},z_{k},\bft_{k})=\left(1-\sqrt{\mu s}\right)^{-k}\left(\frac{\mu}{2}\left\Vert z_{k}-x^{*}\right\Vert ^{2}+f(x_{k})-f\left(x^{*}\right)\right),
\]
where $\bft_{k}=-k\frac{\log(1-\sqrt{\mu s})}{\sqrt{\mu}}$, is monotonically decreasing along the iterates of \nagscsp \citep[see][Section~4.5]{daspremont2021acceleration}.

\paragraph{Bregman Lagrangians.}
We can show that the first Bregman Lagrangian flow \eqref{eq:first_family} and the second Bregman Lagrangian flow \eqref{eq:second_family} achieve an $O(e^{-\beta(t)})$ convergence rate by showing that the energy functional $\mathcal{E}(t)=V(X(t),Z(t),t)$ is monotonically decreasing, where the Lyapunov function $V$ is defined as
\begin{equation}
	\label{eq:first_lagrangian_lyapunov}
	V_{\mathrm{1st}}(X,Z,t):=D_{h}\left(x^{*},Z\right)+e^{\beta(t)}\left(f(X)-f\left(x^{*}\right)\right)
\end{equation}
for the first Bregman Lagrangian flow and
\begin{equation}
	\label{eq:second_lagrangian_lyapunov}
	V_{\mathrm{2nd}}(X,Z,t) :=e^{\beta(t)}\left(\mu D_{h}\left(x^{*},Z\right)+f(X)-f\left(x^{*}\right)\right)
\end{equation}
for the second Bregman Lagrangian flow. See \citep{wibisono2016,wilson2021lyapunov} for the proofs.

\subsection{Hyperbolic functions and their higher-order generalization}
\label{sec:hyperbolic}

\paragraph{Hyperbolic functions.}
We first review the definitions and properties of hyperbolic functions. The $\sinh$, $\cosh$, $\tanh$, $\coth$, $\sech$, and $\csch$ functions are defined as
\begin{equation}
	\label{eq:hyperbolic}
	\begin{aligned}
		\sinh x & =\frac{e^{x}-e^{-x}}{2}, & \sinh x & \sim x\text{ as }x\rightarrow0, & \sinh x & \sim\frac{e^{x}}{2}\text{ as }x\rightarrow\infty\\
		\cosh x & =\frac{e^{x}+e^{-x}}{2}, & \cosh x & \sim1\text{ as }x\rightarrow0, & \cosh x & \sim\frac{e^{x}}{2}\text{ as }x\rightarrow\infty\\
		\tanh x & =\frac{\sinh x}{\cosh x}, & \tanh x & \sim x\text{ as }x\rightarrow0, & \tanh x & \sim1\text{ as }x\rightarrow\infty\\
		\coth x & =\frac{\cosh x}{\sinh x}, & \coth x & \sim\frac{1}{x}\text{ as }x\rightarrow0, & \coth x & \sim1\text{ as }x\rightarrow\infty\\
		\sech x & =\frac{1}{\cosh x}, & \sech x & \sim1\text{ as }x\rightarrow0, & \sech x & \sim2e^{-x}\text{ as }x\rightarrow\infty\\
		\csch x & =\frac{1}{\sinh x}, & \csch x & \sim\frac{1}{x}\text{ as }x\rightarrow0, & \csch x & \sim2e^{-x}\text{ as }x\rightarrow\infty.
	\end{aligned}
\end{equation}
Furthermore, the $\sinhc$, $\tanhc$, $\cothc$, and $\cschc$ functions are defined as follows \citep[see][]{ten2012extension}:
\begin{equation}
	\label{eq:sinhc}
	\begin{aligned}
		\sinhc x & :=\begin{cases}
			\frac{\sinh x}{x}, & \text{if }x\neq0\\
			1, & \text{if }x=0
		\end{cases} & \sinhc x & \sim1\text{ as }x\rightarrow0, & \sinhc x & \sim\frac{e^{x}}{2x}\text{ as }x\rightarrow\infty\\
		\tanhc x & :=\frac{\sinhc x}{\cosh x} & \tanhc x & \sim1\text{ as }x\rightarrow0, & \tanhc x & \sim\frac{1}{x}\text{ as }x\rightarrow\infty\\
		\cothc x & :=\frac{1}{\tanhc x}, & \cothc x & \sim1\text{ as }x\rightarrow0, & \cothc x & \sim x\text{ as }x\rightarrow\infty\\
		\cschc x & :=\frac{1}{\sinhc x}, & \cschc x & \sim1\text{ as }x\rightarrow0, & \cschc x & \sim2xe^{-x}\text{ as }x\rightarrow\infty.
	\end{aligned}
\end{equation}
The graphs of these functions are shown in Figure~\ref{fig:hyp}.

\begin{figure}[ht]
	\centering
	\subfigure[$\sinh$, $\cosh$, $\tanh$]{\includegraphics[width=0.32\textwidth]{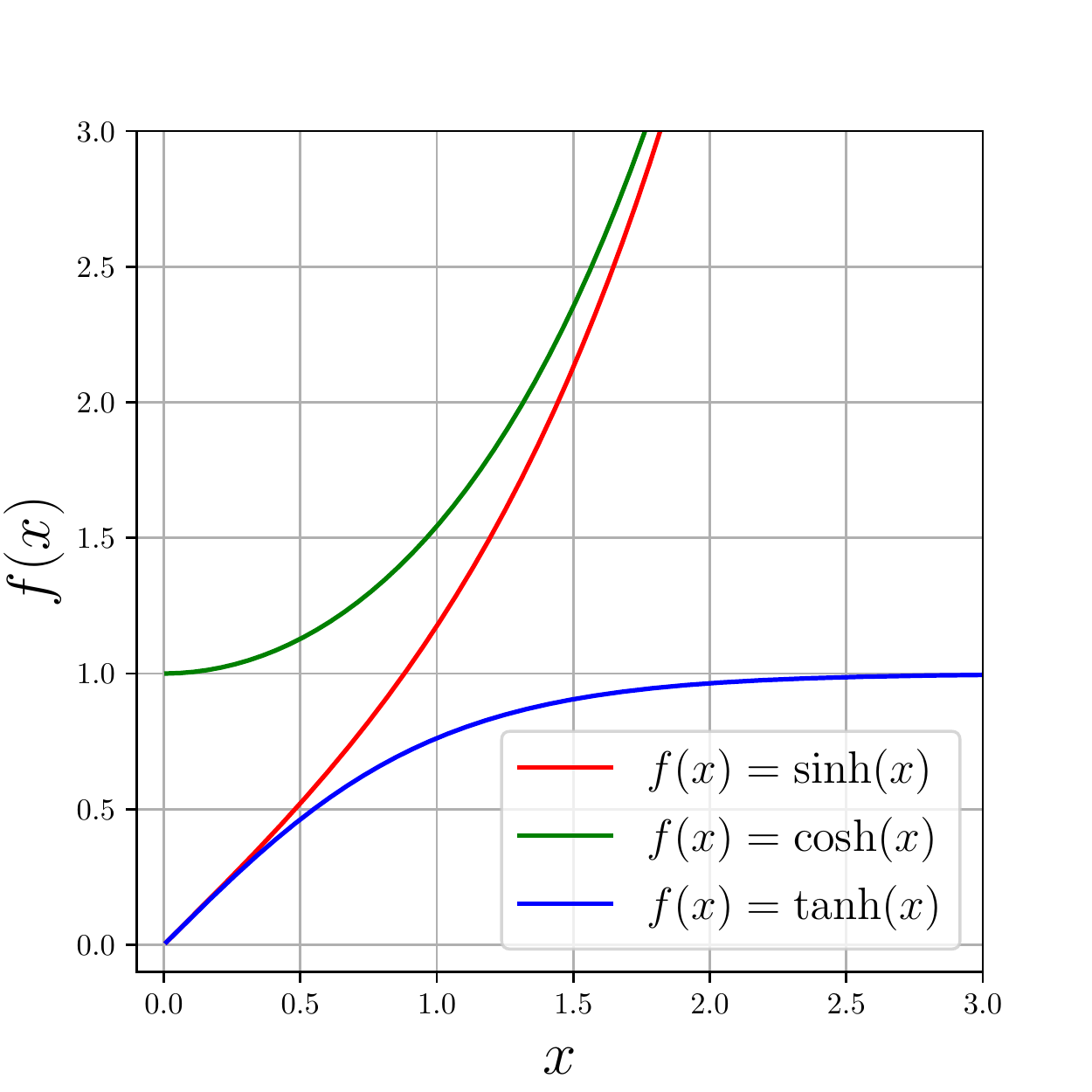}\label{fig:hyp1}}
	\subfigure[$\coth$, $\sech$, $\csch$]{\includegraphics[width=0.32\textwidth]{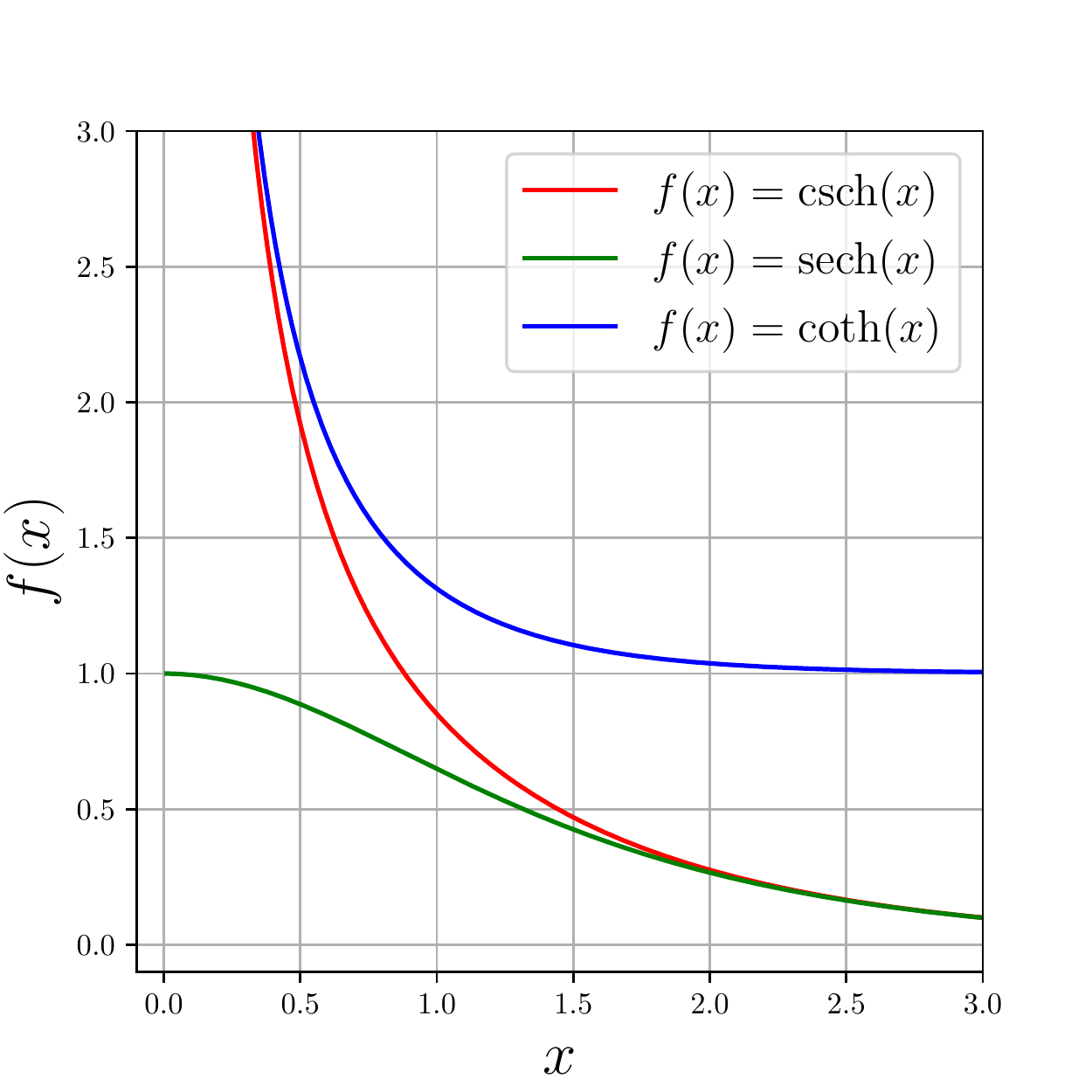}\label{fig:hyp2}}
	\subfigure[$\sinhc$, $\tanhc$, $\cothc$, $\cschc$]{\includegraphics[width=0.32\textwidth]{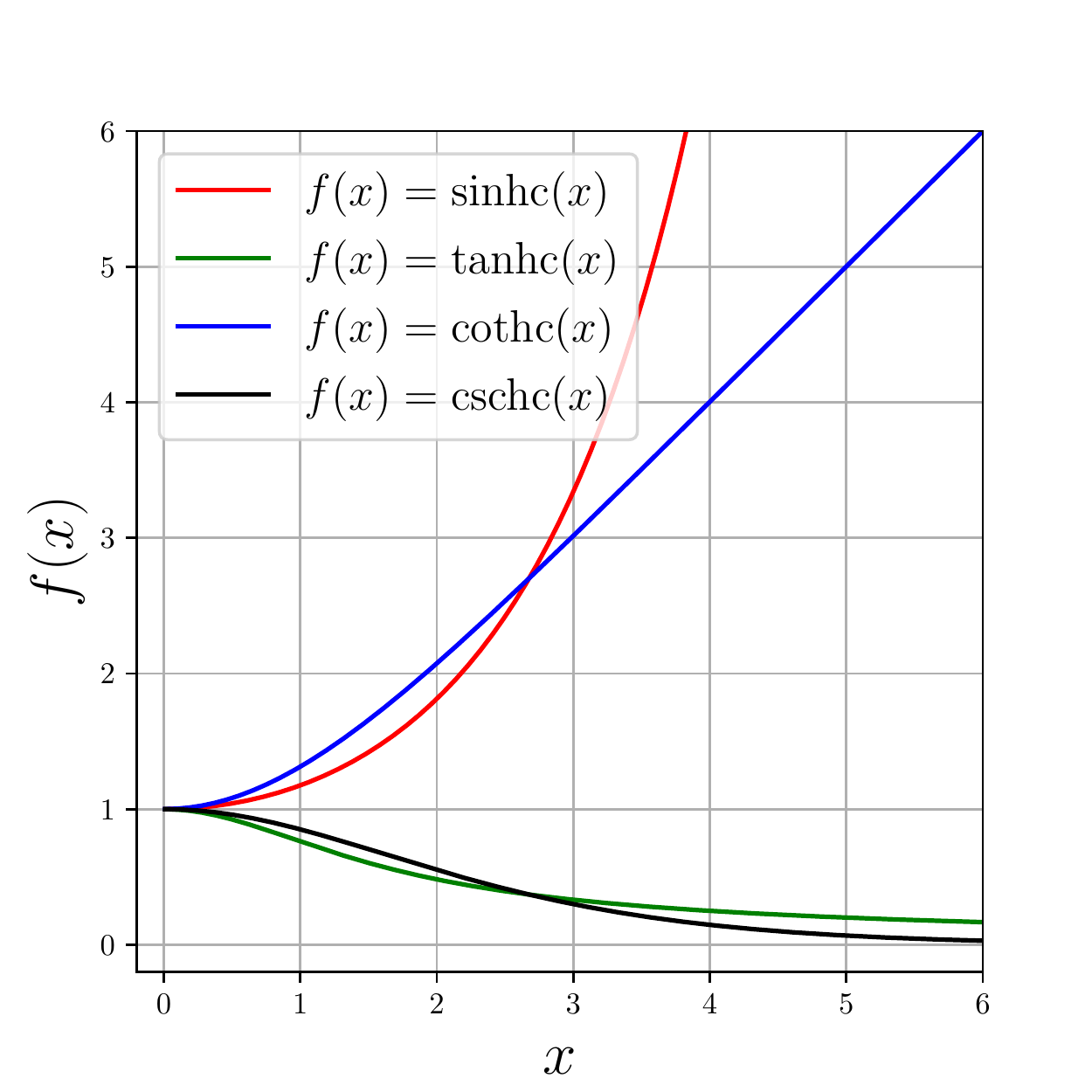}\label{fig:hyp3}}
	\caption{Hyperbolic functions and their variants.}\label{fig:hyp}
\end{figure}

\paragraph{Higher-order hyperbolic functions.}
We now define the higher-order hyperbolic functions that  will be used to design higher-order accelerated optimization algorithms. We define the $p$-th order hyperbolic sine function $\sinh_p:[0,\infty)\rightarrow\mathbb{R}$ as the solution of the initial value problem
\begin{equation}
	\label{eq:def_sinhp}
	\sinh_p'(t)=\cosh_p(t):=\left(1+\sinh_p^{p}(t)\right)^{1/p},\quad \sinh_p(0)=0.
\end{equation}
Furthermore, we define the $\tanh_p$, $\coth_p$, $\sech_p$, and $\csch_p$ functions as
\begin{align*}
	\tanh_{p}(t) & =\frac{\sinh_{p}(t)}{\cosh_{p}(t)}, & \coth_{p}(t) & =\frac{1}{\tanh_{p}(t)}, & \sech_{p}(t) & =\frac{1}{\sinh_{p}(t)}, & \csch_{p}(t) & =\frac{1}{\cosh_{p}(t)}.
\end{align*}
We define the $\sinhc_p$, $\tanhc_p$, $\cothc_p$, and $\cschc_p$ functions as
\begin{align*}
	\sinhc_{p}x & :=\begin{cases}
		\frac{\sinh_{p}x}{x}, & \text{if }x\neq0\\
		1, & \text{if }x=0
	\end{cases} & \tanhc_{p}x & :=\frac{\sinhc_{p}x}{\cosh_{p}x}, \\
	\cothc_{p}x & :=\frac{1}{\tanhc_{p}x}, & \cschc_{p}x & :=\frac{1}{\sinhc_{p}x}.
\end{align*}
Note that the higher-order hyperbolic functions recover the usual hyperbolic functions when $p=2$. The following proposition says that the $\sinh_p$ function grows exponentially.
\begin{proposition}
	\label{prop:sinhp_grows_exponentially}
	There is a constant $C_p>0$ such that $\sinh_p(t) \sim C_p e^t$ as $t\rightarrow\infty$. In particular, we have $C_p=1/2$ for $p=2$.
\end{proposition}
The proof of Proposition~\ref{prop:sinhp_grows_exponentially} can be found in Appendix~\ref{app:prop:sinhp_grows_exponentially}. Using \eqref{eq:def_sinhp} and Proposition~\ref{prop:sinhp_grows_exponentially}, it is straightforward to check the following asymptotic properties:
\begin{align*}
	\sinh_{p}x & \sim x\text{ as }x\rightarrow0, & \sinh_{p}x & \sim C_{p}e^{x}\text{ as }x\rightarrow\infty\\
	\cosh_{p}x & \sim1\text{ as }x\rightarrow0, & \cosh_p x & \sim C_{p}e^{x}\text{ as }x\rightarrow\infty\\
	\tanh_{p}x & \sim x\text{ as }x\rightarrow0, & \tanh_p x & \sim1\text{ as }x\rightarrow\infty.
\end{align*}

\subsection{Limiting arguments and examples}
\label{sec:limiting}

We investigate two additional ways to derive the limiting ODEs of first-order algorithms. The first approach is to write the algorithm as a two-sequence scheme and then derive the limiting ODE via the second-order Taylor series expansion. This argument frequently appears in the literature \citep[see][]{su2016differential,shi2021understanding}. The second approach, which is novel, is to express the algorithm using the \emph{difference matrix} $\mathbf{H}=(h_{ij})$ and then derive the \emph{differential kernel} $H(t,\tau)$ corresponding to the matrix $(h_{ij})$. We only present the results here and defer the detailed computations to Appendices~\ref{app:two-sequence_scheme} and \ref{app:differential_kernel}.

\subsubsection{Limiting ODEs of two-sequence algorithms}
\label{sec:limiting_two}

We consider the following \emph{two-sequence scheme}:
\begin{equation}
	\label{eq:two_sequence_scheme}
	\begin{aligned}
		x_{k+1} & =y_{k}-s\nabla f\left(y_{k}\right)\\
		y_{k+1} & =x_{k+1}+\beta_{k}\left(x_{k+1}-x_{k}\right)+\gamma_{k}\left(x_{k+1}-y_{k}\right).
	\end{aligned}
\end{equation}
If we have
\begin{equation}
	\label{eq:two_seq_assumptions}
	\lim_{s\rightarrow0}\frac{1-\beta_{t/\sqrt{s}}}{\sqrt{s}}=b(t) \textrm{ and }\lim_{s\rightarrow0}\gamma_{t/\sqrt{s}}=c(t) \textrm{ for all }t>0
\end{equation}
for some smooth functions $b,c:(0,\infty)\rightarrow\mathbb{R}$, then under the identification $X(t_k)=x_k$ with $t_k=k\sqrt{s}$, the two-sequence scheme \eqref{eq:two_sequence_scheme} converges to the ODE
\begin{equation}
	\label{eq:two_sequence_limiting_ODE}
	\ddot{X}(t)+b(t)\dot{X}(t)+(1+c(t))\nabla f(X(t))=0
\end{equation}
as $s\rightarrow0$.

\paragraph{Recovering the limiting ODE of three-sequence scheme.}
We can write the three-sequence scheeme \eqref{eq:nag} as the two-sequence scheme \eqref{eq:two_sequence_scheme} with the following parameters \citep[see][Appendix~B]{lee2021geometric}:
\begin{equation}
	\label{eq:beta_three}
	\begin{aligned}
		\beta_{k} & =\frac{\left(1-\tau_{k}\right)\tau_{k+1}\left(1-\mu\delta_{k}\right)}{\tau_{k}}\\
		\gamma_{k} & =\frac{\tau_{k+1}\left((1/s-\mu)\delta_{k}\tau_{k}-1+\mu\delta_{k}\right)}{\tau_{k}}.
	\end{aligned}
\end{equation}
If the limits \eqref{eq:taut_deltat} with $\bft_k=k\sqrt{s}$ exist, then we have
\begin{align*}
	\lim_{s\rightarrow0}\frac{1-\beta_{t/\sqrt{s}}}{\sqrt{s}} & =\tau(t)-\frac{\dot{\tau}(t)}{\tau(t)}+\mu\delta(t)\\
	\lim_{s\rightarrow0}\gamma_{t/\sqrt{s}} & =\tau(t)\delta(t)-1
\end{align*}
for all $t>0$. Therefore, we recover the limiting ODE \eqref{eq:limiting_ode_oneline} of the three-sequence scheme. In particular, if the algorithmic parameters $(\tau_k)$ and $(\delta_k)$ satisfy the collinearity condition \eqref{eq:collinear}, then we have $\gamma_k=0$ for all $k\geq0$, and thus $c(t)=0$.

\paragraph{Two-sequence form of \nagc.}
Because \nagcsp is the three-sequence scheme \eqref{eq:nag} with $\tau_{k}=\frac{2}{k+1}$, $\delta_{k}=\frac{s(k+1)}{2}$, and $\mu=0$, we can rewrite it as the two-sequence scheme \eqref{eq:two_sequence_scheme} with
\begin{align*}
	\beta_{k} & =\frac{\left(1-\frac{2}{k+1}\right)\frac{2}{k+2}}{\frac{2}{k+1}}=\frac{k-1}{k+2}\\
	\gamma_{k} & =\frac{\frac{2}{k+2}\cdot\frac{s(k+1)}{2}}{s}-\frac{\frac{2}{k+2}}{\frac{2}{k+1}}=0.
\end{align*}
Thus, \nagcsp converges to the ODE \eqref{eq:two_sequence_limiting_ODE} with
\begin{align*}
	b(t) & =\lim_{s\rightarrow0}\frac{1-\frac{t/\sqrt{s}-1}{t/\sqrt{s}+2}}{\sqrt{s}}=\frac{3}{t}\\
	c(t) & =0,
\end{align*}
which recovers \nagcsp ODE \eqref{eq:c_ode}.

\paragraph{Two-sequence form of \nagsc.}
Because \nagscsp is the three-sequence scheme \eqref{eq:nag} with $\tau_{k}=\frac{\sqrt{\mu s}}{1+\sqrt{\mu s}}$ and $\delta_{k}=\sqrt{\frac{s}{\mu}}$, it can be written as the two-sequence scheme \eqref{eq:two_sequence_scheme} with
\begin{align*}
	\beta_{k} & =\frac{\left(1-\frac{\sqrt{\mu s}}{1+\sqrt{\mu s}}\right)\frac{\sqrt{\mu s}}{1+\sqrt{\mu s}}\left(1-\mu\sqrt{\frac{s}{\mu}}\right)}{\frac{\sqrt{\mu s}}{1+\sqrt{\mu s}}}=\frac{1-\sqrt{\mu s}}{1+\sqrt{\mu s}}\\
	\gamma_{k} & =\frac{\frac{\sqrt{\mu s}}{1+\sqrt{\mu s}}\sqrt{\frac{s}{\mu}}}{s}-\left(1-\mu\sqrt{\frac{s}{\mu}}+\mu\frac{\sqrt{\mu s}}{1+\sqrt{\mu s}}\sqrt{\frac{s}{\mu}}\right)=0.
\end{align*}
Thus \nagscsp converges to the ODE \eqref{eq:two_sequence_limiting_ODE} with
\begin{align*}
	b(t) & =\lim_{s\rightarrow0}\frac{1-\frac{1-\sqrt{\mu s}}{1+\sqrt{\mu s}}}{\sqrt{s}}=2\sqrt{\mu}\\
	c(t) & =0,
\end{align*}
which recovers \nagscsp ODE \eqref{eq:sc_ode}.

\subsubsection{Difference matrices and differential kernels}
\label{sec:limiting_kernel}

We can formulate most of the practical first-order momentum methods as the following  \emph{fixed-step first-order scheme} \citep[see][]{drori2014performance}:
\begin{equation}
	\label{eq:discrete_FSFO}
	y_{i+1}=y_{i}-s\sum_{j=0}^{i}h_{ij}\nabla f\left(y_{j}\right)\textrm{ for }i=0,\ldots,N-1,
\end{equation}
where $N$ is the number of iterations. We can write this scheme equivalently as 
\[
\left[\begin{array}{c}
	y_{1}-y_{0}\\
	y_{2}-y_{1}\\
	\vdots\\
	y_{N}-y_{N-1}
\end{array}\right]=-s\left[\begin{array}{cccc}
	h_{0,0} & 0 & \cdots & 0\\
	h_{1,0} & h_{1,1} & \cdots & 0\\
	\vdots & \vdots & \ddots & \vdots\\
	h_{N-1,0} & h_{N-1,2} & \cdots & h_{N-1,N-1}
\end{array}\right]\left[\begin{array}{c}
	\nabla f\left(y_{0}\right)\\
	\nabla f\left(y_{1}\right)\\
	\vdots\\
	\nabla f\left(y_{N-1}\right)
\end{array}\right]
\]
Here, we call the lower triangular matrix $\mathbf{H}=(h_{ij})$ the \emph{difference matrix} for the algorithm~\eqref{eq:discrete_FSFO}.

To derive the limiting ODE of the algorithm \eqref{eq:discrete_FSFO}, we introduce a smooth curve $X:[0,T]\rightarrow\mathbb{R}^n$ with the identifications $X(k\sqrt{s})= y_k$ and $T=N\sqrt{s}$. As a continuous-time analog of the difference matrix $(h_{ij})$, we intoduce a continuously differentiable function $H$ (possibly depending on $s$ now) defined on $\{(t,\tau)\in\mathbb{R}^{2}:0<\tau\leq t< T\}$ with the identification $H(t_i,\tau_j)=h_{ij}$, where $t_i=i\sqrt{s}$ and $\tau_j=j\sqrt{s}$. Substituting $X(t_i)=y_i$ in \eqref{eq:discrete_FSFO} yields
\begin{equation}
	\label{eq:discrete_FSFO2}
	\frac{X\left(t_{i+1}\right)-X\left(t_{i}\right)}{\sqrt{s}}=-\left(\tau_{j+1}-\tau_j\right)\sum_{j=0}^{i}H\left(t_{i},\tau_{j}\right)\nabla f\left(X\left(\tau_{j}\right)\right).
\end{equation}
Then, we can observe that the right-hand side of \eqref{eq:discrete_FSFO2} is a Riemann sum of the function $\tau\mapsto -H(t_i,\tau)\nabla f(X(\tau))$ over $[0,t_{i+1}]$. Thus, taking the limit $s\rightarrow0$ yields
\begin{equation}
	\label{eq:continuous_FSFO}
	\dot{X}(t)=-\int_{0}^{t}H(t,\tau)\nabla f(X(\tau))\,d\tau, \textrm{ where }H(t,\tau)=\lim_{s\rightarrow0}h_{\frac{t}{\sqrt{s}},\frac{\tau}{\sqrt{s}}}
\end{equation}
as the limiting ODE of the fixed-step first-order scheme \eqref{eq:discrete_FSFO}. Note that the form of this equation clearly reflects the momentum effect because the gradient $\nabla f(X(\tau))$ at time $\tau$ affects the velocity $\dot{X}(t)$ at all times $t$ after $\tau$. Inspired by the observation that the function $H(t,\tau)$ plays a role similar to the \emph{kernel function} in the integral transform, we call it the \emph{differential kernel} (or the \emph{H-kernel}) corresponding to the difference matrix $(h_{ij})$.

\paragraph{From differential kernels to second-order ODEs.}
Differentiating both sides of \eqref{eq:continuous_FSFO} and applying the Leibniz integral rule, we obtain
\begin{equation}
	\label{eq:2nd_ode_kernel}
	\ddot{X}(t)=-H(t,t)\nabla f(X(t))-\int_{0}^{t}\frac{\partial H(t,\tau)}{\partial t}\nabla f(X(\tau))\,d\tau.
\end{equation}
If there exists a function $b(t)$ such that 
\[
\frac{\partial H(t,\tau)}{\partial t}=-b(t)H(t,\tau),
\]
then it follows from \eqref{eq:continuous_FSFO} that the equation \eqref{eq:2nd_ode_kernel} is expressed as the following second-order ODE:
\begin{equation}
	\label{eq:2nd_ode_kernel2}
	\ddot{X}(t)+b(t)\dot{X}+H(t,t)\nabla f(X(t))=0.
\end{equation}

\paragraph{Recovering the limiting ODE of two-sequence scheme.}
We can write the two-sequence scheme \eqref{eq:two_sequence_scheme} as the fixed-step first-order scheme with
\[
h_{ij}=\left(\beta_{j}+\gamma_{j}\right)\prod_{\nu=j+1}^{i}\beta_{\nu}+\delta_{ij},
\]
where $\delta_{ij}$ is the Kronecker delta funciton. For $i>j$,\footnote{We exclude the case $i=j$ because the difference matrix $h_{ij}$ has singularities at these points due to the Kronecker delta function.} we have
\[
h_{i+1,j}-h_{i,j}=\left(\beta_{i+1}-1\right)h_{ij}.
\]
Under the identification $H(t_i,\tau_j)=h_{ij}$, we have
\[
h_{i+1,j}-h_{i,j}=H\left(t_{i+1},\tau_{j}\right)-H\left(t_{i},\tau_{j}\right)=\frac{\partial H\left(t_{i},\tau_{j}\right)}{\partial t}\sqrt{s}+o\left(\sqrt{s}\right).
\]
Thus, when the limits \eqref{eq:two_seq_assumptions} exist, taking the limit $s\rightarrow0$ yields
\begin{equation}
	\label{eq:partial_derivative}
	\frac{\partial H\left(t,\tau\right)}{\partial t}=-b(t)H\left(t,\tau\right).
\end{equation}
Also, because $h_{k+1,k}=\beta_{k+1}+\gamma_{k}$ and $\lim_{s\rightarrow0}\beta_{t/\sqrt{s}}=1$ by \eqref{eq:two_seq_assumptions}, we have $H(t,t)=1+c(t)$ for all $t\in(0,T)$. Therefore, the ODE \eqref{eq:2nd_ode_kernel2} recovers the limiting ODE \eqref{eq:two_sequence_limiting_ODE} of the two-sequence scheme. Moreover, we can explicitly write the differential kernel $H$ as
\begin{equation}
	\label{eq:kernel_for_two}
	H(t,\tau)=\left(1+c(\tau)\right)e^{-\int_{\tau}^{t}b(s)\,ds}.
\end{equation}

\paragraph{Difference matrix for \nagc.}

Because we can write \nagcsp as the two-sequence scheme \eqref{eq:two_sequence_scheme} with $\beta_k=\frac{k-1}{k+2}$ and $\gamma_{k}=0$, we can rewrite it as the fixed-step first-order scheme \eqref{eq:discrete_FSFO} with
\[
h_{ij}=\prod_{\nu=j}^{i}\frac{\nu-1}{\nu+2}+\delta_{ij}=\frac{(j-1)j(j+1)}{i(i+1)(i+2)}+\delta_{ij}.
\]
By definition, the differential kernel corresponding to this matrix $(h_{ij})$ is
\begin{equation}
	\label{eq:kernel_c}
	H(t,\tau)=\lim_{s\rightarrow0}\frac{\left(\frac{\tau}{\sqrt{s}}-1\right)\frac{\tau}{\sqrt{s}}\left(\frac{\tau}{\sqrt{s}}+1\right)}{\frac{t}{\sqrt{s}}\left(\frac{t}{\sqrt{s}}+1\right)\left(\frac{t}{\sqrt{s}}+2\right)}=\frac{\tau^{3}}{t^{3}}.
\end{equation}
This can be also obtained by substituting $b(t)=3/t$ and $c(t)=0$ into \eqref{eq:kernel_for_two}:
\[
H(t,\tau)=e^{-\int_{\tau}^{t}\frac{3}{s}\,ds}=e^{-3(\log(t)-\log(\tau))}=\frac{\tau^{3}}{t^{3}}.
\]
Because
\[
\frac{\partial H(t,\tau)}{\partial t}=-\frac{3\tau^{3}}{t^{4}}=-\frac{3}{t}H(t,\tau),
\]
the ODE \eqref{eq:2nd_ode_kernel} with \eqref{eq:kernel_c} recovers \nagcsp ODE \eqref{eq:c_ode}.

\paragraph{Difference matrix for \nagsc.}

Because we can write \nagscsp as the two-sequence scheme \eqref{eq:two_sequence_scheme} with $\beta_k=\frac{1-\sqrt{\mu s}}{1+\sqrt{\mu s}}$ and $\gamma_{k}=0$, we can rewrite it as the fixed-step first-order scheme \eqref{eq:discrete_FSFO} with
\[
h_{ij}=\prod_{\nu=j}^{i}\frac{1-\sqrt{\mu s}}{1+\sqrt{\mu s}}+\delta_{ij}=\left(\frac{1-\sqrt{\mu s}}{1+\sqrt{\mu s}}\right)^{i-j+1}+\delta_{ij}.
\]
By definition, the differential kernel corresponding to this matrix $(h_{ij})$ is
\begin{equation}
	\label{eq:kernel_sc}
	H(t,\tau)=\lim_{s\rightarrow0}\left(\frac{1-\sqrt{\mu s}}{1+\sqrt{\mu s}}\right)^{\frac{t}{\sqrt{s}}-\frac{\tau}{\sqrt{s}}+1}=\frac{e^{2\sqrt{\mu}\tau}}{e^{2\sqrt{\mu}t}}.
\end{equation}
This can be also obtained by substituting $b(t)=2\sqrt{\mu}$ and $c(t)=0$ into \eqref{eq:kernel_for_two}:
\[
H(t,\tau)=e^{-\int_{\tau}^{t}2\sqrt{\mu}\,ds}=e^{-2\sqrt{\mu}(t-\tau)}=\frac{e^{2\sqrt{\mu}\tau}}{e^{2\sqrt{\mu}t}}.
\]
It follows from
\[
\frac{\partial H(t,\tau)}{\partial t}=-2\sqrt{\mu}e^{2\sqrt{\mu}(\tau-t)}=-2\sqrt{\mu}H(t,\tau)
\]
that the ODE \eqref{eq:2nd_ode_kernel} with \eqref{eq:kernel_sc} recovers \nagscsp ODE \eqref{eq:c_ode}.

\section{Unified Bregman Lagrangian}
\label{sec:unified_lagrangian}

In this section, we address the inconsistency between the first Bregman Lagrangian \eqref{eq:first_lagrangian} and the second Bregman Lagrangian \eqref{eq:second_lagrangian}. For a continuously differentiable strictly convex function $h$, we define the \emph{unified Bregman Lagrangian} as
\begin{equation}
	\label{eq:unified_lagrangian}
	\begin{aligned}
		\mathcal{L}\left(X,\dot{X},t\right) & =\mathcal{L}_{\text{1st}}\left(X,\dot{X},t\right)+\mathcal{L}_{\text{2nd}}\left(X,\dot{X},t\right)-\left[\mathcal{L}_{\text{2nd}}\left(X,\dot{X},t\right)\right]_{\mu=0}\\
		& =e^{\alpha+\gamma}\left(\left(1+\mu e^{\beta}\right)D_{h}\left(X+e^{-\alpha}V,X\right)-e^{\beta}f(X)\right).
	\end{aligned}
\end{equation}
Then by construction, this Lagrangian recovers the first Bregman Lagrangian \eqref{eq:first_lagrangian} when $\mu=0$. Because the Lagrangian \eqref{eq:unified_lagrangian} is continuous in the strong convexity parameter $\mu$, it is a continuous extension of the first Bregman Lagrangian to the strongly convex case.
\begin{proposition}
	\label{prop:e-l-eq}
	Under the ideal scaling condition \eqref{eq:ideal_scaling_a}, the Euler--Lagrange equation \eqref{eq:e-l-eq} for the unified Bregman Lagrangian \eqref{eq:unified_lagrangian} reduces to the following system of ODEs:
	\begin{subequations}
		\label{eq:unified_family}
		\begin{align}
			\label{eq:unified_family_a}
			\dot{X} & =e^{\alpha}(Z-X)\\
			\label{eq:unified_family_b}
			\frac{d}{dt}\nabla h(Z) & =\frac{\mu\dot{\beta}e^{\beta}}{1+\mu e^{\beta}}\left(\nabla h(X)-\nabla h(Z)\right)-\frac{e^{\alpha+\beta}}{1+\mu e^{\beta}}\nabla f(X).
		\end{align}
	\end{subequations}
\end{proposition}
The proof of Proposition~\ref{prop:e-l-eq} can be found in Appendix~\ref{app:prop:e-l-eq}.  
To analyze the convergence rate of this dynamics, we define the time-dependent Lyapunov function $V:\mathbb{R}^n\times\mathbb{R}^n\times[\tst,\infty)\rightarrow\mathbb{R}$ as
\begin{equation}
	\label{eq:lyapunov_lagrangian}
	V(X,Z,t)=\left(1+\mu e^{\beta(t)}\right)D_{h}\left(x^{*},Z\right)+e^{\beta(t)}\left(f(X)-f\left(x^{*}\right)\right).
\end{equation}

\begin{theorem}
	\label{thm:mainthm_lagrangian}
	Suppose that the ideal scaling condition \eqref{eq:ideal_scaling_b} holds. Let $f$ be a $\mu$-uniformly (possibly with $\mu=0$) convex function with respect to $h$. Then, for any solution $(X,Z)$ to the unified Bregman Lagrangian flow \eqref{eq:unified_family}, the continuous-time energy function
	\[
	\mathcal{E}(t)=V(X(t),Z(t),t)=\left(1+\mu e^{\beta(t)}\right)D_{h}\left(x^{*},Z(t)\right)+e^{\beta(t)}\left(f(X(t))-f\left(x^{*}\right)\right)
	\]
	is monotonically decreasing on $[\tst,\infty)$.
\end{theorem}
The proof of Theorem~\ref{thm:mainthm_lagrangian} can be found in Appendix~\ref{app:thm:mainthm_lagrangian}. Writing $\mathcal{E}(t)\leq\mathcal{E}(\tst)$ explicitly, we obtain an $O(e^{-\beta(t)})$ convergence rate for the dynamics \eqref{eq:unified_family}.
\begin{corollary}
	\label{cor:maincor_Lagrangian}
	Suppose that the ideal scaling condition \eqref{eq:ideal_scaling_b} holds. Let $f$ be a $\mu$-uniformly (possibly with $\mu=0$) convex function with respect to $h$. Then, any solution $(X,Z)$ to the unified Bregman Lagrangian flow \eqref{eq:unified_family} satisfies the inequality
	\[
	f(X(t))-f\left(x^{*}\right)\leq e^{-\beta(t)}\left(\left(1+\mu e^{\beta(\tst)}\right)D_{h}\left(x^{*},Z(\tst)\right)+e^{\beta(\tst)}\left(f\left(X(\tst)\right)-f\left(x^{*}\right)\right)\right)
	\]
	for all $t>0$.
\end{corollary}
Similarly to the first Bregman Lagrangian flow \eqref{eq:first_family} and the second Bregman Lagrangian flow \eqref{eq:second_family} \citep[see][]{wibisono2016,wilson2021lyapunov}, the dynamical system \eqref{eq:unified_family} is closed under time-dilation.
\begin{theorem}
	\label{thm:time_dilation}
	Let $\bfT:I_2\rightarrow I_1$ be an increasing continuously differentiable bijective function, where $I_1$ and $I_2$ are intervals in $[0,\infty)$. If $(X_1,Z_1)$ is a solution to the unified Bregman Lagrangian flow \eqref{eq:unified_family} on $I_1$ with parameters $\alpha_1,\beta_1:I_1\rightarrow\mathbb{R}$, then the reparametrized curves $X_2(t)=X_1(\bfT(t))$ and $Z_2(t)=Z_1(\bfT(t))$ is a solution to the unified Bregman Lagrangian flow on $I_2$ with the parameters $\alpha_2,\beta_2:I_2\rightarrow\mathbb{R}$ defined by
	\begin{align*}
		{\alpha_2}(t) & =\alpha_1(\bfT(t))+\log\dot{\bfT}(t)\\
		{\beta_2}(t) & =\beta_1(\bfT(t)).
	\end{align*}
\end{theorem}
The proof of Theorem~\ref{thm:time_dilation} can be found in Appendix~\ref{app:thm:time_dilation}.

\paragraph{Recovering the first and second Bregman Lagrangians.}

We now discuss how the first Bregman Lagrangian flow \eqref{eq:first_family}, the second Bregman Lagrangian flow \eqref{eq:second_family}, and the corresponding Lyapunov analyses can be recovered from the proposed unified Bregman Lagrangian flow \eqref{eq:unified_family} and the corresponding Lyapunov analysis (Theorem~\ref{thm:mainthm_lagrangian}). 
When $\mu=0$, it is easy to check that the unified Bregman Lagrangian flow and the corresponding Lyapunov function \eqref{eq:lyapunov_lagrangian} recover the first Bregman Lagrangian flow and the corresponding Lyapunov function \eqref{eq:first_lagrangian_lyapunov}. When the limits $\alpha(\infty):=\lim_{t\rightarrow\infty}\alpha(t)$ and $\dot{\beta}(\infty):=\lim_{t\rightarrow\infty}\dot{\beta}(t)>0$ exist, the second Bregman Lagrangian flow with ${\alpha}_{\mathrm{2nd}}(t):\equiv\alpha(\infty)$ and ${\beta}_{\mathrm{2nd}}(t):=\dot{\beta}(\infty)t$ is the \emph{asymptotic version} of the unified Bregman Lagrangian flow with $\alpha(t)$ and $\beta(t)$ in the sense that the coefficients of \eqref{eq:unified_family} converge to the ones of \eqref{eq:second_family} as $t\rightarrow\infty$. In Appendix~\ref{app:recover_lagrangian}, we show that the Lyapunov analysis for the second Bregman Lagrangian flow with $\tilde{\alpha}$ and $\tilde{\beta}$ can be recovered from Theorem~\ref{thm:mainthm_lagrangian} by taking the limit $t\rightarrow\infty$ of some inequalities.

\section{Unified Methods for Minimizing Convex and Strongly Convex  Functions}
\label{sec:unified_sections}

In Section~\ref{sec:unified_ode}, we address the inconsistency between \nagcsp ODE \eqref{eq:c_ode} and \nagscsp ODE \eqref{eq:sc_ode} by proposing an ODE model that unifies \nagcsp ODE and \nagscsp ODE. In Section~\ref{sec:unified_nag}, we address the inconsistency between \nagcsp \eqref{eq:nag-c} and \nagscsp \eqref{eq:nag-sc} by proposing novel algorithms that can be viewed as a discrete-time counterpart of the unified NAG ODE. Throughout this section, we assume the standard smooth strongly convex setting in Section~\ref{sec:convexity}.

\subsection{Proposed dynamics: Unified NAG ODE}
\label{sec:unified_ode}

We consider the unified Bregman Lagrangian flow \eqref{eq:unified_family} with $\alpha(t)=\log(\frac{2}{t}\cothc(\frac{\sqrt{\mu}}{2}t))$, $\beta(t)=\log(\frac{t^{2}}{4}\sinhc^{2}(\frac{\sqrt{\mu}}{2}t))$,\footnote{We can constructively choose these functions (see Appendix~\ref{app:alphabeta}). Note that when $\mu=0$, we have $\alpha(t)=\log\frac{2}{t}$ and $\beta(t)=\log\frac{t^{2}}{4}$, which recover \nagcsp ODE \eqref{eq:c_ode} from the first Bregman Lagrangian flow \eqref{eq:first_family}. Also, as $t\rightarrow\infty$, we have $\alpha(t)\sim\log\sqrt{\mu}$ and $\beta(t)\sim\sqrt{\mu}t-\log\left(4\mu\right)$, which recover \nagscsp ODE \eqref{eq:sc_ode} from the second Bregman Lagrangian flow \eqref{eq:second_family}.} $h(x)=\frac{1}{2}\left\Vert x\right\Vert ^{2}$, and the initial conditions $X(\tst)=Z(\tst)=x_0$, which we call the \emph{unified NAG system}:
\begin{equation}
	\label{eq:u_ode}
	\begin{aligned}
		\dot{X} & =\frac{2}{t}\cothc\left(\frac{\sqrt{\mu}}{2}t\right)(Z-X)\\
		\dot{Z} & =\frac{t}{2}\tanhc\left(\frac{\sqrt{\mu}}{2}t\right)\left(\mu X-\mu Z-\nabla f(X)\right).
	\end{aligned}
\end{equation}
Writing this system in a single equation, we obtain the \emph{unified NAG ODE} (see Appendix~\ref{app:equivalent_of_ode}):
\begin{equation}
	\label{eq:u_ode_oneline}
	\ddot{X}+\left(\frac{\sqrt{\mu}}{2}\tanh\left(\frac{\sqrt{\mu}}{2}t\right)+\frac{3}{t}\cothc\left(\frac{\sqrt{\mu}}{2}t\right)\right)\dot{X}+\nabla f(X)=0
\end{equation}
with $X(\tst)=x_0$ and $\dot{X}(\tst)=0$. 

\paragraph{Existence and uniqueness of the solution.}
To prove the existence and uniqueness of solution to the unified NAG system \eqref{eq:u_ode}, we cannot directly apply the classical existence and uniqueness theorem because the coefficient $\frac{2}{t}\cothc\left(\frac{\sqrt{\mu}}{2}t\right)$ has a singularity at $t=0$. Thus, we follow the arguments in \citep{krichene2015accelerated,su2016differential}.
\begin{theorem}
	\label{thm:existence_unique_unified}
	The unified NAG system \eqref{eq:u_ode} has a unique solution $(X,Z)$ in $C^1([0,\infty),\mathbb{R}^n\times\mathbb{R}^n)$.
\end{theorem}
The proof of Theorem~\ref{thm:existence_unique_unified} can be found in Appendix~\ref{app:existence_uniqueness_unified_nag}.

\begin{figure}[ht]
	\centering
	\subfigure[$\mu=0.3$]{\includegraphics[width=0.32\textwidth]{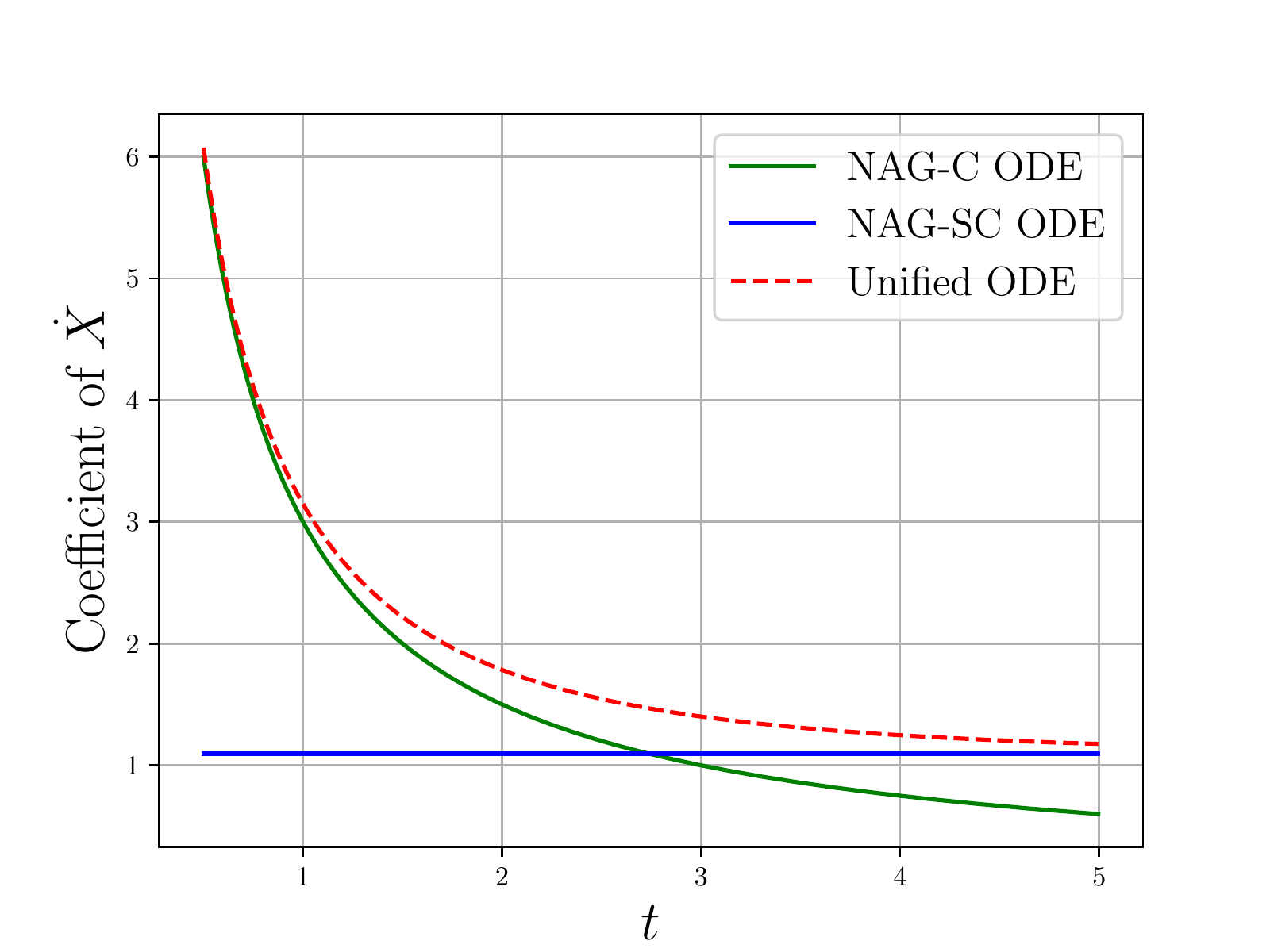}\label{fig:ode_coeff1}}
	\subfigure[$\mu=1$]{\includegraphics[width=0.32\textwidth]{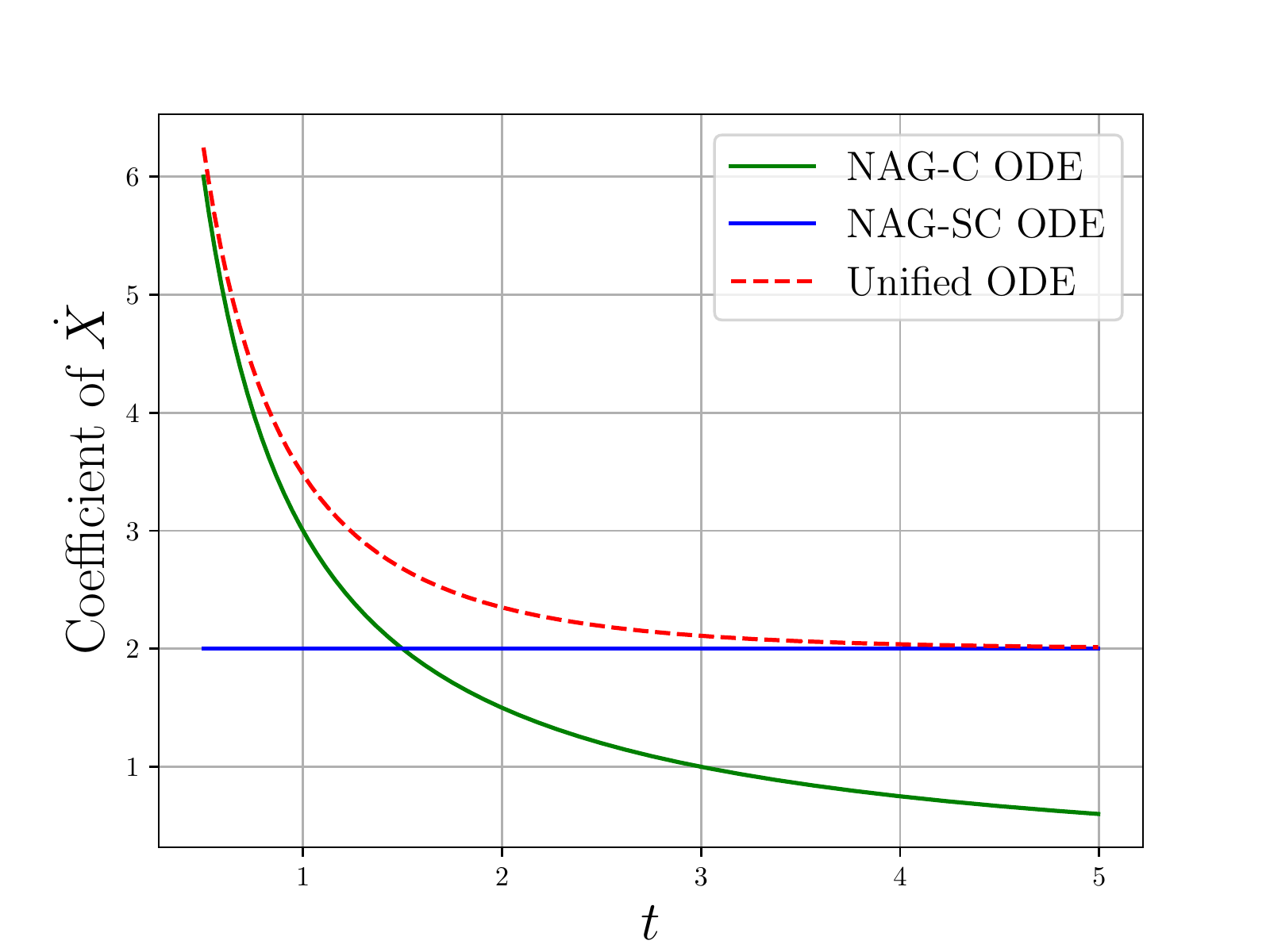}\label{fig:ode_coeff2}}
	\subfigure[$\mu=5$]{\includegraphics[width=0.32\textwidth]{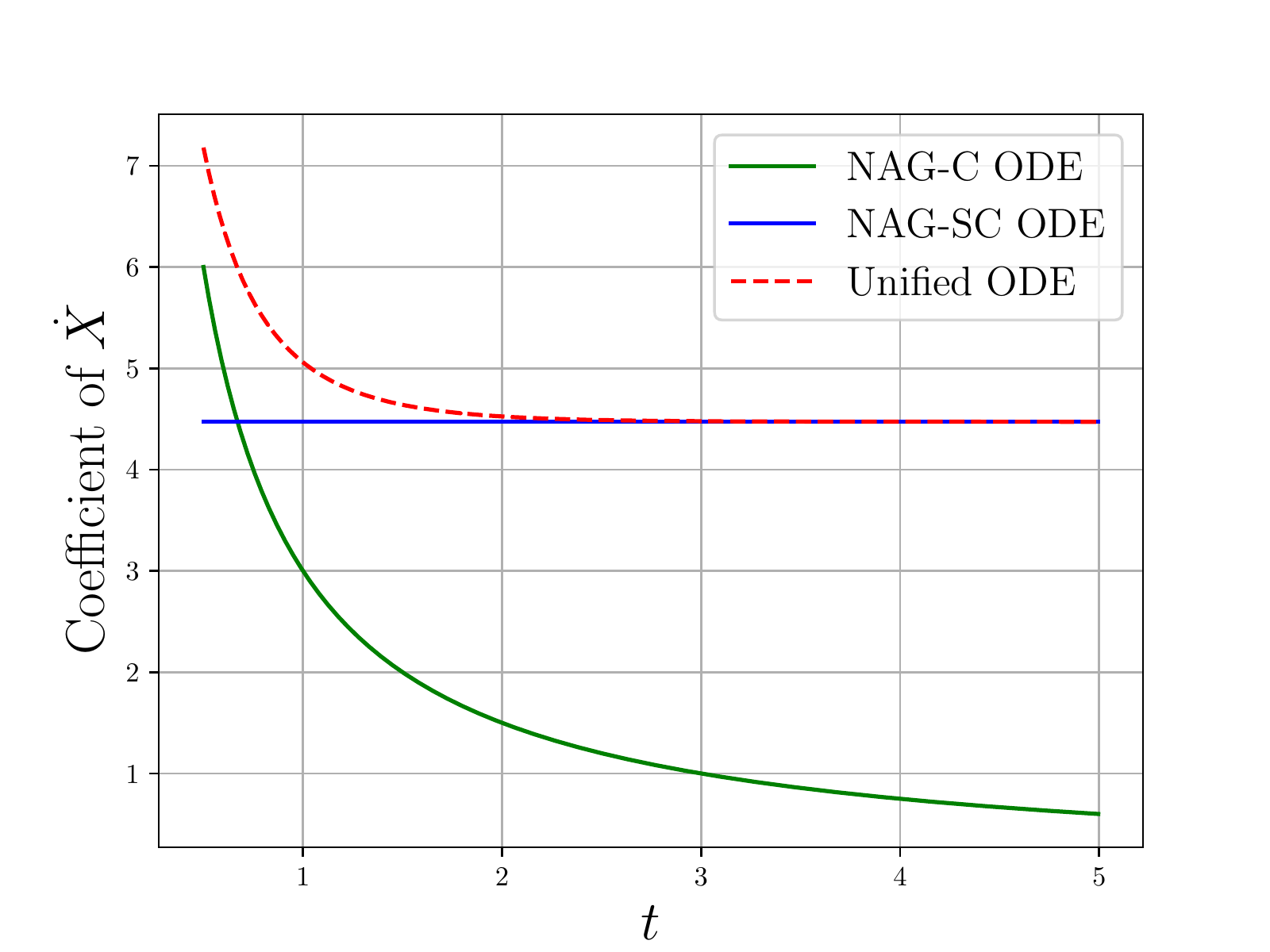}\label{fig:ode_coeff3}}
	\caption{Plots for the coefficient of $\dot{X}$, which can be interpreted as a measure of friction.}\label{fig:ode_coeff}
\end{figure}
\paragraph{Damping system interpretation.}
As mentioned in \citep{su2014}, the second-order ODE \eqref{eq:u_ode_oneline} can be viewed as a damping system, and the coefficient of $\dot{X}$ can be viewed as a measure of friction. Because the coefficient of $\dot{X}$ in \nagscsp ODE \eqref{eq:sc_ode} is $2\sqrt{\mu}$, \nagscsp ODE behaves like an \emph{underdamped system} when $\mu$ is small. Thus, the flow generated by \nagscsp ODE may present excessive oscillatory behaviors (see Figure~\ref{fig:experiments2}). In the unified NAG ODE \eqref{eq:u_ode_oneline}, the coefficient of $\dot{X}$ is large when $t$ is small  and converges to $2\sqrt{\mu}$ as $t\rightarrow\infty$ (see Figure~\ref{fig:ode_coeff}). Thus, the unified NAG ODE behaves like an \emph{overdamped system} (which displays less severe oscillations) when $t$ is small, regardless of the value of $\mu$.

\paragraph{Convergence analysis.}
For the unified NAG system,  the Lyapunov function \eqref{eq:lyapunov_lagrangian} can be written as
\begin{equation}
	\label{eq:lyapunov_continuous}
	V(X,Z,t)=\frac{1}{2}\cosh^{2}\left(\frac{\sqrt{\mu}}{2}t\right)\left\Vert Z-x^{*}\right\Vert ^{2}+\frac{t^{2}}{4}\sinhc^{2}\left(\frac{\sqrt{\mu}}{2}t\right)\left(f(X)-f\left(x^{*}\right)\right).
\end{equation}
Furthermore, we can rewrite Theorem~\ref{thm:mainthm_lagrangian} and Corollary~\ref{cor:maincor_Lagrangian} for this ODE model as follows:
\begin{theorem}
	\label{thm:mainthm_continuous}
	For the solution $(X,Z)$ to the unified NAG system \eqref{eq:u_ode}, the continuous-time energy functional
	\begin{align*}
		\mathcal{E}(t) & =V(X(t),Z(t),t)\\
		& =\frac{1}{2}\cosh^{2}\left(\frac{\sqrt{\mu}}{2}t\right)\left\Vert Z(t)-x^{*}\right\Vert ^{2}+\frac{t^{2}}{4}\sinhc^{2}\left(\frac{\sqrt{\mu}}{2}t\right)\left(f(X(t))-f\left(x^{*}\right)\right)
	\end{align*}
	is monotonically decreasing on $[\tst,\infty)$.
\end{theorem}
\begin{corollary}
	\label{cor:maincor_continuous}
	The solution $(X,Z)$ to the unified NAG  system \eqref{eq:u_ode} satisfies the inequality
	\begin{equation}
		\label{eq:guarantee_continuous}
		f(X(t))-f\left(x^{*}\right)\leq\frac{2}{t^{2}}\cschc^{2}\left(\frac{\sqrt{\mu}}{2}t\right)\left\Vert x_{0}-x^{*}\right\Vert ^{2}
	\end{equation}
	for all $t>\tst$.
\end{corollary}

Since $\cschc^2$ is decreasing on $[0,\infty)$, Corollary~\ref{cor:maincor_continuous} implies that the unified NAG ODE \eqref{eq:u_ode_oneline} achieves an $O(\|x_0-x^*\|^2/t^2)$ convergence rate regardless of the value of $\mu \geq 0$. When $\mu>0$, since $\frac{1}{t^{2}}\cschc^{2}\left(\frac{\sqrt{\mu}}{2}t\right)\sim\mu e^{-\sqrt{\mu}t}$ as $t\rightarrow\infty$, the unified NAG ODE achieves an $O(e^{-\sqrt{\mu}t}\|x_0-x^*\|^2)$ convergence rate. Combining these bounds, we conclude that the unified NAG ODE achieves an
\[
O\left(\min\left\{ 1/t^{2},e^{-\sqrt{\mu t}}\right\} \left\Vert x_{0}-x^{*}\right\Vert ^{2}\right)
\]
convergence rate.

\paragraph{Advantages of the unified NAG ODE compared to \nagscsp ODE.}
We now remark that our novel ODE model resolves the three drawbacks of \nagscsp ODE \eqref{eq:sc_ode} discussed in Section~\ref{sec:inconsistency}.
\begin{itemize}
	\item While the solution to \nagscsp ODE may not converge to the minimizer of $f$ when $\mu=0$, the solution to the unified NAG ODE always converges to the minimizer regardless of the value of $\mu$. 
	\item While the convergence guarantee for \nagscsp ODE may be worse than that for \nagcsp ODE in early stages, the convergence guarantee \eqref{eq:guarantee_continuous} for the unified NAG ODE is always better than that for \nagcsp ODE because $\cschc^2$ is decreasing on $[0,\infty)$ and the rate \eqref{eq:guarantee_continuous} recovers the exact convergence guarantee of \nagcsp ODE when $\mu=0$.
	\item While the convergence rate of \nagscsp ODE involves both the initial squared distance $\|x_0-x^*\|^2$ and the initial function value accuracy $f(x_0)-f(x^*)$, the convergence rate of the unified NAG ODE involves only the initial squared distance $\|x_0-x^*\|^2$.
\end{itemize}

\paragraph{Recovering \nagcsp ODE and \nagscsp ODE.}
We now discuss how \nagcsp ODE \eqref{eq:c_ode}, \nagscsp ODE \eqref{eq:sc_ode}, and their convergence analyses can be recovered from the proposed unified NAG ODE \eqref{eq:u_ode_oneline}. When $\mu=0$, it is easy to check that the unified ODE  recovers \nagcsp ODE and that the Lyapunov function \eqref{eq:lyapunov_continuous} recovers \eqref{eq:c_ode_lyapunov} for \nagcsp ODE. In the unified NAG ODE, because the coefficient of $\dot{X}$ converges to $2\sqrt{\mu}$ as $t\rightarrow\infty$ (see Figure~\ref{fig:ode_coeff}), \nagscsp ODE is the \emph{asymptotic version} of the unified NAG ODE. In Appendix~\ref{app:recover_lagrangian}, we show that the Lyapunov analysis for \nagscsp ODE can be recovered from Theorem~\ref{thm:mainthm_continuous} by taking the limit $t\rightarrow\infty$ of some inequalities.

\subsection{Proposed family of algorithms: Unified NAG family}
\label{sec:unified_nag}

Given the algorithmic stepsize $s$ and a strictly increasing sequence $(\bft_{k})_{k=0}^{\infty}$ (depending on $s$) in $[0,\infty)$ satisfying $\lim_{s\rightarrow0}\bft_0=0$, we consider the three-sequence scheme \eqref{eq:nag} with the algorithmic paramameters\footnote{We can constructively choose these sequences: First, we observe the relationship $\delta_k=\sqrt{s}\delta(\bft_{k+1})$, where $\bft_k=k\sqrt{s}$, between the algorithmic parameter $\delta_k=\frac{s(k+1)}{2}$ of \nagcsp and the coefficient $\delta(t)=\frac{t}{2}$ of \nagcsp system. Inspired by this relationship, for our algorithm, we define the sequence $\delta_k$ as $\delta_k=\sqrt{s}\delta(\bft_{k+1})$, where $\delta(t)=\frac{t}{2}\tanhc\left(\frac{\sqrt{\mu}}{2}t\right)$, and then set the sequence $\tau_k$ so that the collinearity condition \eqref{eq:collinear} holds.}
\begin{equation}
	\label{eq:tau_delta}
	\begin{aligned}
		\tau_{k} & =\frac{\frac{2\sqrt{s}}{\bft_{k+1}}\cothc\left(\frac{\sqrt{\mu}}{2}\bft_{k+1}\right)-\mu s}{1-\mu s}\\
		\delta_{k} & =\frac{\sqrt{s}\bft_{k+1}}{2}\tanhc\left(\frac{\sqrt{\mu}}{2}\bft_{k+1}\right),
	\end{aligned}
\end{equation}
that is, we consider the following \emph{unified NAG family}:
\begin{equation}
	\label{eq:nag-u}
	\begin{aligned}
		y_{k} & =x_{k}+\frac{\frac{2\sqrt{s}}{\bft_{k+1}}\cothc\left(\frac{\sqrt{\mu}}{2}\bft_{k+1}\right)-\mu s}{1-\mu s}\left(z_k-x_k\right)\\
		x_{k+1} & =y_{k}-s\nabla f\left(y_{k}\right)\\
		z_{k+1} & =z_{k}+\frac{\sqrt{s}\bft_{k+1}}{2}\tanhc\left(\frac{\sqrt{\mu}}{2}\bft_{k+1}\right)\left(\mu y_{k}-\mu z_{k}-\nabla f\left(y_{k}\right)\right).
	\end{aligned}
\end{equation}
Then, it is straightforward to check that the sequences $(\tau_k)$ and $(\delta_k)$ satisfy the collinearity condition \eqref{eq:collinear}. The following remark indicates that this algorithm can be regarded as a discretized version of the unified NAG system \eqref{eq:u_ode}.
\begin{remark}
	When the sequence $(\bft_{k})_{k=0}^{\infty}$ in $[0,\infty)$ satisfies the conditions \eqref{eq:tk_condition0} and \eqref{eq:tk_condition}, we have
	\begin{align*}
		\lim_{s\rightarrow0}\frac{\tau_{\bfk(t)}}{\sqrt{s}} & =\frac{2}{t}\cothc\left(\frac{\sqrt{\mu}}{2}t\right)\\
		\lim_{s\rightarrow0}\frac{\delta_{\bfk(t)}}{\sqrt{s}} & =\lim_{s\rightarrow0}\frac{t}{2}\tanhc\left(\frac{\sqrt{\mu}}{2}t\right)=\frac{t}{2}\tanhc\left(\frac{\sqrt{\mu}}{2}t\right)
	\end{align*}
	for all $t>\tst$, where $\bfk$ is the inverse function of the sequence $\bft$. Thus, the result in Section~\ref{sec:inconsistency} implies that the unified NAG family \eqref{eq:nag-u} converges to the unified NAG system \eqref{eq:u_ode} as $s\rightarrow 0$.
\end{remark}
When $\mu>0$ and $\lim_{k\rightarrow\infty}\bft_k=\infty$, we have $\lim_{k\rightarrow\infty}\tau_k=\frac{\sqrt{q}}{1+\sqrt{q}}$ and $\lim_{k\rightarrow\infty}\delta_k=\sqrt{\frac{s}{\mu}}$. Thus, \nagscsp \eqref{eq:nag-sc} is the \emph{asymptotic version} of the unified NAG family \eqref{eq:nag-u} in the sense that the coefficients of the unified NAG family converge to the coefficients of \nagsc. To obtain the convergence rate of the unified NAG family, we introduce the following assumptions on the sequence $(\bft_{k})$:\footnote{These assumption is purely inspired from the proof of Theorem~\ref{thm:mainthm_discrete}. Note that the assumptions \eqref{eq:tk_condition0} and \eqref{eq:tk_condition} are not required for the convergence analysis. Note that when $\mu=0$, under the identification $\theta_k=\frac{2\sqrt{s}}{t_k}$, the unified NAG family is equivalent to \citep[Algorithm~1]{tseng2008accelerated} and the condition \eqref{eq:tk_condition1} is equivalent to \citep[Equation~15]{tseng2008accelerated}.} 
\begin{equation}
	\label{eq:tk_condition2}
	\frac{2\sqrt{s}}{\bft_{k}}\cothc\left(\frac{\sqrt{\mu}}{2}\bft_{k}\right)\leq1\text{ for }k\geq 2
\end{equation}
and
\begin{equation}
	\label{eq:tk_condition1}
	\left(1-\frac{2\sqrt{s}}{\bft_{k+1}}\cothc\left(\frac{\sqrt{\mu}}{2}\bft_{k+1}\right)\right)\frac{\bft_{k+1}^{2}}{4}\sinhc^{2}\left(\frac{\sqrt{\mu}}{2}\bft_{k+1}\right)\leq\frac{\bft_{k}^{2}}{4}\sinhc^{2}\left(\frac{\sqrt{\mu}}{2}\bft_{k}\right)\text{ for }k\geq 0.
\end{equation}
The following results are the discrete-time analogs of Theorem~\ref{thm:mainthm_continuous} and Corollary~\ref{cor:maincor_continuous}.
\begin{theorem}
	\label{thm:mainthm_discrete}
	For the iterates of the unified NAG family \eqref{eq:nag-u} with $\left(\bft_{k}\right)_{k=0}^{\infty}$ satisfying the conditions \eqref{eq:tk_condition2} and \eqref{eq:tk_condition1}, the following discrete-time energy function is monotonically decreasing:
	\begin{equation}
		\label{eq:energy_unag_family}
		\mathcal{E}_{k}=V\left(x_{k},z_{k},\bft_{k}\right)=\frac{1}{2}\cosh^{2}\left(\frac{\sqrt{\mu}}{2}\bft_{k}\right)\left\Vert z_{k}-x^{*}\right\Vert ^{2}+\frac{\bft_{k}^{2}}{4}\sinhc^{2}\left(\frac{\sqrt{\mu}}{2}\bft_{k}\right)\left(f\left(x_{k}\right)-f\left(x^{*}\right)\right),
	\end{equation}
	where the Lyapunov function $V$ is defined in \eqref{eq:lyapunov_continuous}.
\end{theorem}
The proof of Theorem~\ref{thm:mainthm_discrete} can be found in Appendix~\ref{app:thm:mainthm_discrete}. Writing $\mathcal{E}_k\leq\mathcal{E}_0$ explicitly, we obtain the following result.
\begin{corollary}
	\label{cor:maincor_discrete}
	For the iterates of the unified NAG family \eqref{eq:nag-u} with $\left(\bft_{k}\right)_{k=0}^{\infty}$ satisfying the conditions \eqref{eq:tk_condition2} and \eqref{eq:tk_condition1}, the  following inequality holds for all $k\geq0$:
	\begin{multline}
		\label{eq:nag-u_convergence}
		f\left(x_{k}\right)-f\left(x^{*}\right)\leq			\frac{4}{\bft_{k}^{2}}\cschc^{2}\left(\frac{\sqrt{\mu}}{2}\bft_{k}\right)\\
		\times \left(\frac{1}{2}\cosh^{2}\left(\frac{\sqrt{\mu}}{2}\bft_{0}\right)\left\Vert x_{0}-x^{*}\right\Vert ^{2}+\frac{\bft_{0}^{2}}{4}\sinhc^{2}\left(\frac{\sqrt{\mu}}{2}\bft_{0}\right)\left(f\left(x_{0}\right)-f\left(x^{*}\right)\right)\right).
	\end{multline}
\end{corollary}

In the following subsections, we propose two concrete algorithms with specific choices of the sequence $(\bft_k)$. In \cref{sec:constant}, we propose the \emph{unified NAG}, a simple unified algorithm which continuously extend \nagcsp \eqref{eq:nag-c} to the strongly convex setting. In \cref{sec:variable}, we constructively recover the original NAG \eqref{eq:original_nag} and its convergence rate from the unified NAG family \eqref{eq:nag-u}.

\subsubsection{Constant timestep scheme: Unified NAG}
\label{sec:constant}

First, we set the constant timestep $\delta$ as
\begin{equation}
	\label{eq:delta_stepsize}
	\delta=\begin{cases}
		-\frac{\log\left(1-\sqrt{\mu s}\right)}{\sqrt{\mu}}, & \textrm{if }\mu>0\\
		\sqrt{s}, & \textrm{if }\mu=0,
	\end{cases}
\end{equation}
and then define the sequence $\bft_k=k\delta$, that is,
\begin{equation}
	\label{eq:constant_tk}
	\bft_{k}=\begin{cases}
		-\frac{\log(1-\sqrt{\mu s})}{\sqrt{\mu}}k, & \textrm{if }\mu>0\\
		\sqrt{s}k, & \textrm{if }\mu=0.
	\end{cases}
\end{equation}
Note that this choice is same as the previous choices of $\bft_k$ for \nagcsp and \nagscsp in Section~\ref{sec:lyapunov}. For this specific sequence $(\bft_k)$, the unified NAG family \eqref{eq:nag-u} can be written simply as
\begin{equation}
	\label{eq:unified_nag_specific}
	\begin{aligned}
		y_{k} & =x_{k}+\frac{\frac{2}{\iota(k+1)}\cothc\left(\frac{k+1}{2}\iota\sqrt{\mu s}\right)-\mu s}{1-\mu s}\left(z_{k}-x_{k}\right)\\
		x_{k+1} & =y_{k}-s\nabla f\left(y_{k}\right)\\
		z_{k+1} & =z_{k}+\frac{\iota s(k+1)}{2}\tanhc\left(\frac{k+1}{2}\iota\sqrt{\mu s}\right)\left(\mu y_{k}-\mu z_{k}-\nabla f\left(y_{k}\right)\right),
	\end{aligned}
\end{equation}
where $\iota=-\frac{\log(1-\sqrt{\mu s})}{\sqrt{\mu s}}$ for $\mu>0$ and $\iota=1$ for $\mu=0$. We refer to this algorithm as the \emph{unified NAG}. 

The sequence $(\bft_k)$ in \eqref{eq:constant_tk} can be shown to satisfy the conditions \eqref{eq:tk_condition2} and \eqref{eq:tk_condition1} (see \cref{app:constant}), and thus the convergence guarantee \eqref{eq:nag-u_convergence} holds for this specific algorithm. Also it is straightforward to check that the conditions \eqref{eq:tk_condition0} and \eqref{eq:tk_condition} hold, and thus the unified NAG \eqref{eq:unified_nag_specific} converges to the unified NAG system \eqref{eq:u_ode} as $s\to0$. Because $\cschc^2$ is decreasing on $[0,\infty)$ and $\delta\geq\sqrt{s}$, we have
\[
\frac{4}{\bft_{k}^{2}}\cschc^{2}\left(\frac{\sqrt{\mu}}{2}\bft_{k}\right)\leq\frac{4}{\bft_{k}^{2}}=\frac{4}{\delta^2k^2}\leq\frac{4}{sk^{2}}.
\]
This implies that the convergence guarantee of the unified NAG is always better than that of \nagcsp and that the unified NAG achieves an $O(\|x_0-x^*\|^2/k^2)$ convergence rate, regardless of the value of $\mu$. When $\mu>0$, since
\[
\frac{4}{\bft_{k}^{2}}\cschc^{2}\left(\frac{\sqrt{\mu}}{2}\bft_{k}\right)\sim4\mu e^{-\sqrt{\mu}\bft_{k}}=4\mu\left(1-\sqrt{\mu s}\right)^{k}\text{ as }k\rightarrow\infty,
\]
the unified NAG achieves an $O((1-\sqrt{\mu s})^k\|x_0-x^*\|^2)$ convergence rate. Combining these two guarantees, we conclude that the unified NAG achieves an $$O\left(\min\left\{ 1/k^{2},\left(1-\sqrt{\mu s}\right)^{k}\right\} \left\Vert x_{0}-x^{*}\right\Vert ^{2}\right)$$ convergence rate.

\paragraph{Advantages of the unified NAG compared to \nagsc.}
We now highlight that the unified NAG resolves the three drawbacks of \nagscsp \eqref{eq:nag-c} discussed in Section~\ref{sec:introduction}.
\begin{itemize}
	\item While \nagscsp cannot handle the non-strongly convex case, the unified NAG can handle the case $\mu=0$. Moreover, when $\mu=0$, the unified NAG and its convergence rate \eqref{eq:nag-u_convergence} recover \nagcsp and its convergence rate.
	\item While the convergence guarantee for \nagscsp may be worse than that for \nagcsp in early stages, the convergence guarantee for the unified NAG is always better than that for \nagc.
	\item While the convergence rate of \nagscsp involves both the initial squared distance $\|x_0-x^*\|^2$ and the initial function value accuracy $f(x_0)-f(x^*)$, the convergence rate of the unified NAG involves only the initial squared distance $\|x_0-x^*\|^2$.
\end{itemize}

\subsubsection{Adaptive timestep scheme: Recovering the original NAG}
\label{sec:variable}

The constant timestep scheme (unified NAG) in the previous section can be improved in terms of the convergence rate by defining the sequence $(\bft_{k})_{k=0}^{\infty}$ more aggressively as
\begin{equation}
	\label{eq:adaptive_timestep}
	\bft_{k+1}=\begin{cases}
		\textrm{Given constant }\bft_0>0 \textrm{ (possibly depending on } s), & k+1=0\\
		\textrm{The largest real number satisfying \eqref{eq:tk_condition1}}, & k+1\geq1.
	\end{cases}
\end{equation}
Then, it is easy to check that the sequence $(\bft_{k})_{k=0}^{\infty}$ is well-defined and strictly increasing. We refer to the unified NAG family \eqref{eq:nag-u} with this sequence as the \emph{adaptive timestep scheme}.

Note that the conditions \eqref{eq:tk_condition2} and \eqref{eq:tk_condition1} hold by construction.\footnote{The first condition follows from the facts that \eqref{eq:tk_condition2} holds for the sequence $\bft_k=k\delta$ (see Section~\ref{sec:constant}) and we have $\bft_k>2\delta$ for $k\geq2$ for the sequence \eqref{eq:adaptive_timestep}.} Therefore, the convergence guarantee \eqref{eq:nag-u_convergence} holds for the adaptive timestep scheme. In \cref{app:variable}, we show that if $\bft_0\rightarrow0$ as $s\rightarrow0$, then the conditions \eqref{eq:tk_condition0} and \eqref{eq:tk_condition} hold, and thus the adaptive timestep scheme converges to the unified NAG system \eqref{eq:u_ode} as $s\to0$.
By construction, we have $\bft_{k+1}-\bft_{k}>\delta$, where $\delta$ is defined in \eqref{eq:delta_stepsize}, which implies that  $\bft_{k}>\bft_{0}+k\delta$ for all $k\geq0$. Thus, the adaptive timestep scheme has a (slightly) better convergence rate than the unified NAG. Surprisingly, our new algorithm, which is purely obtained from the unified Lagrangian framework, is equivalent to the original Nesterov's method \eqref{eq:original_nag}.  
\begin{proposition}
	\label{prop:equivalence}
	The  adaptive timestep scheme is equivalent to the original NAG \eqref{eq:original_nag} with $\gamma_{0}=\frac{4}{\bft_{0}^{2}}\cothc^{2}\left(\frac{\sqrt{\mu}}{2}\bft_{0}\right)>\mu$. Moreover, the sequence $\gamma_k$ and $\alpha_k$ in the original NAG can be written as $\gamma_{k}=\frac{4}{\bft_{k}^{2}}\cothc^{2}\left(\frac{\sqrt{\mu}}{2}\bft_{k}\right)$ and $\alpha_{k}=\frac{2\sqrt{s}}{\bft_{k+1}}\cothc\left(\frac{\sqrt{\mu}}{2}\bft_{k+1}\right)$. Conversely, the original NAG \eqref{eq:original_nag} with $\gamma_0>\mu$ is equivalent to the adaptive timestep scheme, where $\bft_{0}$ satisfies $\gamma_{0}=\frac{4}{\bft_{0}^{2}}\cothc^{2}\left(\frac{\sqrt{\mu}}{2}\bft_{0}\right)$.
\end{proposition}

The proof of Proposition~\ref{prop:equivalence} can be found in \cref{app:prop:equivalence}. The following remark shows that under the identification in \cref{prop:equivalence}, the convergence rate \eqref{eq:nag-u_convergence} of the adaptive timestep scheme is equivalent to the convergence rate \eqref{eq:original_nag_convergence} of the original NAG obtained by \citet{nesterov2018lectures}.

\begin{remark}
	By Corollary~\ref{cor:maincor_discrete}, the iterates of the adaptive timestep scheme satisfy 
	{\begin{equation}
			\label{eq:recovering_original_convergnece}
			\begin{aligned}
				& f\left(x_{k}\right)-f\left(x^{*}\right)\\
				& \leq\frac{4}{\bft_{k}^{2}}\cschc^{2}\left(\frac{\sqrt{\mu}}{2}\bft_{k}\right)\\
				& \quad \times \left(\frac{1}{2}\cosh^{2}\left(\frac{\sqrt{\mu}}{2}\bft_{0}\right)\left\Vert x_{0}-x^{*}\right\Vert ^{2}+\frac{\bft_{0}^{2}}{4}\sinhc^{2}\left(\frac{\sqrt{\mu}}{2}\bft_{0}\right)\left(f\left(x_{0}\right)-f\left(x^{*}\right)\right)\right)\\
				& =\frac{4}{\bft_{k}^{2}}\cschc^{2}\left(\frac{\sqrt{\mu}}{2}\bft_{k}\right)\frac{\bft_{0}^{2}}{4}\sinhc^{2}\left(\frac{\sqrt{\mu}}{2}\bft_{0}\right)\\
				& \quad \times \left(\frac{2}{\bft_{0}^{2}}\cothc^{2}\left(\frac{\sqrt{\mu}}{2}\bft_{0}\right)\left\Vert x_{0}-x^{*}\right\Vert ^{2}+\left(f\left(x_{0}\right)-f\left(x^{*}\right)\right)\right)\\
				& =\prod_{i=0}^{k-1}\left(1-\frac{2\sqrt{s}}{\bft_{i+1}}\cothc\left(\frac{\sqrt{\mu}}{2}\bft_{i+1}\right)\right)   \left(\frac{2}{\bft_{0}^{2}}\cothc^{2}\left(\frac{\sqrt{\mu}}{2}\bft_{0}\right)\left\Vert x_{0}-x^{*}\right\Vert ^{2}+\left(f\left(x_{0}\right)-f\left(x^{*}\right)\right)\right),
			\end{aligned}
	\end{equation}}
	where the last equality follows from our updating rule \eqref{eq:adaptive_timestep} of the sequence $(\bft_k)$. Therefore, we recover the convergence rate \eqref{eq:original_nag_convergence} of the original NAG with $\gamma_{k}=\frac{4}{\bft_{k}^{2}}\cothc^{2}\left(\frac{\sqrt{\mu}}{2}\bft_{k}\right)$ and $\alpha_{k}=\frac{2\sqrt{s}}{\bft_{k+1}}\cothc\left(\frac{\sqrt{\mu}}{2}\bft_{k+1}\right)$. 
\end{remark}

\section{Extension to Higher-order Non-Euclidean Setting}
\label{sec:higher-order}

Based on the first Bregman Lagrangian \eqref{eq:first_lagrangian} and the prior work \citep{baes2009estimate}, \citet{wibisono2016} proposed the \emph{accelerated tensor flow} and \emph{accelerated tensor method} for convex objective functions to achieve a polynomial $O(1/t^p)$ or $O(1/k^p)$ convergence rate. They also tried to design accelerated tensor methods for uniformly convex objective functions to achieve an \emph{exponential convergence rate}. They were able to obtain an exponential convergence rate for continuous-time flows obtained from the first Bregman Lagrangian, but  a rate-matching discretization was not identified. Instead, they showed that the accelerated tensor method (convex case) with a restart scheme achieves an exponential convergence rate for uniformly convex objective functions. However, as they admitted, understanding the connection between the discrete-time algorithm and the continuous-time flow is unclear and remains as an open problem.

In this section, using the unified Bregman Lagrangian~\eqref{eq:unified_lagrangian}, we continuously extend to the accelerated tensor flow and the accelerated tensor method in \citep{wibisono2016} to the strongly convex case. Our novel dynamics and algorithm achieve exponential convergence rates without using a restarting technique.

We make the following assumptions throughout this section:
\begin{itemize}
	\item The distance-generating function $h$ is $1$-uniformly convex \eqref{eq:higher_order_uniform_convexity} of order $p\geq2$.
	\item The objective function $f$ is $\mu$-uniformly (possibly with $\mu=0$) convex \eqref{eq:uniform_convexity} with respect to the distance-generating function $h$. 
	\item The objective function $f$ is $\frac{(p-1)!}{s}$-smooth of order $p-1$ \eqref{eq:higher_order_smoothness}, where $s$ is the algorithmic stepsize. 
\end{itemize}
These assumptions are standard in the literature of higher-order optimization~\citep[see][]{nesterov2008accelerating,baes2009estimate,wibisono2016,gasnikov2019optimal,wilson2021lyapunov}. In particular, when $p=2$ and $h(x)=\Vert x\Vert^{2}$, these assumptions recover the standard smooth strongly convex setting in Section~\ref{sec:convexity}.

Following \citep{wibisono2016}, we define the \emph{tensor update operator} $G_{p,s,N}:\mathbb{R}^n\rightarrow\mathbb{R}^n$ as
\begin{equation}
	\label{eq:higher_order_gradient_update}
	G_{p,s,N}(y)=\arg\min_{x}\left\{ f_{p-1}(x;y)+\frac{N}{ps}\left\Vert x-y\right\Vert ^{p}\right\} ,
\end{equation}
where the function $x\mapsto f_{p-1}(x,y)=\sum_{i=0}^{p-1}\frac{1}{i!}\nabla^i f(y)(x-y)^i$ is the $(p-1)$-st order Taylor approximation of the objective function $f$ at $y\in\mathbb{R}^n$.
\citet[Lemma~2.2]{wibisono2016} showed that one can choose $N>0$ so that there exists a constant $M>0$ for which the inequality
\begin{equation}
	\label{eq:M_ineq}
	\left\langle \nabla f\left(x\right),y-x\right\rangle \geq Ms^{\frac{1}{p-1}}\left\Vert \nabla f\left(x\right)\right\Vert ^{\frac{p}{p-1}}
\end{equation}
holds for $x=G_{p,s,N}(y)$. From now on, we denote the tensor update operator satisfying the inequality \eqref{eq:M_ineq} by $G_{p,M}$. As a special case, when $p=2$, the operator \eqref{eq:higher_order_gradient_update} with $N=1$ satisfies the inequality \eqref{eq:M_ineq} with $M=1/2$.\footnote{See, for example, the proof of Lemma~6 in the arXiv version of \citep{wilson2021lyapunov}: \texttt{arXiv:1611.02635v4}.}

\subsection{Proposed dynamics: Unified accelerated tensor flow}
\label{sec:unified_higher_flow}

We consider the unified Bregman Lagrangian flow \eqref{eq:unified_family} with the parameters
\begin{equation}
	\label{eq:alpha_beta_for_unified}
	\begin{aligned}
		\alpha(t) & =\log p-\log t+\log\left(\cothc_p\left(C^{1/p}\mu^{1/p}t\right)\right)\\
		\beta(t) & =p\log t+\log C+p\log\left(\sinhc_{p}\left(C^{1/p}\mu^{1/p}t\right)\right)
	\end{aligned}
\end{equation}
and the initial conditions $X(\tst)=Z(\tst)=x_0$, where $C>0$ is a constant. It is straightforward to check that the ideal scaling condition \eqref{eq:ideal_scaling_b} holds. This dynamical system can be written as
\begin{equation}
	\label{eq:hyperbolic_flow}
	\begin{aligned}
		\dot{X} & =\frac{p}{t}\cothc_p\left(C^{1/p}\mu^{1/p}t\right)(Z-X)\\
		\frac{d}{dt}\nabla h(Z) & =Cpt^{p-1}\tanhc_{p}^{p-1}\left(C^{1/p}\mu^{1/p}t\right)\left(\mu\nabla h(X)-\mu\nabla h(Z)-\nabla f(X)\right).
	\end{aligned}
\end{equation}
From now on, we refer to this system of ODEs  as the \emph{unified accelerated tensor flow}. Using the existence and uniqueness of solution to the unified NAG system (Theorem~\ref{thm:existence_unique_unified}) and the time-dilation property (Theorem~\ref{thm:time_dilation}), we can prove the following theorem (see Appendix~\ref{app:existence_uniqueness_hyperbolic}).
\begin{theorem}
	\label{thm:existence_unique_tensor}
	The unified accelerated tensor flow \eqref{eq:u_ode} has a unique solution $(X,Z)$ in $C^1([0,\infty),\mathbb{R}^n\times\mathbb{R}^n)$.
\end{theorem}

For this dynamical system, the Lyapunov function \eqref{eq:lyapunov_lagrangian} can be expressed as
\begin{equation}
	\label{eq:lyapunov_higher}
	V(X,Z,t)=\cosh_p^{p}\left(C^{1/p}\mu^{1/p}t\right)D_{h}\left(x^{*},Z\right)+Ct^{p}\sinhc_{p}^{p}\left(C^{1/p}\mu^{1/p}t\right)\left(f(X)-f\left(x^{*}\right)\right).
\end{equation}
We can rewrite Theorem~\ref{thm:mainthm_lagrangian} and Corollary~\ref{cor:maincor_Lagrangian} for the unified accelerated tensor flow \eqref{eq:hyperbolic_flow} as follows:
\begin{theorem}
	\label{thm:mainthm_continuous_higher}
	For the solution $(X,Z)$ to the unified accelerated tensor flow \eqref{eq:hyperbolic_flow}, the continuous-time energy function
	\begin{align*}
		\mathcal{E}(t) & =V(X(t),Z(t),t)\\
		& =\cosh_{p}^{p}\left(C^{1/p}\mu^{1/p}t\right)D_{h}\left(x^{*},Z(t)\right)+Ct^{p}\sinhc_{p}^{p}\left(C^{1/p}\mu^{1/p}t\right)\left(f(X(t))-f\left(x^{*}\right)\right)
	\end{align*}
	is monotonically decreasing on $[\tst,\infty)$.
\end{theorem}
\begin{corollary}
	\label{cor:maincor_continuous_higher}
	The solution $(X,Z)$ to the unified accelerated tensor flow \eqref{eq:hyperbolic_flow} satisfies the inequality
	\begin{equation}
		\label{eq:guarantee_continuous_higher}
		f(X(t))-f\left(x^{*}\right)\leq\frac{1}{Ct^{p}\sinhc_{p}^{p}\left(C^{1/p}\mu^{1/p}t\right)}D_{h}\left(x^{*},x_0\right)
	\end{equation}
	for all $t>\tst$.
\end{corollary}

Since $\sinhc_p(0)=1$ and $\sinhc_p$ is increasing on $[0,\infty)$ (see Appendix~\ref{app:sinhc_is_increasing}), Corollary~\ref{cor:maincor_continuous_higher} implies that the unified accelerated tensor flow \eqref{eq:hyperbolic_flow} achieves an $O\left(D_{h}\left(x^{*},x_0\right)/t^p\right)$ convergence rate regardless of the value of $\mu \geq 0$. On the other hand, when $\mu>0$, it follows from Proposition~\ref{prop:sinhp_grows_exponentially} that
\[
\frac{1}{Ct^{p}\sinhc_{p}^{p}\left(C^{1/p}\mu^{1/p}t\right)}=O\left(e^{-pC^{1/p}\mu^{1/p}t}\right) \quad \mbox{as $t\rightarrow\infty$}.
\]
 Therefore, the unified accelerated tensor flow achieves an $O(e^{-pC^{1/p}\mu^{1/p}t}D_{h}\left(x^{*},x_0\right))$ convergence rate. Combining these bounds, we conclude that the unified accelerated tensor flow achieves an $$O\left(\min\left\{ 1/t^{p},e^{-pC^{1/p}\mu^{1/p}t}\right\} D_{h}\left(x^{*},x_{0}\right)\right)$$ convergence rate.

\subsection{Proposed algorithm: Unified accelerated tensor method}
\label{sec:unified_higher}

As a discretization scheme for the unified accelerated tensor flow \eqref{eq:hyperbolic_flow}, we propose the following \emph{unified accelerated tensor method family}:
\begin{subequations}
	\label{eq:hyperbolic_scheme}
	\begin{align}
		\label{eq:hyperbolic_scheme1}
		A_k & = C\bft_{k}^{p}\sinhc_{p}^{p}\left(C^{1/p}\mu^{1/p}\bft_{k}\right)\\
		\label{eq:hyperbolic_scheme2}
		y_{k} & =x_{k}+\frac{A_{k+1}-A_{k}}{A_{k+1}}\left(z_{k}-x_{k}\right)\\
		\label{eq:hyperbolic_scheme3}
		x_{k+1} & =G_{p,M}\left(y_{k}\right)\\
		\label{eq:hyperbolic_scheme4}
		z_{k+1} & =\arg\min_{z}\left\{ \frac{A_{k+1}-A_{k}}{1+\mu A_{k}}\left(\left\langle \nabla f\left(x_{k+1}\right),z\right\rangle +\mu D_{h}\left(z,x_{k+1}\right)\right)+D_{h}\left(z,z_{k}\right)\right\} ,
	\end{align}
\end{subequations}
where $(\bft_k)$ is a strictly increasing sequence (depending on the algorithmic stepsize $s$) in $[0,\infty)$ and $G_{p,M}$ is the tensor update operator satisfying \eqref{eq:M_ineq}. Because the algorithm \eqref{eq:hyperbolic_scheme} is continuous in the strong convexity parameter $\mu$, it handles the convex case and the strongly convex case in a unified way. By the first-order optimality condition, {the step \eqref{eq:hyperbolic_scheme4}} is equivalent to
\begin{equation}
	\label{eq:our_zk}
	\nabla h\left(z_{k+1}\right)-\nabla h\left(z_{k}\right)=\frac{A_{k+1}-A_k}{1+\mu A_{k}}\left(\mu\nabla h\left(x_{k+1}\right)-\mu\nabla h\left(z_{k+1}\right)-\nabla f\left(x_{k+1}\right)\right).
\end{equation}
Although   the scheme~\eqref{eq:hyperbolic_scheme} cannot be written in the three-sequence form \eqref{eq:nag}, we observe that the step \eqref{eq:hyperbolic_scheme2} plays a role of \eqref{eq:nag1} (updating $y_k$ as a convex combination of $x_k$ and $z_k$), the step \eqref{eq:hyperbolic_scheme4} plays a role similar to \eqref{eq:nag3} (updating $z_k$ by gradient/mirror step), and that the tensor update step \eqref{eq:hyperbolic_scheme3} corresponds to the gradient update step \eqref{eq:nag2}.

\paragraph{Limiting ODE.}
In Appendix~\ref{app:limiting_ode_higher}, we show that if
\begin{equation}
	\label{eq:tk_condition0_higher}
	\lim_{s\rightarrow0}\bft_{0}=\tst
\end{equation}
and the timesteps are asymptotically equivalent to $s^{1/p}$ as $s\rightarrow0$ in the sense that
\begin{equation}
	\label{eq:tk_condition_higher}
	\lim_{s\rightarrow0}\frac{\bft_{\bfk(t)+1}-t}{s^{1/p}}=1\text{ for all }t\in\left(\tst,\infty\right),
\end{equation}
where $\bfk$ is the inverse of $\bft$, then the unified accelerated tensor method family \eqref{eq:hyperbolic_scheme} converges to the unified accelerated tensor flow \eqref{eq:hyperbolic_flow} when letting  $x_k=X(\bft_k)$ and $z_k=Z(\bft_k)$.

\paragraph{Convergence analysis.}
To prove the convergence rate, we introduce the following assumption on the sequence $(\bft_k)$ (note that $A_k$ is uniquely determined by $\bft_k$ and vice versa):
\begin{equation}
	\label{eq:our_Ak}
	\left(A_{k+1}-A_{k}\right)^{p}-Cp^{p}sA_{k+1}^{p-1}\left(1+\mu A_{k}\right)\leq0 \textrm{ with }C=\frac{1}{p}\left(\frac{M}{p-1}\right)^{p-1},
\end{equation}
where $M$ is the constant involved in \eqref{eq:M_ineq}. The following results are the discrete-time analogs of Theorem~\ref{thm:mainthm_continuous_higher} and Corollary~\ref{cor:maincor_continuous_higher}.
\begin{theorem}
	\label{thm:mainthm_discrete_higher}
	For the iterates of the unified accelerated tensor method family \eqref{eq:hyperbolic_scheme} with $(\bft_k)$ satisfying the condition \eqref{eq:our_Ak}, the discrete-time energy function
	\begin{equation}
		\label{eq:discrete_energy_higher}
		\mathcal{E}_{k}=\left(1+\mu A_{k}\right)D_{h}\left(x^{*},z_{k}\right)+A_{k}\left(f\left(x_{k}\right)-f\left(x^{*}\right)\right)
	\end{equation}
	is monotonically decreasing.
\end{theorem}
The proof of Theorem~\ref{thm:mainthm_discrete_higher} can be found in Appendix~\ref{app:thm:mainthm_discrete_higher}. Writing $\mathcal{E}_k\leq\mathcal{E}_0$ explicitly, we obtain the following result.
\begin{corollary}
	\label{cor:maincor_discrete_higher}
	For the iterates of the unified accelerated tensor method family \eqref{eq:hyperbolic_scheme} with $(\bft_k)$ satisfying the condition \eqref{eq:our_Ak}, the  following inequality holds for all $k\geq0$:
	\begin{equation}
		\label{eq:nag-u_convergence_higher}
		f\left(x_{k}\right)-f\left(x^{*}\right)\leq\frac{1}{A_{k}}\left(\left(1+\mu A_{0}\right)D_{h}\left(x^{*},x_{0}\right)+A_{0}\left(f\left(x_{0}\right)-f\left(x^{*}\right)\right)\right).
	\end{equation}
\end{corollary}

\paragraph{Specific algorithm: Unified accelerated tensor method.}

We now consider the following specific choice of sequence $(\bft_k)$:
\begin{equation}
	\label{eq:our_specific_Ak}
	\bft_{k+1}=\begin{cases}
		0, & k+1=0\\
		\textrm{The largest real number satisfying \eqref{eq:our_Ak}}, & k+1\geq1.
	\end{cases}
\end{equation}
Then, the condition \eqref{eq:our_Ak} clearly holds, and thus the convergence results hold. In addition, we can show that this sequence satisfies the conditions \eqref{eq:tk_condition0_higher} and \eqref{eq:tk_condition_higher} (see Appendix~\ref{app:limiting_ode_higher}). Hence, the algorithm converges to the unified accelerated tensor flow \eqref{eq:hyperbolic_flow} as $s\to0$. Furthermore, we can show that the inequalities
\[
A_{k} \geq O\left(k^{p}\right),\quad A_{k} \geq O\left(\left(1+C^{1/p}p\mu^{1/p}s^{1/p}\right)^{k}\right)
\]
hold (see Appendix~\ref{app:lower_bound_Ak}). Therefore, Corollary~\ref{cor:maincor_discrete_higher} implies the following convergence rate:
\[
f\left(x_{k}\right)-f\left(x^{*}\right)\leq O\left(\min\left\{ 1/k^{p},\left(1+C^{1/p}p\mu^{1/p}s^{1/p}\right)^{-k}\right\} \right).
\]

\subsection{Recovering the non-strongly convex case}
When $\mu=0$, the  system  of ODEs~\eqref{eq:hyperbolic_flow} recovers the following \emph{accelerated tensor flow (convex case)} given in \citep{wibisono2016}:\footnote{This flow can be obtained by putting $\alpha(t)=\log p-\log t$ and $\beta(t)=p\log t+\log C$ (Equation~\ref{eq:alpha_beta_for_unified} with $\mu=0$) in the first Bregman Lagrangian flow \eqref{eq:first_family}.}
\begin{equation}
	\label{eq:polynomial_flow}
	\begin{aligned}
		\dot{X} & =\frac{p}{t}(Z-X)\\
		\frac{d}{dt}\nabla h(Z) & =-Cpt^{p-1}\nabla f(X).
	\end{aligned}
\end{equation}
Moreover, the unified accelerated tensor method family \eqref{eq:hyperbolic_scheme} becomes the following family:
\begin{equation}
	\label{eq:polynomial_scheme}
	\begin{aligned}
		A_k & =C\bft_{k}^{p}\\
		y_{k} & =x_{k}+\frac{A_{k+1}-A_{k}}{A_{k+1}}\left(z_{k}-x_{k}\right)\\
		x_{k+1} & =G_{p,M}\left(y_{k}\right)\\
		z_{k+1} & =\arg\min_{z}\left\{ \left(A_{k+1}-A_{k}\right)\left\langle \nabla f\left(x_{k+1}\right),z\right\rangle +D_{h}\left(z,z_{k}\right)\right\}.
	\end{aligned}
\end{equation}
This recovers the \emph{accelerated tensor method (convex case)} in \citep{wibisono2016} if the sequence $(\bft_k)$ is chosen as
\begin{equation}
	\label{eq:wibisono_Ak}
	\bft_{k}=s^{1/p}k^{1/p}(k+1)^{1/p}\cdots(k+p-1)^{1/p},
\end{equation}
for which the inequality \eqref{eq:our_Ak} holds with $\mu=0$. 

\section{Further Exploration: ODE Model for Minimizing Gradient Norms of Strongly Convex Functions}
\label{sec:minimizing_gradient}

So far, we have focused on ODEs and algorithms that  achieve a fast convergence rate for the accuracy of objective function values $f(X(t))-f(x^*)$ or $f(x_k)-f(x^*)$. Typically, the goal of numerically solving a convex optimization problem is to reduce the deviation from the minimum value. 
Alternatively, the gradient norm $\|\nabla f(x_k)\|^2$ can be used as a performance measure. 
This criterion is often reasonable for both theoretical and practical purposes \citep[see][]{nesterov2012make,diakonikolas2022potential}. Recently, \citet{kim2021optimizing} proposed OGM-G, which is a method that achieves the optimal convergence rate (up to a constant factor) for {minimizing} the gradient norm $\|\nabla f(x_N)\|^2$ of non-strongly convex functions. Recently, this method has attracted some attention: \citet{lee2021geometric} provided a Lyapunov argument for its convergence analysis. \citet{suh2022continuous} derived and analyzed the limiting ODE of OGM-G. However, most studies on OGM-G have focused only on the non-strongly convex case.

In this section, we propose a novel continuous-time dynamical system that reduces the squared gradient norm $\Vert\nabla f(X(T))\Vert^2$ of strongly convex objective functions $f$ with an
\[
O\left(\min\left\{ 1/T^{2},e^{-\sqrt{\mu}T}\right\} \left(f(x_{0})-f\left(x^*\right)+\frac{\mu}{2}\left\Vert x_{0}-X(T)\right\Vert ^{2}\right)\right)
\]
convergence rate. Interestingly, the ODE model presented in this section and the unified NAG ODE \eqref{eq:u_ode_oneline} have an \emph{anti-transpose} relationship between the corresponding differential kernels.

\subsection{Motivation: Symmetric relationship between OGM ODE and OGM-G ODE}
For non-strongly convex objective functions, \citet{suh2022continuous} proposed OGM-G ODE, an ODE model whose solution $X:[0,T]\rightarrow\mathbb{R}^n$ reduces the squared gradient norm $\|\nabla f(X(T))\|^2$ with an $O((f(x_0)-f(x^*))/T^2)$ convergence rate. In this section, we investigate a symmetric relationship between OGM ODE (which we will discuss later) and OGM-G ODE. This relationship will give us a hint for designing our novel ODE model.

\paragraph{Anti-transpose relationship between OGM and OGM-G.}

We first review a symmetric relationship between OGM \citep{kim2016optimized}, an algorithm for reducing the function value accuracy $f(x_N)-f(x^*)$, and OGM-G \citep{kim2021optimizing}, an algorithm for reducing the squared gradient norm $\|\nabla f(x_N)\|^2$. Given the number $N$ of total iterations, define a sequence $(\theta_k)_{k=0}^N$ as
\begin{equation}
	\label{eq:theta_sequence}
	\theta_{k}=\begin{cases}
		1 & \textrm{if }k=0\\
		\frac{1+\sqrt{4\theta_{k-1}^{2}+1}}{2} & \textrm{if }1\leq k\leq N-1\\
		\frac{1+\sqrt{8\theta_{k-1}^{2}+1}}{2} & \textrm{if }k=N.
	\end{cases}
\end{equation}
Then, OGM is equivalent to the fixed-step first-order scheme \eqref{eq:discrete_FSFO} with the difference matrix $\mathbf{H}^{\mathrm{F}}$, and OGM-G is equivalent to the fixed-step first-order scheme \eqref{eq:discrete_FSFO} with the difference matrix $\mathbf{H}^{\mathrm{G}}$, where the entries of $\mathbf{H}^{\mathrm{F}}$ and $\mathbf{H}^{\mathrm{G}}$ are defined as
\begin{equation}
	\label{eq:HFHG}
	\begin{aligned}
		h_{ij}^{\mathrm{F}} & =\begin{cases}
			\frac{\theta_{i}-1}{\theta_{i+1}}h_{i-1,j} & \mbox{if }k=0,\ldots,i-2,\\
			\frac{\theta_{i}-1}{\theta_{i+1}}\left(h_{i-1,i-1}-1\right) & \mbox{if }k=i-1,\\
			1+\frac{2\theta_{i}-1}{\theta_{i+1}} & \mbox{if }k=i,
		\end{cases}\\
		h_{ij}^{\mathrm{G}} & =\begin{cases}
			\frac{\theta_{N-i-1}-1}{\theta_{N-i}}h_{i,j+1} & \mbox{if }k=0,\ldots,i-2,\\
			\frac{\theta_{N-i-1}-1}{\theta_{N-i}}\left(h_{i,i}-1\right) & \mbox{if }k=i-1,\\
			1+\frac{2\theta_{N-i-1}-1}{\theta_{N-i}} & \mbox{if }k=i.
		\end{cases}
	\end{aligned}
\end{equation}
\citet{kim2021optimizing} observed the following relationship between the difference kernels for OGM and OGM-G:
\begin{equation}
	\label{eq:symmetry_beta_gamma}
	h^{\mathrm{F}}_{ij}=h^{\mathrm{G}}_{N-1-j,N-1-i} \textrm{ for all }i\textrm{ and }j.
\end{equation}
When the condition \eqref{eq:symmetry_beta_gamma} holds, we say there is an \emph{anti-transpose} relationship between $\mathbf{H}^{\mathrm{F}}$ and $\mathbf{H}^{\mathrm{G}}$ because the matrix $\mathbf{H}^{\mathrm{F}}$ can be obtained by reflecting $\mathbf{H}^{\mathrm{G}}$ about its anti-diagonal and vice versa.

\paragraph{A (naive) symmetric relationship between OGM ODE and OGM-G ODE.}
Next, we look at the relationship between the limiting ODEs of OGM and OGM-G. When letting $T=N\sqrt{s}$ and $x_k=X(t_k)$, OGM converges to the ODE
\begin{equation}
	\label{eq:ogm_ode}
	\ddot{X}+\frac{3}{t}\dot{X}+2\nabla f(X)=0
\end{equation}
with $X(0)=x_0$ and $\dot{X}(0)=0$ (see Appendix~\ref{app:limiting_ogm}). Because this ODE is equivalent to the first Bregman Lagrangian flow \eqref{eq:first_family} with $\alpha(t)=\log\frac{2}{t}$ and $\beta(t)=\log\frac{t^2}{2}$, its solution reduces the function value accuracy $f(X(T))-f(x^*)$ with an $O(\|x_0-x^*\|^2/T^2)$ convergence rate. Under the same setting, \citet{suh2022continuous} showed that OGM-G converges to the ODE
\begin{equation}
	\label{eq:ogm-g_ode}
	\ddot{X}+\frac{3}{T-t}\dot{X}+2\nabla f(X)=0
\end{equation}
with $X(0)=x_0$ and $\dot{X}(0)=0$, and showed that the solution to this ODE reduces the squared gradient norm $\Vert\nabla f(X(T))\Vert^2$ with an $O(f(x_0)-f(x^*)/T^2)$ convergence rate. We can observe that the coefficients in \eqref{eq:ogm-g_ode} can be obtained by substituting $t$ with $T-t$ into the coefficient in \eqref{eq:ogm_ode} and vice versa. 

Based on the symmetric relationship between OGM ODE and OGM-G ODE, one might intuitively think that ``\emph{OGM-G ODE is a time-reversed version of OGM ODE}.'' This interpretation, however, might be misleading because the solution to OGM ODE and the solution to OGM-G ODE do not have a time-reversed relationship. In the following paragraph, using the differential kernel \eqref{eq:continuous_FSFO}, we present a different, conceivably more accurate, symmetrical relationship between the two ODEs.

\paragraph{Anti-transpose relationship between OGM ODE and OGM-G ODE.}
Substituting $b^{\mathrm{F}}(t)=3/t$, $b^{\mathrm{G}}(t)=3/(T-t)$, and $c^{\mathrm{F}}(t)=c^{\mathrm{G}}(t)=1$ in \eqref{eq:kernel_for_two}, the differential kernels $H^{\mathrm{F}}(t,\tau)$ corresponding to OGM ODE and $H^{\mathrm{G}}(t,\tau)$ corresponding to OGM-G ODE can be computed as
\begin{align*}
	H^{\mathrm{F}}(t,\tau) & =\frac{2\tau^{3}}{t^{3}}\\
	H^{\mathrm{G}}(t,\tau) & =\frac{2(T-t)^{3}}{(T-\tau)^{3}}.
\end{align*}
Here, we can observe the following \emph{anti-transpose relationship} between two differential kernels:
\begin{equation}
	\label{eq:symmetry_b_c}
	H^{\mathrm{F}}(t,\tau)=H^{\mathrm{G}}(T-\tau,T-t).
\end{equation}
Note that this can also be obtained by using the definition of the differential kernel and the anti-transpose relationship \eqref{eq:symmetry_beta_gamma} between two matrices $\mathbf{H}^{\mathrm{F}}$ and $\mathbf{H}^{\mathrm{G}}$ defined in \eqref{eq:HFHG}. To summarize, the relationships between OGM, OGM-G, and their limiting ODEs are illustrated in Figure~\ref{fig:chart2}.

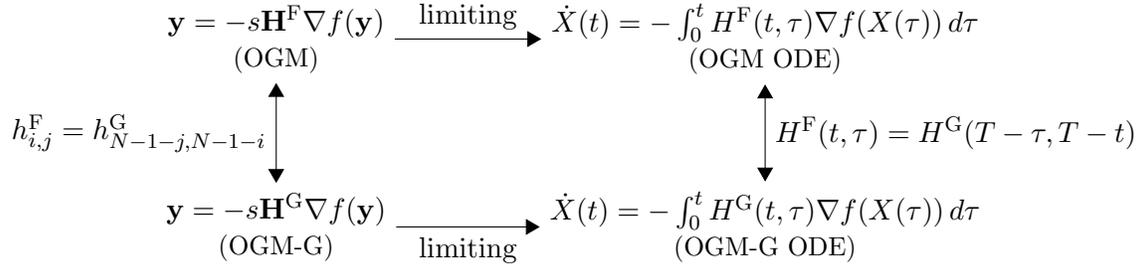
\begin{figure}[!htp]
	\begin{center}
		\begin{tikzpicture}[node distance=6.5cm]
			\node (ogmode) [align=center] {$\dot{X}(t)=-\int_{0}^{t}H^{\mathrm{F}}(t,\tau)\nabla f(X(\tau))\,d\tau$ \\ {\small (OGM ODE)}};
			\node (ogm) [align=center, left of=ogmode] {$\mathbf{y}=-s\mathbf{H}^{\mathrm{F}}\nabla f(\mathbf{y})$ \\ {\small (OGM)}};
			\node (ogmg) [align=center, left of=ogmode, yshift=-2.5cm] {$\mathbf{y}=-s\mathbf{H}^{\mathrm{G}}\nabla f(\mathbf{y})$ \\ {\small (OGM-G)}};
			\node (ogmgode) [align=center, right of=ogmg] {$\dot{X}(t)=-\int_{0}^{t}H^{\mathrm{G}}(t,\tau)\nabla f(X(\tau))\,d\tau$ \\ {\small (OGM-G ODE)}};
			
			\draw [arrow] (ogm) --node [above] {limiting} (ogmode);
			\draw [arrow] (ogmg) --node [below] {limiting} (ogmgode);
			\draw [botharrow] (ogmode) --node [right] {$H^{\mathrm{F}}(t,\tau)=H^{\mathrm{G}}(T-\tau,T-t)$} (ogmgode);
			\draw [botharrow] (ogm) --node [left] {$h_{i,j}^{\mathrm{F}}=h_{N-1-j,N-1-i}^{\mathrm{G}}$} (ogmg);
		\end{tikzpicture}
	\end{center}
	\caption{Anti-transpose relationships between OGM (reducing the function value accuracy), OGM-G (reducing the gradient norm), and their limiting ODEs.}
	\label{fig:chart2}
\end{figure}

\paragraph{A failed attempt to design an ODE that minimizes the gradient norm of strongly convex functions.}
A downside of OGM-G ODE \eqref{eq:ogm-g_ode} is that it exploits only the non-strong convexity of the objective function $f$. Thus, one might want to design an ODE model that minimizes the gradient norm of strongly convex objective functions. Inspired by the symmetric relationship between OGM ODE and OGM-G ODE, one might substitute $t$ with $T-t$ into the coefficients in \nagscsp ODE \eqref{eq:sc_ode} to yield the following ODE:
\begin{equation}
	\label{eq:ogm-g_ode_sc}
	\ddot{X}+2\sqrt{\mu}\dot{X}+\nabla f(X)=0,
\end{equation}
and one might guess that the solution to this ODE reduces the squared gradient norm $\|\nabla f(X(T))\|^2$ with an $O(e^{-\sqrt{\mu}T})$ convergence rate. However, one cannot easily modify the argument in \citep{suh2022continuous} to prove the convergence rate of the gradient norm for \eqref{eq:ogm-g_ode_sc} because their argument depends on the property $\dot{X}(T)=0$, which is not true for the solution to \eqref{eq:ogm-g_ode_sc}. 

\subsection{Proposed dynamics: Unified NAG-G ODE}

In this subsection, we claim that the symmetric counterpart of the unified NAG ODE \eqref{eq:u_ode_oneline} works well for our purpose, unlike the aforementioned failed attempt. The property that the unified NAG ODE is a continuous extension of \nagcsp ODE allows us to use the argument in \citep[Section~4.1]{suh2022continuous}. Substituting $t$ with $T-t$ into the coefficients in the unified NAG ODE \eqref{eq:u_ode_oneline}, we obtain the following ODE:
\begin{equation}
	\label{eq:unified_ogm-g}
	\ddot{X}+\left(\frac{\sqrt{\mu}}{2}\tanh\left(\frac{\sqrt{\mu}}{2}(T-t)\right)+\frac{3}{T-t}\cothc\left(\frac{\sqrt{\mu}}{2}(T-t)\right)\right)\dot{X}+\nabla f(X)=0.
\end{equation}
We refer to this ODE with the initial conditions $X(0)=x_0$ and $\dot{X}(0)=0$ as the  \emph{unified NAG~-~G ODE}. Clearly, this ODE has a unique solution in $C^1([0,T),\mathbb{R}^n)$.\footnote{Sketch of the proof: For any $\epsilon\in(0,T/2)$, the existence and uniqueness of solution on $[0,T-\epsilon]$ follows from Cauchy-Lipschitz theorem \citep[Theorem~25]{teschl2012ordinary}. Paste these solutions on $[0,T)=\cup_{\epsilon\in(0,T/2)}[0,T-\epsilon)$.} We can continuously extend this solution to $t=T$ with $\dot{X}(T)=0$ and $\ddot{X}(T)=\lim_{t\rightarrow T^{-}}\frac{\dot{X}(t)}{t-T}=\frac{1}{2}\nabla f(X(T))$ (see Appendix~\ref{app:dotX_and_ddotX}). To analyze the convergence rate, we use the Lyapunov analysis again.
\begin{theorem}
	\label{thm:mainthm_unified_ogm-g}
	For the solution $X:[0,T]\rightarrow\mathbb{R}^n$ to the unified NAG-G ODE \eqref{eq:unified_ogm-g}, the continuous-time energy function
	\begin{equation}
		\label{eq:energy_ogm-g}
		\begin{aligned}
			\mathcal{E}(t) & =\frac{4}{(T-t)^{2}}\cschc^{2}\left(\frac{\sqrt{\mu}}{2}(T-t)\right)\left(f(X(t))-f(X(T))\right)\\
			& \quad-\frac{8}{(T-t)^{4}}\cschc^{4}\left(\frac{\sqrt{\mu}}{2}(T-t)\right)\left\Vert X(t)-X(T)\right\Vert ^{2}\\
			& \quad+\frac{8}{(T-t)^{4}}\cschc^{2}\left(\frac{\sqrt{\mu}}{2}(T-t)\right)\cothc^{2}\left(\frac{\sqrt{\mu}}{2}(T-t)\right)\\
			& \qquad\times\left\Vert X(t)+\frac{T-t}{2}\tanhc\left(\frac{\sqrt{\mu}}{2}(T-t)\right)\dot{X}(t)-X(T)\right\Vert ^{2}.
		\end{aligned}
	\end{equation}
	is monotonically decreasing on $[\tst,T)$.
\end{theorem}
The proof of Theorem~\ref{thm:mainthm_unified_ogm-g} can be found in Appendix~\ref{app:unified_ogm-g}. By L'H\^{o}pital's rule, we have
\begin{align*}
	\lim_{t\rightarrow T^{-}}\frac{f(X(t))-f(X(T))}{(T-t)^{2}} & =\lim_{t\rightarrow T^{-}}\frac{1}{2}\left\langle \frac{\dot{X}(t)}{t-T},\nabla f(X)\right\rangle =\frac{1}{4}\left\Vert \nabla f(X(T))\right\Vert ^{2}\\
	\lim_{t\rightarrow T^{-}}\frac{X(t)-X(T)}{(T-t)^{2}} & =\lim_{t\rightarrow T^{-}}\frac{\dot{X}(t)}{2(t-T)}=\frac{1}{4}\nabla f(X(T)).
\end{align*}
It follows from $\cschc(0)=\cothc(0)=1$ that
\begin{align*}
	& \lim_{t\rightarrow T^{-}}\mathcal{E}(t)\\
	& =\lim_{t\rightarrow T^{-}}\left(4\cdot\frac{f(X(t))-f(X(T))}{(T-t)^{2}}-8\left\Vert \frac{X(t)-X(T)}{(T-t)^{2}}\right\Vert ^{2}+8\left\Vert \frac{X(t)-X(T)}{(T-t)^{2}}-\frac{\dot{X}(t)}{2(t-T)}\right\Vert ^{2}\right)\\
	& =\left\Vert \nabla f(X(T))\right\Vert ^{2}-\frac{1}{2}\left\Vert \nabla f(X(T))\right\Vert ^{2}+0\\
	& =\frac{1}{2}\left\Vert \nabla f(X(T))\right\Vert ^{2}.
\end{align*}
Writing $\lim_{t\rightarrow T^{-}}\mathcal{E}(t)\leq\mathcal{E}(\tst)$ explicitly, we obtain the following result.
\begin{corollary}
	\label{cor:maincor_unified_ogm-g}
	The solution $X$ to the unified NAG-G ODE \eqref{eq:unified_ogm-g} satisfies the inequality
	\begin{equation}
		\label{eq:guarantee_unified_ogm-g}
		\begin{aligned}
			\left\Vert \nabla f(X(T))\right\Vert ^{2} & \leq \frac{8}{T^{2}}\cschc^{2}\left(\frac{\sqrt{\mu}}{2}T\right)\left(f(x_{0})-f\left(X(T)\right)+\frac{\mu}{2}\left\Vert x_{0}-X(T)\right\Vert ^{2}\right)\\
			& \leq \frac{8}{T^{2}}\cschc^{2}\left(\frac{\sqrt{\mu}}{2}T\right)\left(f(x_{0})-f\left(x^*\right)+\frac{\mu}{2}\left\Vert x_{0}-X(T)\right\Vert ^{2}\right).
		\end{aligned}
	\end{equation}
\end{corollary}
Since $\cschc^2$ is decreasing on $[0,\infty)$, Corollary~\ref{cor:maincor_unified_ogm-g} implies that the unified NAG-G ODE \eqref{eq:u_ode} reduces the squared gradient norm with an $O\left(1/T^2\right)$ convergence rate regardless of the value of $\mu \geq 0$. When $\mu>0$, since $\frac{1}{T^{2}}\cschc^{2}\left(\frac{\sqrt{\mu}}{2}T\right)\sim\mu e^{-\sqrt{\mu}T}$ as $T\rightarrow\infty$, the unified NAG-G ODE reduces the squared gradient norm with an $O\left(e^{-\sqrt{\mu}T}\right)$ convergence rate. Combining these bounds, we conclude that the unified NAG-G ODE reduces the squared gradient norm with the following convergence rate:
\[
\left\Vert \nabla f(X(T))\right\Vert ^{2}\leq O\left(\min\left\{ 1/T^{2},e^{-\sqrt{\mu}T}\right\} \left(f(x_{0})-f\left(x^*\right)+\frac{\mu}{2}\left\Vert x_{0}-X(T)\right\Vert ^{2}\right)\right).
\]

\paragraph{Anti-transpose relationship between the unified NAG ODE and the unified NAG-G ODE.}

The differential kernels $H^{\mathrm{F}}(t,\tau)$ corresponding to the unified NAG ODE and $H^{\mathrm{G}}(t,\tau)$ corresponding to the unified NAG-G ODE can be computed as (see Appendix~\ref{app:equivalent_of_ode})
\begin{align*}
	H^{\mathrm{F}}(t,\tau) & =\frac{\tau^{3}\sinhc^{3}\left(\frac{\sqrt{\mu}}{2}\tau\right)\cosh\left(\frac{\sqrt{\mu}}{2}\tau\right)}{t^{3}\sinhc^{3}\left(\frac{\sqrt{\mu}}{2}t\right)\cosh\left(\frac{\sqrt{\mu}}{2}t\right)}\\
	H^{\mathrm{G}}(t,\tau) & =\frac{(T-t)^{3}\sinhc^{3}\left(\frac{\sqrt{\mu}}{2}(T-t)\right)\cosh\left(\frac{\sqrt{\mu}}{2}(T-t)\right)}{(T-\tau)^{3}\sinhc^{3}\left(\frac{\sqrt{\mu}}{2}(T-\tau)\right)\cosh\left(\frac{\sqrt{\mu}}{2}(T-\tau)\right)}.
\end{align*}
Remarkably, there is an anti-transpose relationship \eqref{eq:symmetry_b_c} between these differential kernels, like the one between the differential kernels corresponding to OGM ODE (which minimizes the function value accuracy, similarly to what the unified NAG ODE does) and OGM-G ODE (which minimizes the gradient norm, similarly to what the unified NAG-G ODE does).

\section{Numerical Experiments}
\label{sec:experiments}

In this section, we validate the performance of the unified NAG \eqref{eq:unified_nag_specific} for a toy problem and the logistic regression problem, and we also compare our method with \nagcsp \eqref{eq:nag-c} and \nagscsp \eqref{eq:nag-sc}. For each problem, we empirically observed that {the unified NAG attains the advantages of both \nagcsp and \nagsc}. 

\paragraph{Toy problem.}
We consider the problem
\begin{equation}
	\label{eq:toy}
	\min_{(x,y)\in\mathbb{R}^2} \; f(x,y)=\frac{\mu}{2}x^{2}+0.005y^{2}.
\end{equation}
This problem is strongly convex with parameter $\min\left\{ \mu,0.01\right\} $. We set the initial point and the algorithmic stepsize as $\left(x_0,y_0\right)=(1,1)$ and $s=1$. When $\mu$ is large ($\mu=10^{-3}$), \cref{fig:exp2a} shows that \nagscsp outperforms \nagcsp and that the unified NAG behaves like \nagsc. When $\mu=10^{-4}$, \cref{fig:exp2b} shows that the unified NAG behaves like \nagcsp in the early stages and behaves like \nagscsp in the late stages. When $\mu$ is small ($\mu=10^{-7}$), \cref{fig:exp2c} shows that \nagcsp outperforms \nagscsp at least in the early stages and that the unified NAG behaves like \nagc. In each case, the performance of the unified NAG is comparable to the better choice between \nagcsp and \nagsc. The trajectories of the algorithms are shown in Figures~\ref{fig:exp2d}, \ref{fig:exp2e}, and \ref{fig:exp2f}. We can see that \nagscsp converges with more severe oscillation compared to \nagcsp and the unified NAG, particularly when the strong convexity parameter $\mu$ is small. This result matches the damping system interpretation in \cref{sec:unified_ode}: \nagscsp behaves like an underdamped system when $\mu$ is small, while our unified NAG always behaves like an overdamped system in the early stages.

\begin{figure}[ht]
	\centering
	\subfigure[Errors $f-f^*$, $\mu=10^{-3}$]{\includegraphics[width=0.32\textwidth]{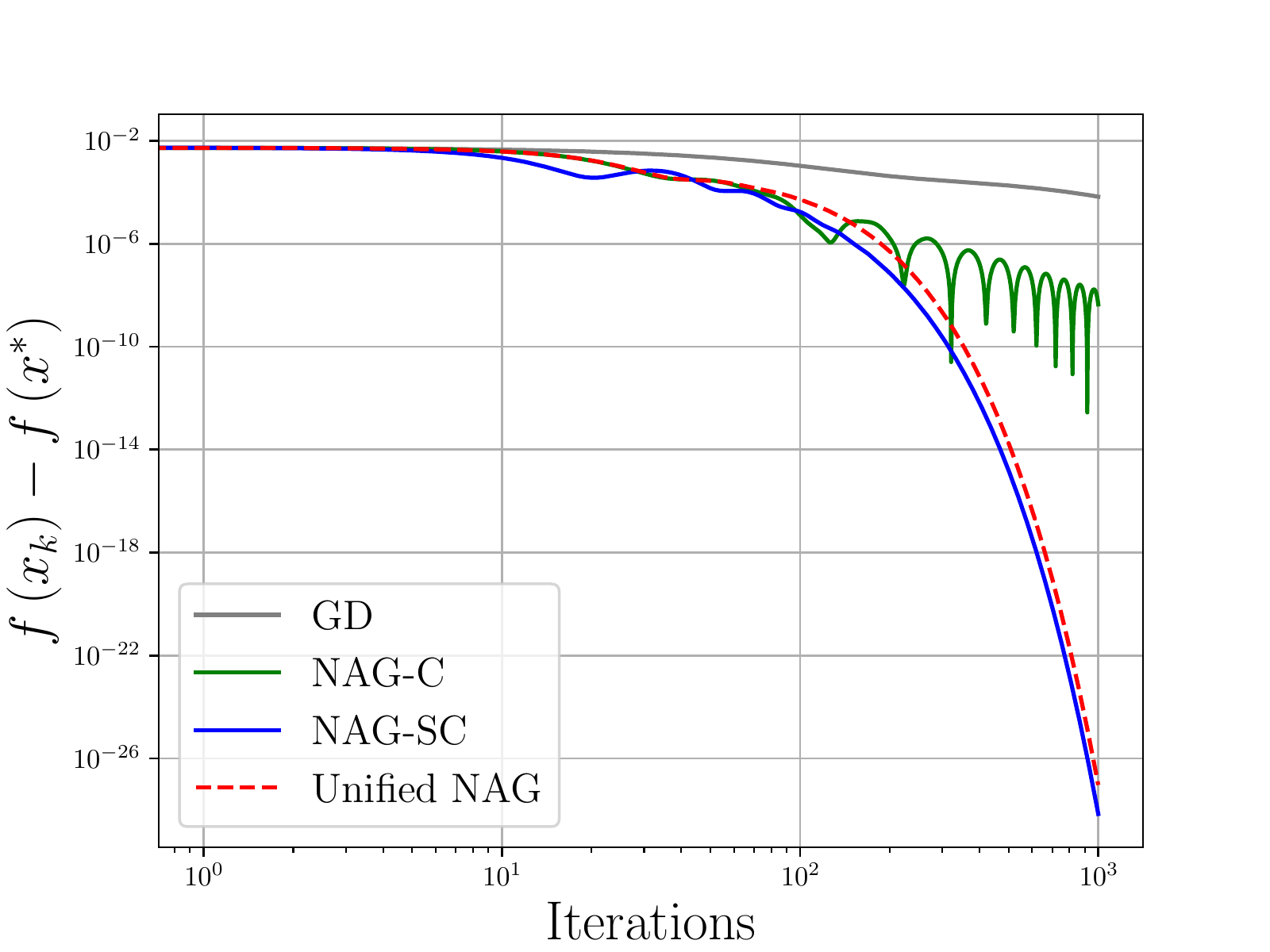}\label{fig:exp2a}}
	\subfigure[Errors $f-f^*$, $\mu=10^{-4}$]{\includegraphics[width=0.32\textwidth]{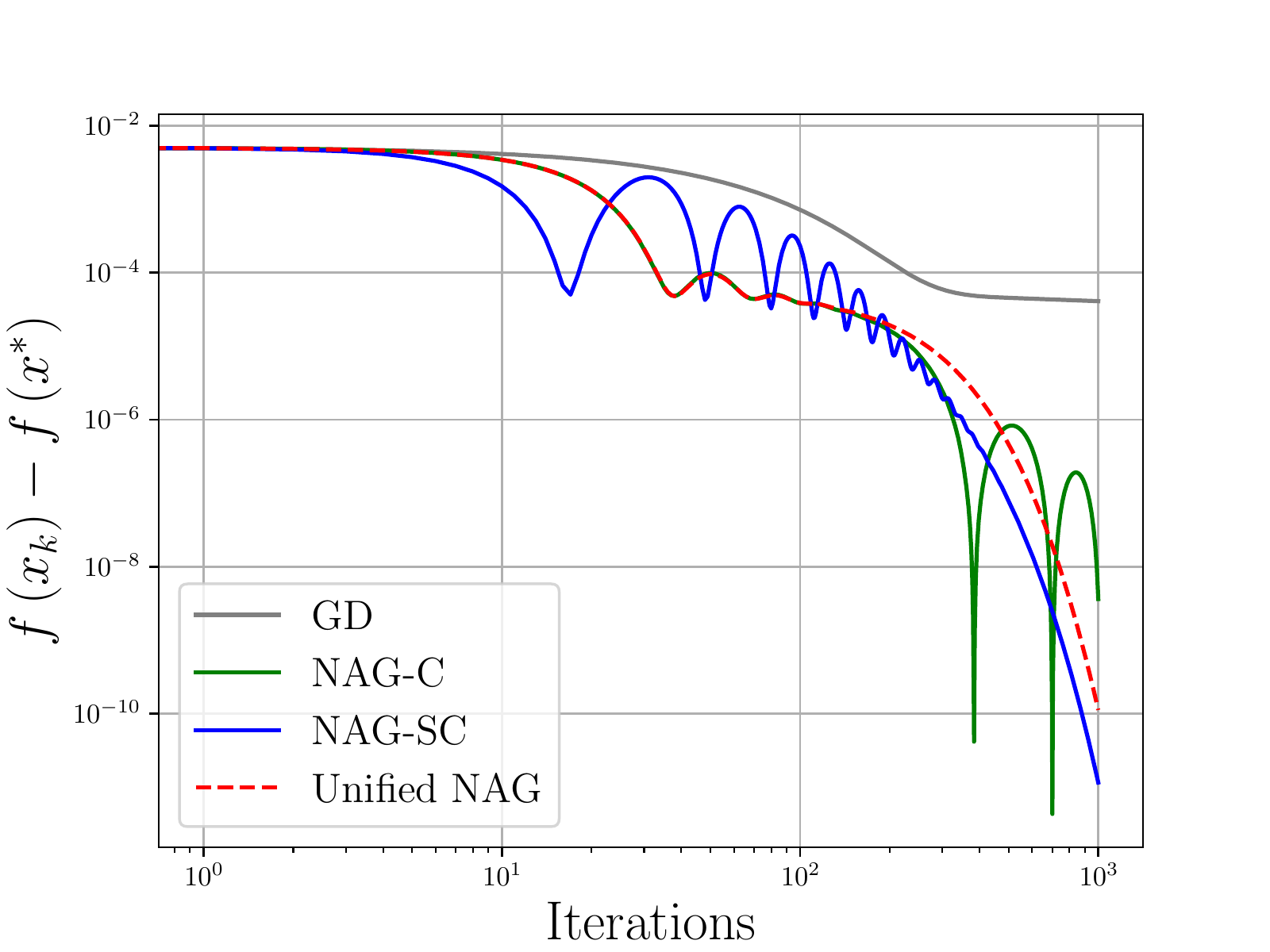}\label{fig:exp2b}}
	\subfigure[Errors $f-f^*$, $\mu=10^{-7}$]{\includegraphics[width=0.32\textwidth]{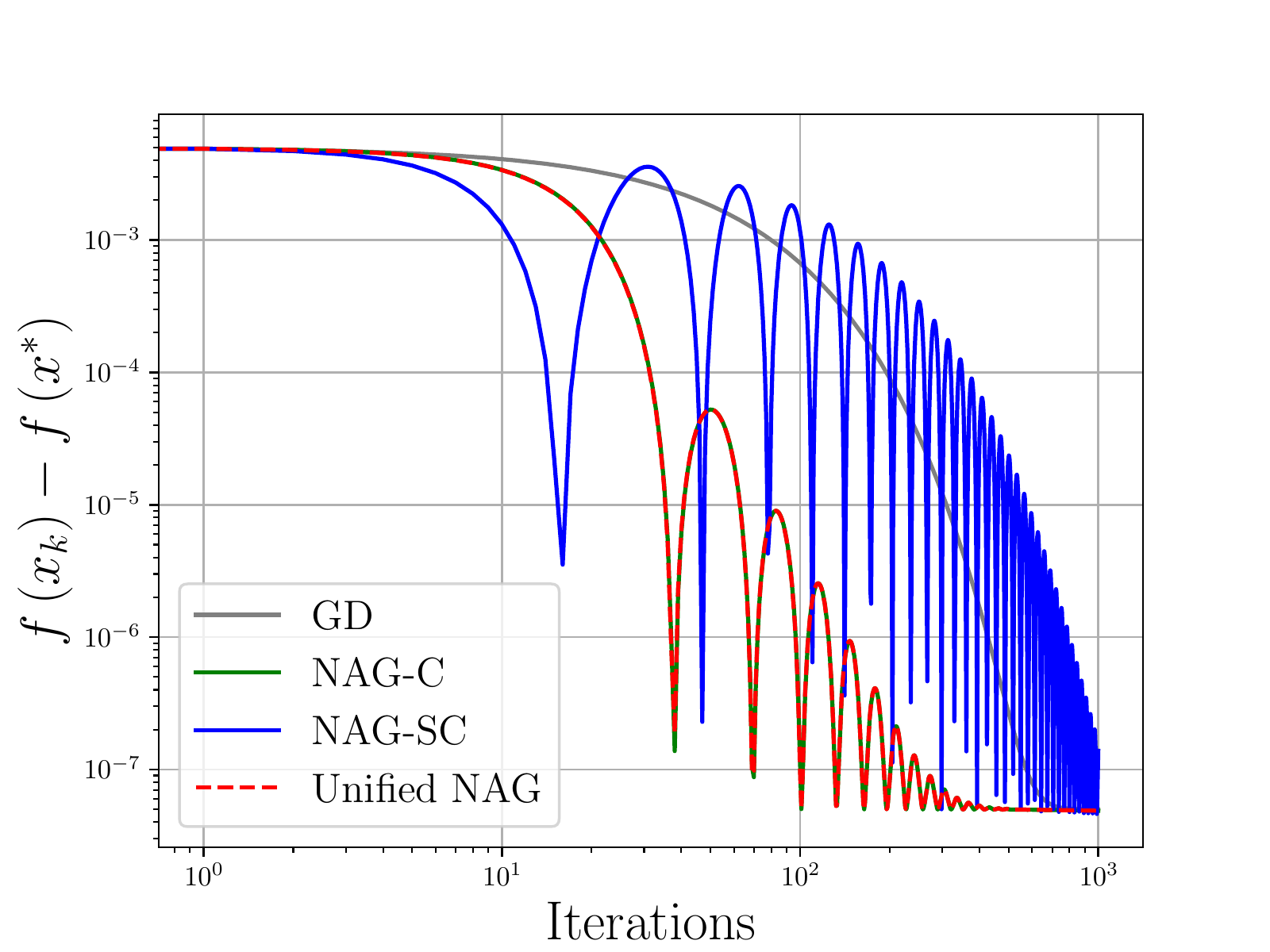}\label{fig:exp2c}}
	\newline
	\subfigure[Trajectories, $\mu=10^{-3}$]{\includegraphics[width=0.32\textwidth]{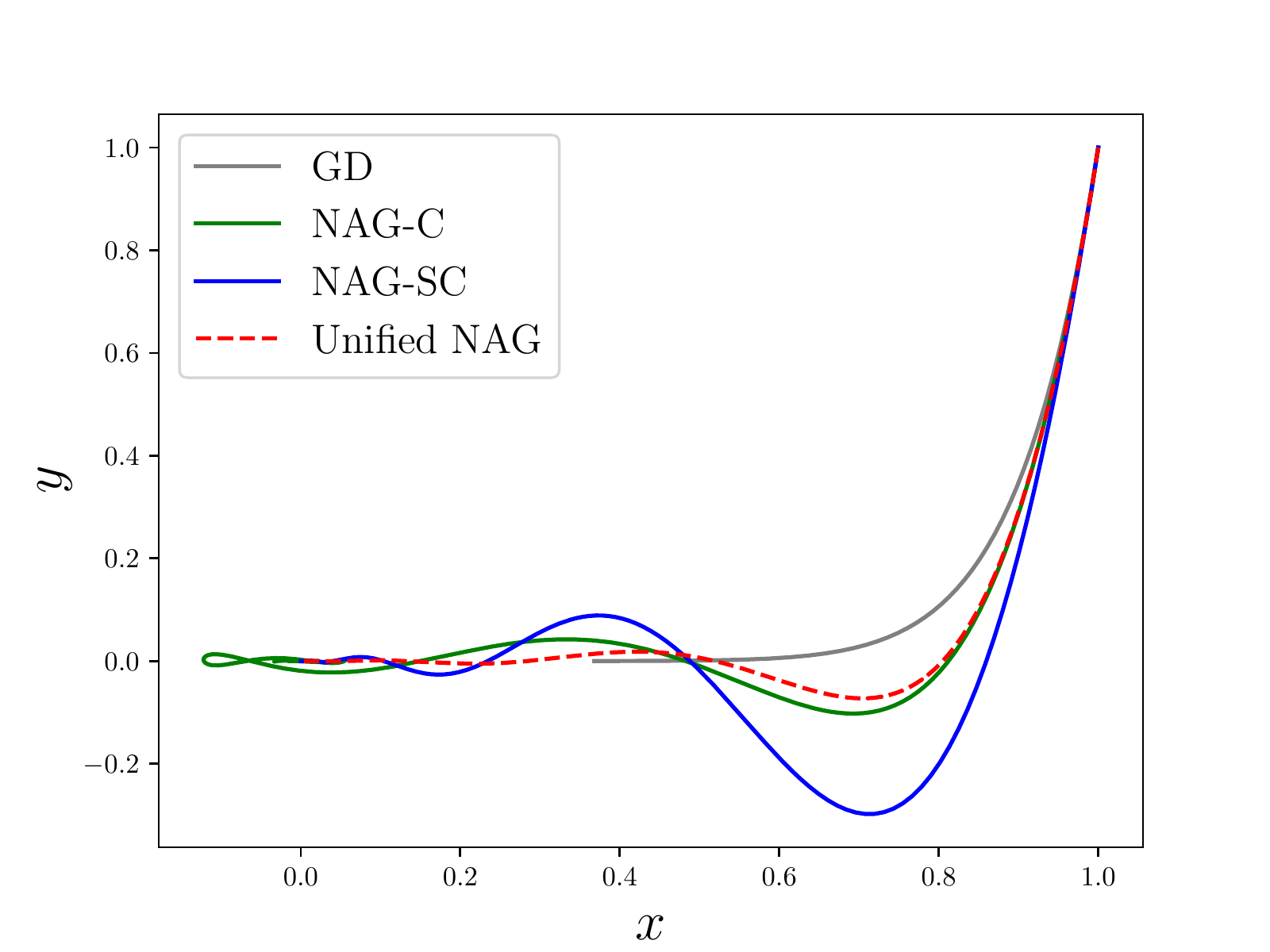}\label{fig:exp2d}}
	\subfigure[Trajectories, $\mu=10^{-4}$]{\includegraphics[width=0.32\textwidth]{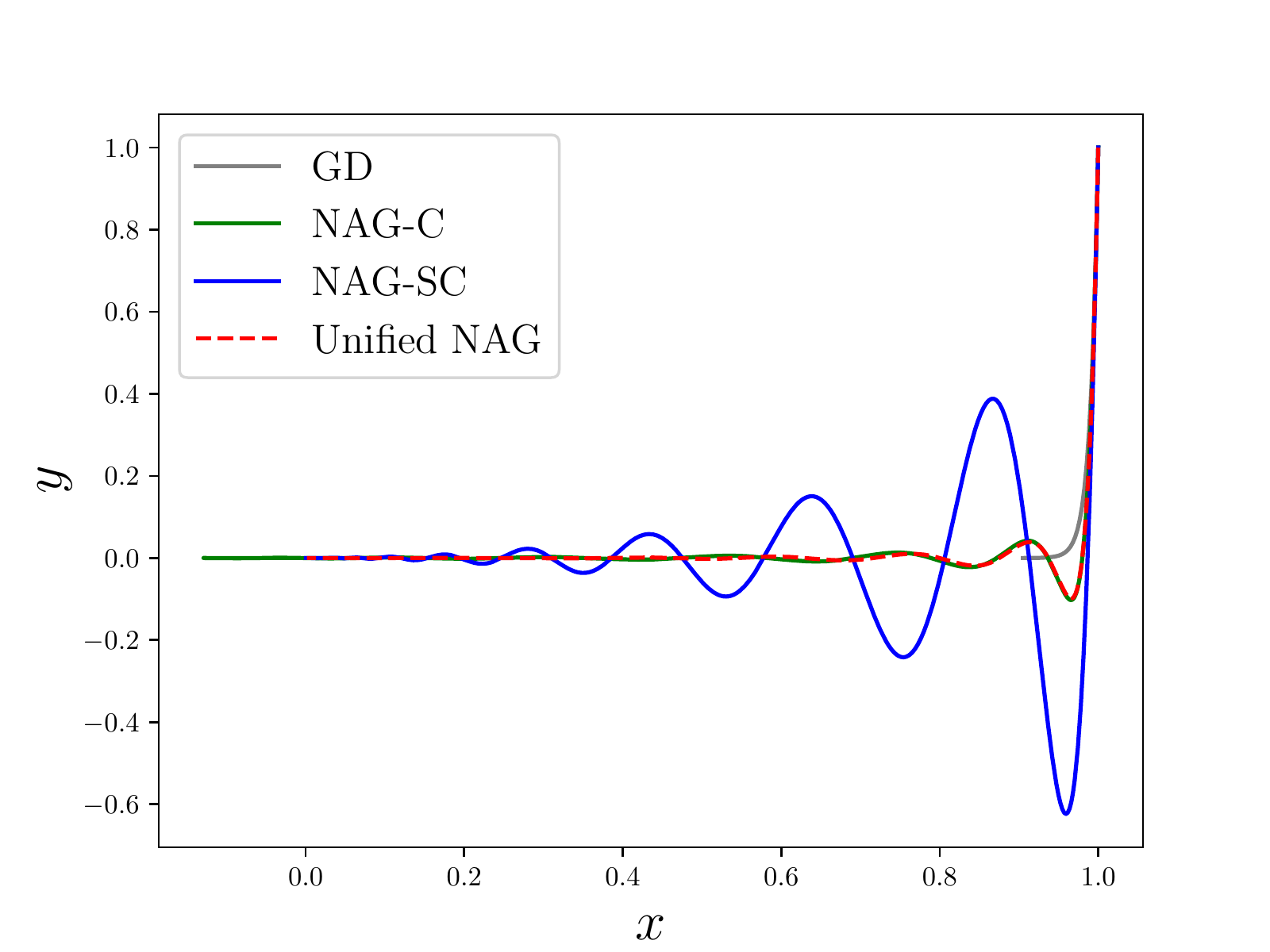}\label{fig:exp2e}}
	\subfigure[Trajectories, $\mu=10^{-7}$]{\includegraphics[width=0.32\textwidth]{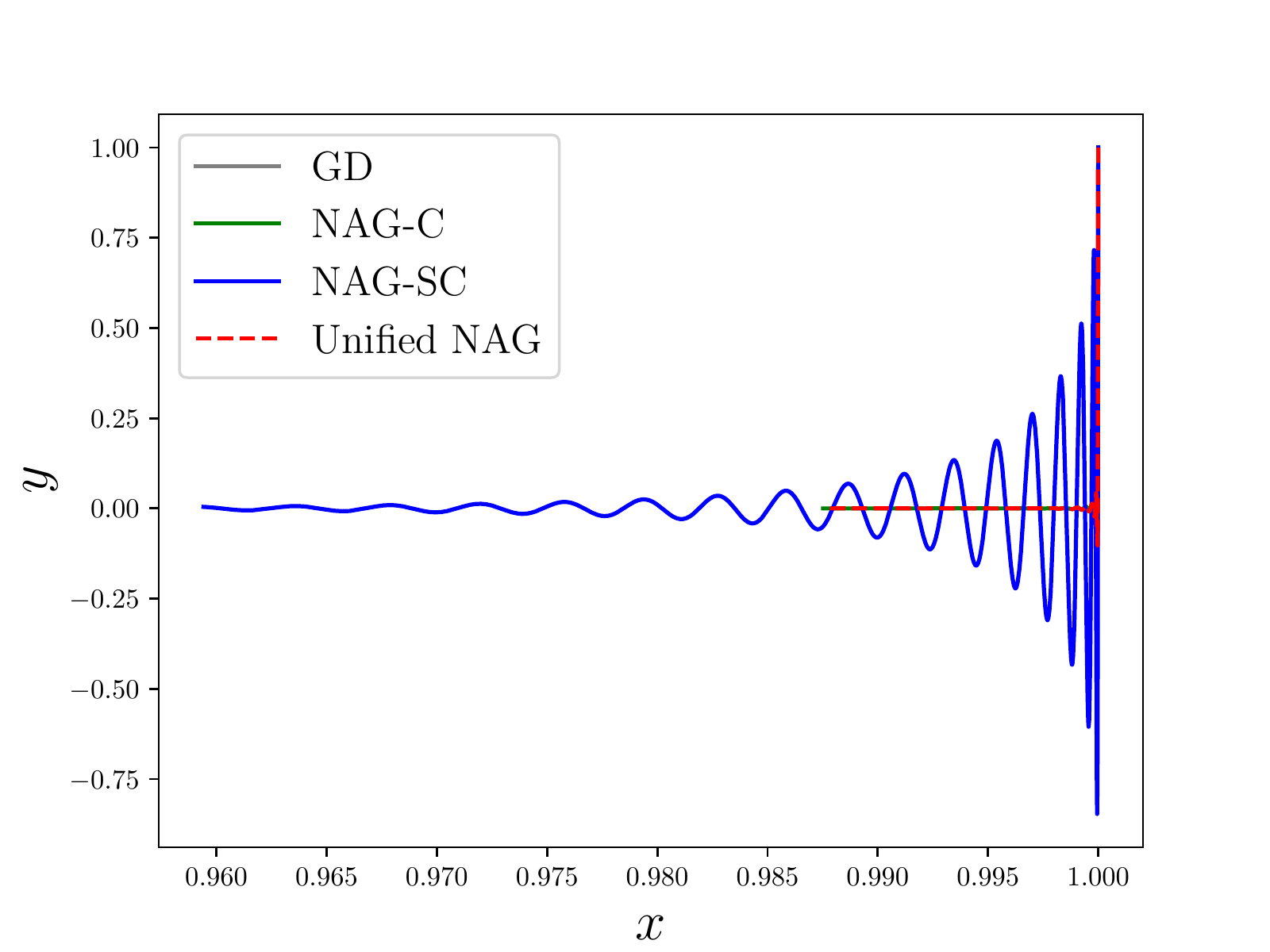}\label{fig:exp2f}}
	\caption{Results for the problem with the objective function $f(x,y)=\frac{\mu}{2}x^{2}+0.005y^{2}$ and the initial state $x_0=(1,1)$.}\label{fig:experiments2}
\end{figure}

\paragraph{$\ell_2$-regularized logistic regression.}

We now consider the $\ell_2$-regularized logistic regression problem
\begin{equation}
	\label{eq:logistic}
	\min_{x\in\mathbb{R}^{n}}\;f(x)=\frac{1}{m}\left(\sum_{i=1}^{m}\left(-y_{i}a_{i}^{T}x+\log\left(1+e^{a_{i}^{T}x}\right)\right)+\lambda\left\Vert x\right\Vert ^{2}\right),
\end{equation}
where $a_i\in\mathbb{R}^n$ and $y_i\in\{0,1\}$ for $i=1,2,\ldots,m$. Then, \eqref{eq:logistic} is the problem \eqref{eq:ml_problem} with the convex functions $f_i(x)=-y_{i}a_{i}^{T}x+\log(1+e^{a_{i}^{T}x})$ and the $\ell_2$-regularization term $R(x)=\left\Vert x\right\Vert ^{2}$. As mentioned in \cref{sec:motivation}, the function $f$ is $\mu$-strongly convex with $\mu=\frac{2\lambda}{m}$. We set $s=0.01$ and choose  the sample size and the dimension as $m=100$ and $n=20$, respectively. Following \citep{su2014}, we use a synthetically generated data set: the entries of $a_i$ are generated by the Gaussian distribution $\mathcal{N}(0,1)$, and the labels $y_i\in\{ 0,1\}$ are generated by the logistic model $P\left(y_{i}\right)=1=\frac{1}{1+e^{-a_{i}^{T}x^0}}$, where the entries of $x^0$ are generated by the Gaussian distribution $\mathcal{N}(0,1/100)$. The results are shown in \cref{fig:experiments3}. Again, we can observe that \nagscsp outperforms \nagcsp when $\mu$ is large and underperforms \nagcsp when $\mu$ is small. 
In each case, the performance of unified NAG is on par with the better one among \nagcsp and \nagsc.

\begin{figure}[ht]
	\centering
	\subfigure[Errors $f-f^*$, $\lambda=5$]{\includegraphics[width=0.32\textwidth]{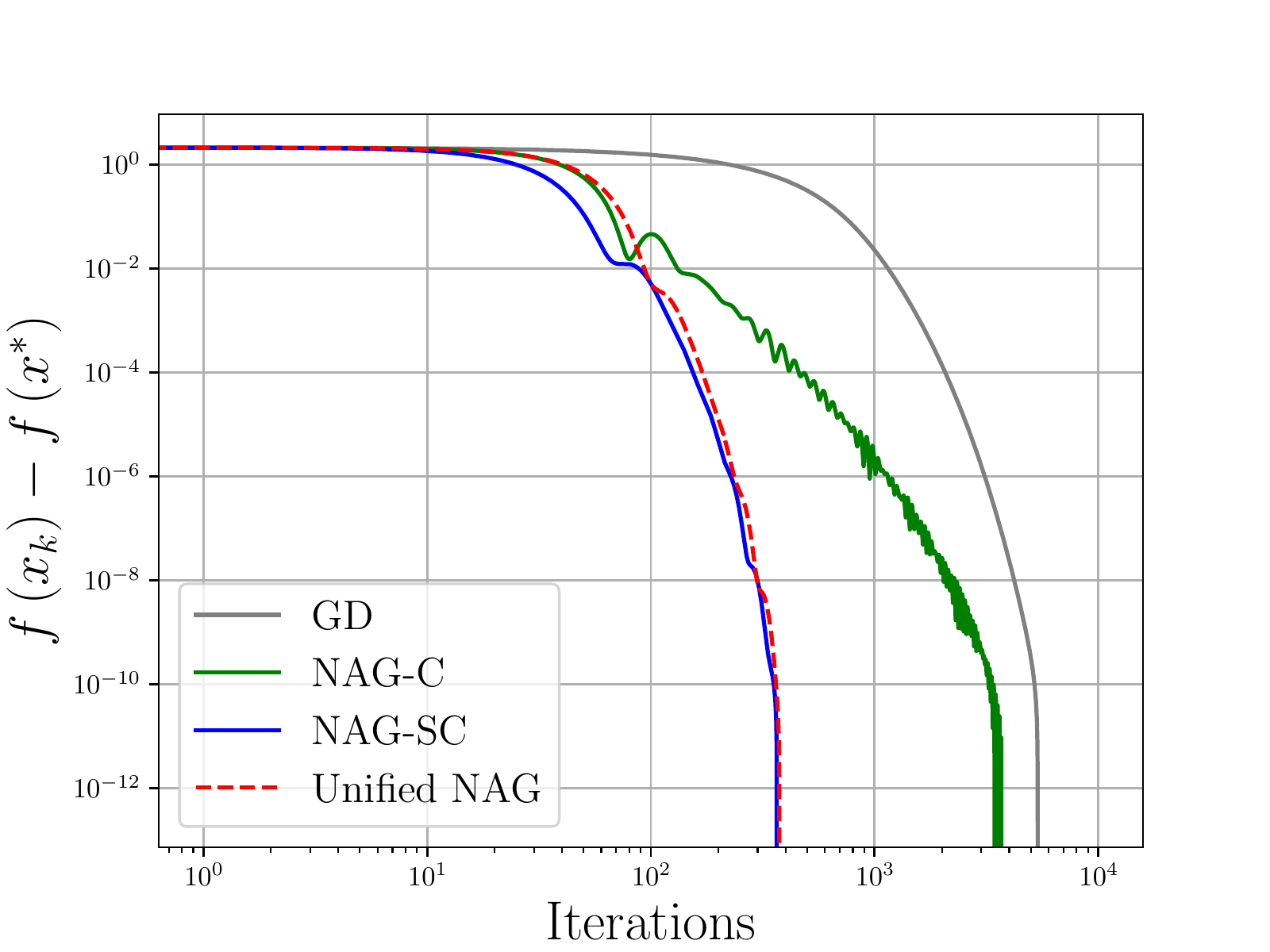}\label{fig:exp3a}}
	\subfigure[Errors $f-f^*$, $\lambda=5\cdot10^{-2}$]{\includegraphics[width=0.32\textwidth]{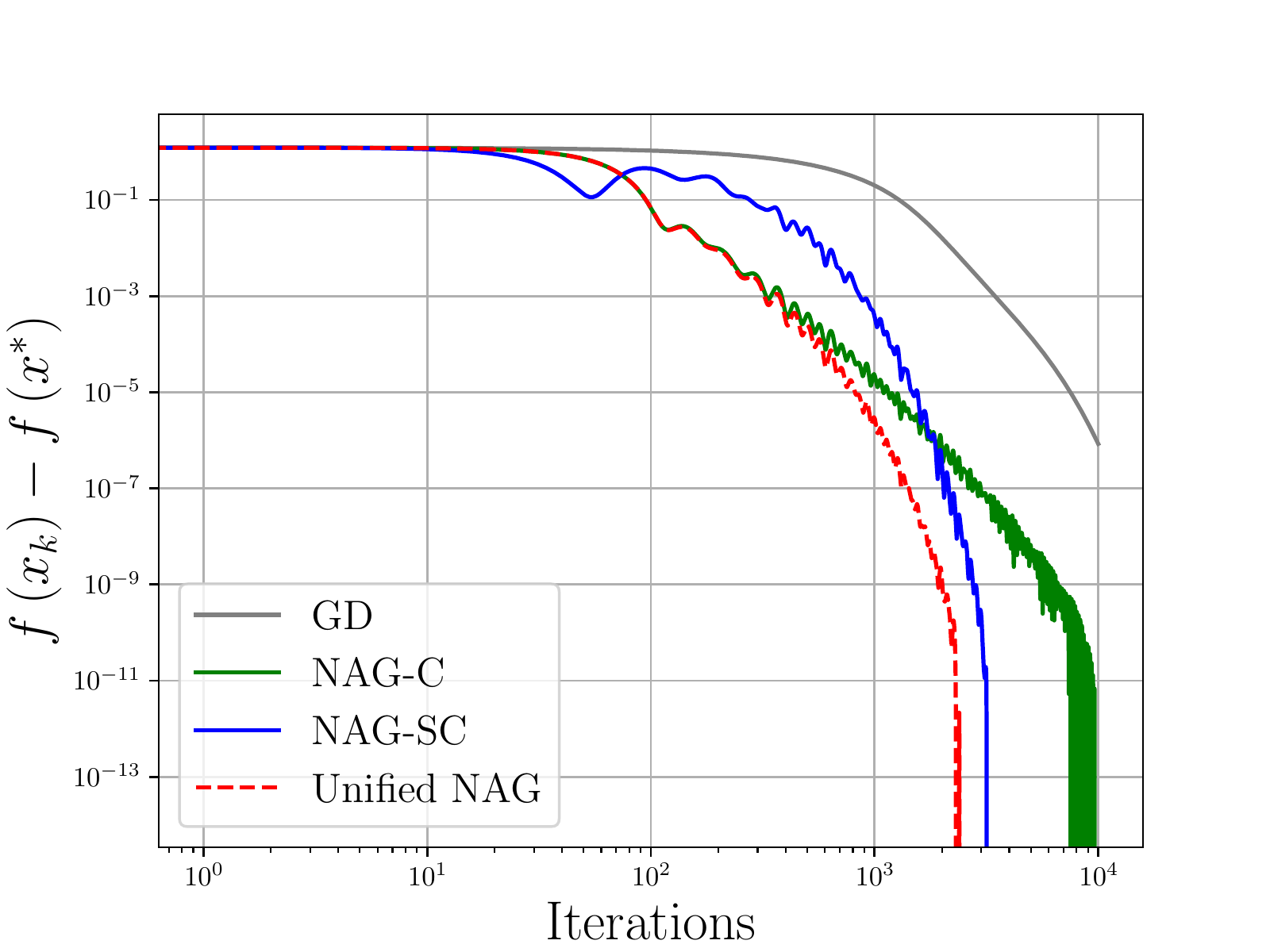}\label{fig:exp3b}}
	\subfigure[Errors $f-f^*$, $\lambda=5\cdot10^{-4}$]{\includegraphics[width=0.32\textwidth]{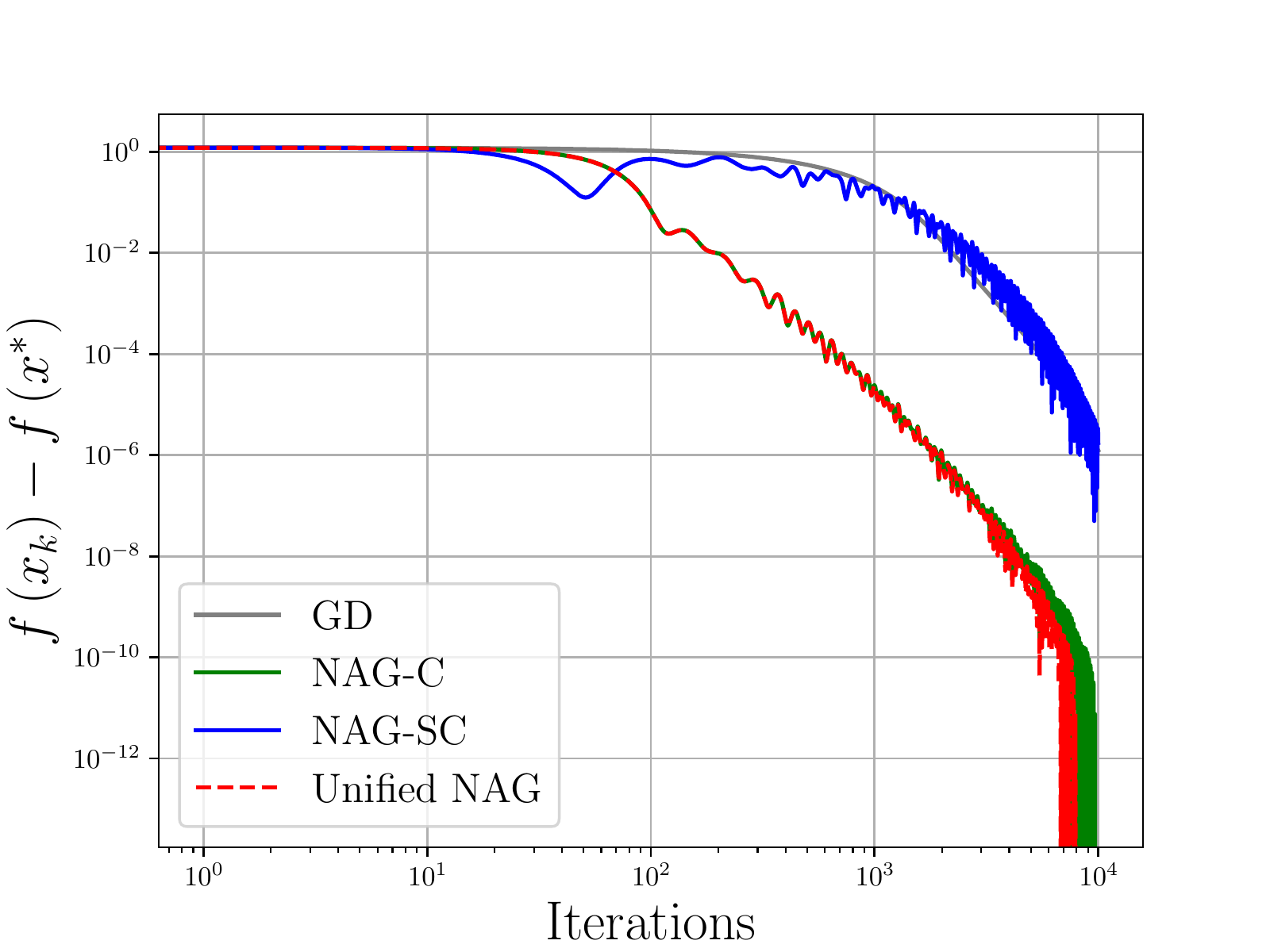}\label{fig:exp3c}}
	\caption{Results for the $\ell_2$-regularized logistic regression problem.}\label{fig:experiments3}
\end{figure}

\section{Conclusions}

In this paper, we examined and resolved inconsistencies between the momentum algorithms and ODE models for convex  and strongly convex cases. To bridge the gap between the two cases,  we proposed the unified Bregman Lagrangian \eqref{eq:unified_lagrangian}, the unified NAG ODE \eqref{eq:u_ode_oneline}, and the unified NAG \eqref{eq:unified_nag_specific}. Because our algorithm, ODE model and Lagrangian are continuous in $\mu$ and recover the corresponding counterparts for non-strongly convex cases (see Figure~\ref{fig:chart}), they can be viewed as continuous extensions of the \nagc, \nagcsp ODE, and the first Bregman Lagrangian. We theoretically and empirically showed that unlike \nagsc, the unified NAG has a better convergence rate compared to \nagcsp regardless of the values of $\mu$, which is quite significant in practice, as mentioned in Section~\ref{sec:motivation}. Based on the Lagrangian formalism, we proposed the unified accelerated tensor flow \eqref{eq:hyperbolic_flow} and scheme \eqref{eq:hyperbolic_scheme}, achieving exponential convergence rates in the higher-order setting. Lastly, hinted from the unified NAG ODE, we designed the unified NAG-G ODE \eqref{eq:unified_ogm-g}, a novel dynamical system that minimizes the gradient norm of strongly convex functions. Using our novel tool, the differential kernel \eqref{eq:continuous_FSFO}, we discovered an anti-transpose relationship \eqref{eq:symmetry_b_c} between OGM ODE and OGM-G ODE. Surprisingly, such relationship can also be found between the unified NAG ODE and the unified NAG-G ODE.


\acks{We thank Prof. Ernest K. Ryu at Seoul National University for providing feedback on this work. This work was supported in part by Samsung Electronics, the National Research Foundation of Korea funded by MSIT(2020R1C1C1009766), and the Information and Communications Technology Planning and Evaluation (IITP) grant funded by MSIT(2022-0-00124, 2022-0-00480).}


\newpage

\appendix

\section{Existing Unified Dynamics}
\label{app:existing_unified}

\subsection{Relationship between the rescaled original NAG flow and the unified Bregman Lagranfian flow}
\label{app:existing_unified1}

First, we show that the rescaled original NAG flow \eqref{eq:original_nag_flow_scaled} can be expressed as the unified Bregman Lagrangian flow \eqref{eq:unified_family}. Given the parameter function $a(t)$ and the constant $\gamma_0$ of the rescaled original NAG flow, we can write the functions $\gamma(t)$ and $b(t)$ involved in \eqref{eq:original_nag_flow} and \eqref{eq:original_nag_flow_scaled} as
\begin{align*}
	\gamma(t) & =\mu+\left(\gamma_{0}-\mu\right)e^{-t}\\
	b(t) & =\mu+\left(\gamma_{0}-\mu\right)e^{-\int_{0}^{t}a(s)\,ds}.
\end{align*}
We define the functions $\alpha(t)$ and $\beta(t)$ as
\begin{equation}
	\label{eq:alpha_beta_for_original}
	\begin{aligned}
		\alpha(t) & =\log a(t)\\
		\beta(t) & =\log\left(\frac{1}{\gamma_{0}-\mu}\right)+\int_{0}^{t}a(s)\,ds.
	\end{aligned}
\end{equation}
Then, we have
\[
\frac{\dot{\beta}e^{\beta}}{1+\mu e^{\beta}}=\frac{e^{\alpha+\beta}}{1+\mu e^{\beta}}=\frac{e^{\alpha}}{\mu+e^{-\beta}}=\frac{a(t)}{\mu+\left(\gamma_{0}-\mu\right)e^{-\int_{0}^{t}a(s)\,ds}}=\frac{a(t)}{b(t)}.
\]
Thus, the rescaled original NAG flow is equivalent to the unified Bregman Lagrangian flow with the parameter functions \eqref{eq:alpha_beta_for_original} and the Euclidean distance-generating function $h(x)=\frac{1}{2}\|x\|^2$.

Conversely, we show that if the ideal scaling conditon \eqref{eq:ideal_scaling_b} holds with equality and the distance-generating function $h$ is Euclidean, then the unified Bregman Lagrangian flow can be written as the rescaled original NAG flow. Given the parameter functions $\alpha(t)$ and $\beta(t)$ of the unified Bregman Lagrangian flow, we define the function $a(t)$ and the constant $\gamma_0$ as 
\begin{align*}
	a(t) & =e^{\alpha(t)}\\
	\gamma_{0} & =\mu+e^{-\beta(0)}.
\end{align*}
Then, because
\[
b(t)=\mu+\left(\gamma_{0}-\mu\right)e^{-\int_{0}^{t}a(s)\,ds}=\mu+e^{-\beta(t)},
\]
we can write the rescaled original NAG flow as
\begin{align*}
	\dot{X}(t) & =e^{\alpha(t)}(Z(t)-X(t))\\
	\left(\mu+e^{-\beta(t)}\right)\dot{Z}(t) & =e^{\alpha(t)}(\mu X(t)-\mu Z(t)-\nabla f(X(t))),
\end{align*}
which is equivalent to the unified Bregman Lagrangian flow if the ideal scaling conditon \eqref{eq:ideal_scaling_b} holds with equality and $h(x)=\frac{1}{2}\|x\|^2$.

\subsection{Relationship between the rescaled original NAG flow with specific parameters and the unified NAG system}
\label{app:existing_unified2}

In particular, given $\gamma>0$, one can choose the function $a(t)$ in the rescaled original NAG flow as \citep[see][Equation~70]{luo2021differential}
\begin{equation}
	\label{eq:rescaled_at}
	a(t)=\begin{cases}
		\frac{2\sqrt{\gamma_{0}}}{\sqrt{\gamma_{0}}t+2}, & \textrm{if }\mu=0,\\
		\sqrt{\mu}\cdot\frac{e^{\sqrt{\mu}t}-\frac{\sqrt{\mu}-\sqrt{\gamma_{0}}}{\sqrt{\mu}+\sqrt{\gamma_{0}}}}{e^{\sqrt{\mu}t}+\frac{\sqrt{\mu}-\sqrt{\gamma_{0}}}{\sqrt{\mu}+\sqrt{\gamma_{0}}}}, & \textrm{if }\mu>0.
	\end{cases}
\end{equation}
In this case, we have $b(t)=(a(t))^2$. Thus, the rescaled original flow with these functions can be written as
\begin{align*}
	\dot{X}(t) & =\frac{2\sqrt{\gamma_{0}}}{\sqrt{\gamma_{0}}t+2}(Z(t)-X(t))\\
	\dot{Z}(t) & =-\frac{\sqrt{\gamma_{0}}t+2}{2\sqrt{\gamma_{0}}}-\nabla f(X(t))
\end{align*}
when $\mu=0$, and
\begin{align*}
	\dot{X}(t) & =\sqrt{\mu}\cdot\frac{e^{\sqrt{\mu}t}-\frac{\sqrt{\mu}-\sqrt{\gamma_{0}}}{\sqrt{\mu}+\sqrt{\gamma_{0}}}}{e^{\sqrt{\mu}t}+\frac{\sqrt{\mu}-\sqrt{\gamma_{0}}}{\sqrt{\mu}+\sqrt{\gamma_{0}}}}(Z(t)-X(t))\\
	\dot{Z}(t) & =\frac{1}{\sqrt{\mu}}\cdot\frac{e^{\sqrt{\mu}t}+\frac{\sqrt{\mu}-\sqrt{\gamma_{0}}}{\sqrt{\mu}+\sqrt{\gamma_{0}}}}{e^{\sqrt{\mu}t}-\frac{\sqrt{\mu}-\sqrt{\gamma_{0}}}{\sqrt{\mu}+\sqrt{\gamma_{0}}}}(\mu X(t)-\mu Z(t)-\nabla f(X(t)))
\end{align*}
when $\mu>0$. In the non-strongly convex case, it is easy to observe that this ODE system converges to \nagcsp system \eqref{eq:c_ode_twoline} as $\gamma_0\to\infty$. In the strongly convex case, because $\frac{\sqrt{\mu}-\sqrt{\gamma_{0}}}{\sqrt{\mu}+\sqrt{\gamma_{0}}}\to-1$ as $\gamma_0\to\infty$ and $\frac{e^{\sqrt{\mu}t}+1}{e^{\sqrt{\mu}t}-1}=\coth(\frac{\sqrt{\mu}}{2}t)$, the ODE system converges to the unified NAG system \eqref{eq:u_ode} as $\gamma_0\to\infty$.

\section{Higher-Order Hyperbolic Functions}
\label{app:hyperbolic}

\subsection{Proof of Proposition~\ref{prop:sinhp_grows_exponentially}}
\label{app:prop:sinhp_grows_exponentially}

Fix $T>0$. We will show that
\begin{equation}
	\label{eq:to_converge_to_Cp}
	\log\left(\sinh_{p}(T+t)\right)-t
\end{equation}
converges to some constant as $t\rightarrow\infty$. We can bound the derivative of \eqref{eq:to_converge_to_Cp} as
\begin{align*}
	\frac{d}{dt}\left\{ \log\left(\sinh_{p}(T+t)\right)-t\right\}  & =\frac{\sinh_{p}'(T+t)}{\sinh_{p}(T+t)}-1\\
	& =\frac{\cosh_{p}(T+t)}{\sinh_{p}(T+t)}-1\\
	& =\left(1+\frac{1}{\sinh_{p}^{p}(T+t)}\right)^{1/p}-1\\
	& \in\left[0,\frac{1}{\sinh_{p}(T+t)}\right],
\end{align*}
where the last line follows from the fact that $1\leq(1+x)^{1/p}\leq1+x^{1/p}$ holds for $x\geq0$.\footnote{To check this basic inequality, one can consider  the $p$-th power of each side.} Thus, if the integral
\begin{equation}
	\label{eq:integral_to_be_finite}
	\int_{0}^{\infty}\frac{1}{\sinh_{p}(T+t)}\,dt
\end{equation}
is finite, then \eqref{eq:to_converge_to_Cp} converges to some constant because it is monotonically increasing and bounded above, and thus this completes the proof. To show that the integral \eqref{eq:integral_to_be_finite} is finite, it is enough to show that the inequality
\[
\sinh_{p}(T+t)\geq\sinh_{p}(T)e^{t}
\]
holds for all $t\geq0$. This can be shown by the following calculation:
\begin{align*}
	\log\left(\sinh_{p}(T+t)\right) & =\log\left(\sinh_{p}(T)\right)+\int_{0}^{t}\frac{d}{ds}\left\{ \log\left(\sinh_{p}(T+s)\right)\right\} \,ds\\
	& =\log\left(\sinh_{p}(T)\right)+\int_{0}^{t}\frac{\sinh_{p}'(T+s)}{\sinh_{p}(T+s)}\,ds\\
	& =\log\left(\sinh_{p}(T)\right)+\int_{0}^{t}\frac{\left(1+\sinh_{p}^{p}(T+s)\right)^{1/p}}{\sinh_{p}(T+s)}\,ds\\
	& \geq\log\left(\sinh_{p}(T)\right)+\int_{0}^{t}1\,ds\\
	& =\log\left(\sinh_{p}(T)\right)+t\\
	& =\log\left(\sinh_{p}(T)e^{t}\right).
\end{align*}

\subsection{The function $\sinhc_p$ is non-decreasing}
\label{app:sinhc_is_increasing}

It is easy to see that $\sinh_p$ and $\cosh_p$ are increasing. Since
\begin{align*}
	\tanh_{p}'(t) & =\frac{d}{dt}\left\{ \frac{\sinh_{p}(t)}{\cosh_{p}(t)}\right\} \\
	& =\frac{\sinh_{p}'(t)\cosh_{p}(t)-\cosh_{p}'(t)\sinh_{p}(t)}{\cosh_{p}^{2}(t)}\\
	& \leq\frac{\sinh_{p}'(t)\cosh_{p}(t)}{\cosh_{p}^{2}(t)}\\
	& =1,
\end{align*}
we have $\tanh_p(t)\leq t$ for all $t\geq0$. Now, we deduce that
\begin{align*}
	\sinhc_{p}'(t) & =\frac{d}{dt}\left\{ \frac{\sinh_{p}(t)}{t}\right\} \\
	& =\frac{t\sinh_{p}'(t)-\sinh_{p}(t)}{t^{2}}\\
	& =\frac{t\cosh_{p}(t)-\sinh_{p}(t)}{t^{2}}\\
	& =\frac{\cosh_{p}(t)}{t^{2}}\left(t-\tanh_{p}(t)\right)\\
	& \geq0,
\end{align*}
and thus $\sinhc$ is non-decreasing.

\section{Limiting Arguments}
\label{app:limiting_arguments}

\subsection{Limiting argument for two-sequence scheme}
\label{app:two-sequence_scheme}

\paragraph{Limiting ODE of two-sequence scheme.}

For the iterates of the two-sequence scheme \eqref{eq:two_sequence_scheme}, we have
\begin{align*}
	\frac{x_{k+1}-x_{k}}{\sqrt{s}} & =\frac{1}{\sqrt{s}}\left(y_{k}-s\nabla f\left(y_{k}\right)-x_{k}\right)\\
	& =\frac{1}{\sqrt{s}}\left(\beta_{k-1}\left(x_{k}-x_{k-1}\right)+\gamma_{k}\left(x_{k}-y_{k-1}\right)-s\nabla f\left(y_{k}\right)\right)\\
	& =\frac{1}{\sqrt{s}}\left(\beta_{k-1}\left(x_{k}-x_{k-1}\right)-s\gamma_{k}\nabla f\left(y_{k-1}\right)-s\nabla f\left(y_{k}\right)\right)\\
	& =\beta_{k-1}\frac{x_{k}-x_{k-1}}{\sqrt{s}}-\sqrt{s}\gamma_{k}\nabla f\left(y_{k-1}\right)-\sqrt{s}\nabla f\left(y_{k}\right).
\end{align*}
Using the Taylor expansions
\begin{align*}
	\frac{x_{k+1}-x_{k}}{\sqrt{s}} & =\dot{X}(t_k)+\frac{1}{2}\ddot{X}(t_k)\sqrt{s}+o\left(\sqrt{s}\right)\\
	\frac{x_{k}-x_{k-1}}{\sqrt{s}} & =\dot{X}(t_k)-\frac{1}{2}\ddot{X}(t_k)\sqrt{s}+o\left(\sqrt{s}\right),
\end{align*}
we obtain
\[
\dot{X}(t_{k})+\frac{1}{2}\ddot{X}(t_{k})\sqrt{s}+o\left(\sqrt{s}\right)=\beta_{k-1}\left(\dot{X}(t_{k})-\frac{1}{2}\ddot{X}(t_{k})\sqrt{s}+o\left(\sqrt{s}\right)\right)-\sqrt{s}\gamma_{k}\nabla f\left(y_{k-1}\right)-\sqrt{s}\nabla f\left(y_{k}\right).
\]
It follows from $\Vert x_{k}-y_{k-1}\Vert=o(\sqrt{s})$ and the Lipschitz continuity of $\nabla f$ that
\begin{align*}
	\sqrt{s}\nabla f\left(y_{k-1}\right) & =\sqrt{s}\nabla f(X(t_{k}))+o\left(\sqrt{s}\right)\\
	\sqrt{s}\nabla f\left(y_{k}\right) & =\sqrt{s}\nabla f\left(y_{k-1}\right)+o\left(\sqrt{s}\right)=\sqrt{s}\nabla f(X(t_{k}))+o\left(\sqrt{s}\right).
\end{align*}
Substituting these into the ODE yields
\[
\frac{1+\beta_{k-1}}{2}\ddot{X}(t_{k})\sqrt{s}+\left(1-\beta_{k-1}\right)\dot{X}(t_{k})+\left(1+\gamma_{k}\right)\nabla f(X(t_{k}))\sqrt{s}+o\left(\sqrt{s}\right)=0.
\]
Dividing both sides by $\sqrt{s}$, substituting $k=t/\sqrt{s}$ and the limits \eqref{eq:two_seq_assumptions}, and then letting $s\rightarrow0$, we obtain (note that $\beta_{t/\sqrt{s}-1}\rightarrow1$ by Equation~\eqref{eq:two_seq_assumptions})
\[
\ddot{X}(t)+b(t)\dot{X}(t)+(1+c(t))\nabla f(X(t))=0.
\]

\paragraph{Recovering the limiting ODE of three-sequence scheme.}

It follows from the Taylor expansion that
\begin{align*}
	\tau_{k} & =\tau\left(t_{k}\right)\sqrt{s}\\
	\tau_{k+1} & =\tau\left(t_{k}\right)\sqrt{s}+\dot{\tau}\left(t_{k}\right)s+\sqrt{s}o\left(\sqrt{s}\right)\\
	\delta_{k} & =\delta\left(t_{k}\right)\sqrt{s}.
\end{align*}
Thus, for the sequences $(\beta_k)$ and $(\gamma_k)$ in \eqref{eq:beta_three}, we have
\begin{align*}
	\frac{1-\beta_{k}}{\sqrt{s}} & =\frac{1}{\sqrt{s}}\left(1-\left(1-\tau_{k}\right)\left(1-\mu\delta_{k}\right)\frac{\tau_{k+1}}{\tau_{k}}\right)\\
	& =\frac{1}{\sqrt{s}}\left(1-\left(1-\sqrt{s}\tau\left(t_{k}\right)\right)\left(1-\mu\sqrt{s}\delta\left(t_{k}\right)\right)\left(1+\frac{\dot{\tau}\left(t_{k}\right)s+\sqrt{s}o\left(\sqrt{s}\right)}{\tau\left(t_{k}\right)\sqrt{s}}\right)\right)\\
	& =\frac{1}{\sqrt{s}}\left(\sqrt{s}\tau\left(t_{k}\right)+\mu\sqrt{s}\delta\left(t_{k}\right)-\sqrt{s}\frac{\dot{\tau}\left(t_{k}\right)}{\tau\left(t_{k}\right)}+o\left(\sqrt{s}\right)\right)\\
	& =\tau\left(t_{k}\right)+\mu\delta\left(t_{k}\right)-\frac{\dot{\tau}\left(t_{k}\right)}{\tau\left(t_{k}\right)}+\frac{o\left(\sqrt{s}\right)}{\sqrt{s}}
\end{align*}
and
\begin{align*}
	\gamma_{k} & =\frac{\tau_{k+1}}{\tau_{k}}\left((1/s-\mu)\delta_{k}\tau_{k}-1+\mu\delta_{k}\right)\\
	& =\left(1+\frac{\dot{\tau}\left(t_{k}\right)\sqrt{s}+o\left(\sqrt{s}\right)}{\tau\left(t_{k}\right)}\right)\left((1-\mu s)\delta\left(t_{k}\right)\tau\left(t_{k}\right)-1+\mu\sqrt{s}\delta\left(t_{k}\right)\right)\\
	& =\delta\left(t_{k}\right)\tau\left(t_{k}\right)-1+o\left(1\right).
\end{align*}
Therefore, we have
\begin{align*}
	\lim_{s\rightarrow0}\frac{1-\beta_{t/\sqrt{s}}}{\sqrt{s}} & =\tau(t)+\mu\delta(t)-\frac{\dot{\tau}(t)}{\tau(t)}\\
	\lim_{s\rightarrow0}\gamma_{t/\sqrt{s}} & =\tau(t)\delta(t)-1,
\end{align*}
which recovers the limiting ODE~\eqref{eq:limiting_ode_oneline} of the three-sequence scheme.

\subsection{Difference matrix and differential kernel}
\label{app:differential_kernel}

\paragraph{From the two-sequence scheme to the difference matrix.}

The iterates of the two-sequence scheme \eqref{eq:two_sequence_scheme} satisfy
\begin{align*}
	y_{k+1}-y_{k} & =x_{k+1}-y_{k}+\beta_{k}\left(x_{k+1}-x_{k}\right)-s\gamma_{k}\nabla f\left(y_{k}\right)\\
	& =\beta_{k}\left(y_{k}-y_{k-1}\right)+s\beta_{k}\nabla f\left(y_{k-1}\right)-s\left(1+\beta_{k}+\gamma_{k}\right)\nabla f\left(y_{k}\right).
\end{align*}
Substituting
\begin{align*}
	y_{k+1}-y_{k} & =-s\sum_{i=0}^{k}h_{k,i}\nabla f\left(y_{i}\right)\\
	y_{k}-y_{k-1} & =-s\sum_{i=0}^{k-1}h_{k-1,i}\nabla f\left(y_{i}\right)
\end{align*}
into the equality and 
comparing the coefficients of each $\nabla f(y_{i})$, we obtain
\[
h_{k,j}=\begin{cases}
	1+\beta_{k}+\gamma_{k}, & \textrm{if }j=k\\
	\beta_{k}\left(h_{k-1,k-1}-1\right), & \textrm{if }j=k-1\\
	\beta_{k}h_{k-1,i}, & \textrm{if }j\leq k-2.
\end{cases}
\]
Using mathematical induction, it is straightforward to show that
\[
h_{ij}=\left(\beta_{j}+\gamma_{j}\right)\prod_{\nu=j+1}^{i}\beta_{\nu}+\delta_{ij}.
\]

\paragraph{Differential kernel for the two-sequence scheme.}
By \eqref{eq:partial_derivative}, we have
\[
\frac{\partial}{\partial s}\log(H(s,t))=\frac{\partial H(s,\tau)}{\partial s}\frac{1}{H(s,\tau)}=-b(s).
\]
Integrating over $s$, we obtain
\[
\log\left(H(t,\tau)\right)-\log\left(H(\tau,\tau)\right)=-\int_{\tau}^{t}b(s)\,ds.
\]
Thus, we have
\[
H(t,\tau)=H(\tau,\tau)e^{-\int_{\tau}^{t}b(s)\,ds}=\left(1+c(\tau)\right)e^{-\int_{\tau}^{t}b(s)\,ds}.
\]

\section{Unified Bregman Lagrangian}
\label{app:proofs}

\subsection{Proof of Proposition~\ref{prop:e-l-eq}}
\label{app:prop:e-l-eq}

For the unified Bregman Lagrangian \eqref{eq:unified_lagrangian}, the partial derivatives $\frac{\partial\mathcal{L}}{\partial \dot{X}}\left(X,\dot{X},t\right)$ and $\frac{\partial\mathcal{L}}{\partial X}\left(X,\dot{X},t\right)$ are given by
\begin{align*}
	\frac{\partial\mathcal{L}}{\partial \dot{X}}\left(X,\dot{X},t\right) & =e^{\gamma}\left(1+\mu e^{\beta}\right)\left(\nabla h\left(X+e^{-\alpha}\dot{X}\right)-\nabla h(X)\right)\\
	\frac{\partial\mathcal{L}}{\partial X}\left(X,\dot{X},t\right) & =e^{\alpha+\gamma}\left(1+\mu e^{\beta}\right)\left(\nabla h\left(X+e^{-\alpha}\dot{X}\right)-\nabla h(X)\right)\\
	& \quad-e^{\gamma}\left(1+\mu e^{\beta}\right)\frac{d}{dt}\nabla h(X)-e^{\alpha+\beta+\gamma}\nabla f(X).
\end{align*}
The time derivative of $\frac{\partial\mathcal{L}}{\partial \dot{X}}$ can be computed as
\begin{align*}
	\frac{d}{dt}\left\{ \frac{\partial\mathcal{L}}{\partial \dot{X}}\left(X,\dot{X},t\right)\right\}  & =\left(\dot{\gamma}e^{\gamma}+\mu\left(\dot{\beta}+\dot{\gamma}\right)e^{\beta+\gamma}\right)\left(\nabla h\left(X+e^{-\alpha}\dot{X}\right)-\nabla h(X)\right)\\
	& \quad+e^{\gamma}\left(1+\mu e^{\beta}\right)\left(\frac{d}{dt}\nabla h\left(X+e^{-\alpha}\dot{X}\right)-\frac{d}{dt}\nabla h(X)\right).
\end{align*}
Thus, the Euler--Lagrange equation \eqref{eq:e-l-eq} can be written as
\begin{multline*}
	e^{\gamma}\left(1+\mu e^{\beta}\right)\frac{d}{dt}\nabla h\left(X+e^{-\alpha}\dot{X}\right)=\left(e^{\alpha+\gamma}\left(1+\mu e^{\beta}\right)-\dot{\gamma}e^{\gamma}-\mu\left(\dot{\beta}+\dot{\gamma}\right)e^{\beta+\gamma}\right)
	\times \left(\nabla h\left(X+e^{-\alpha}\dot{X}\right)-\nabla h(X)\right)-e^{\alpha+\beta+\gamma}\nabla f(X).
\end{multline*}
Substituting $\dot{\gamma}=e^{\alpha}$ \eqref{eq:ideal_scaling_a} into the equation and dividing both sides by $e^{\gamma}\left(1+\mu e^{\beta}\right)>0$, we obtain
\[
\frac{d}{dt}\nabla h\left(X+e^{-\alpha}\dot{X}\right)=-\frac{\mu\dot{\beta}e^{\beta}}{1+\mu e^{\beta}}\left(\nabla h\left(X+e^{-\alpha}\dot{X}\right)-\nabla h(X)\right)-\frac{e^{\alpha+\beta}}{1+\mu e^{\beta}}\nabla f(X).
\]
Letting $Z=X+e^{-\alpha}\dot{X}$ yields the system of ODEs \eqref{eq:unified_family}.

\subsection{Proof of Theorem~\ref{thm:mainthm_lagrangian}}
\label{app:thm:mainthm_lagrangian}

Note that
\begin{align*}
	\frac{d}{dt}D_{h}\left(x^{*},Z\right) & =\frac{d}{dt}\left\{ h\left(x^{*}\right)-h(Z)-\left\langle \nabla h(Z),x^{*}-Z\right\rangle \right\} \\
	& =-\left\langle \nabla h(Z),\dot{Z}\right\rangle -\left\langle \frac{d}{dt}\nabla h(Z),x^{*}-Z\right\rangle +\left\langle \nabla h(Z),\dot{Z}\right\rangle \\
	& =-\left\langle \frac{d}{dt}\nabla h(Z),x^{*}-Z\right\rangle .
\end{align*}
Using this equation, we have
\begin{align*}
	\frac{d}{dt}\left\{ \phi(X(t),Z(t),t)\right\} & =-\left(1+\mu e^{\beta}\right)\left\langle \frac{d}{dt}\nabla h(Z),x^{*}-Z\right\rangle +\mu\dot{\beta}e^{\beta}D_{h}\left(x^{*},Z\right)\\
	& \quad+\dot{\beta}e^{\beta}\left(f(X)-f\left(x^{*}\right)\right)+e^{\beta}\left\langle \nabla f(X),\dot{X}\right\rangle \\
	& =\left\langle \mu\dot{\beta}e^{\beta}\left(\nabla h(Z)-\nabla h(X)\right)+e^{\alpha+\beta}\nabla f(X),x^{*}-Z\right\rangle +\mu\dot{\beta}e^{\beta}D_{h}\left(x^{*},Z\right)\\
	& \quad+\dot{\beta}e^{\beta}\left(f(X)-f\left(x^{*}\right)\right)+e^{\beta}\left\langle \nabla f(X),\dot{X}\right\rangle ,
\end{align*}
where the second equality follows from \eqref{eq:unified_family_b}. It follows from the Bregman three-point identity \eqref{eq:three-point}, the non-negativity of Bregman divergence,
and the $\mu$-uniform convexity of $f$ with respect to $h$ \eqref{eq:uniform_convexity} that
\begin{align*}
	\left\langle \nabla h(Z)-\nabla h(X),x^{*}-Z\right\rangle +D_{h}\left(x^{*},Z\right) & =D_{h}\left(x^{*},X\right)-D_{h}(Z,X)\\
	& \leq D_{h}\left(x^{*},X\right)\\
	& \leq\frac{1}{\mu}D_{f}\left(x^{*},X\right).
\end{align*}
Thus, we have
\begin{align*}
	\frac{d}{dt}\left\{ \phi(X(t),Z(t),t)\right\} & \leq\dot{\beta}e^{\beta}D_{f}\left(x^{*},X\right)+e^{\alpha+\beta}\left\langle \nabla f(X),x^{*}-Z\right\rangle \\
	& \quad+\dot{\beta}e^{\beta}\left(f(X)-f\left(x^{*}\right)\right)+e^{\beta}\left\langle \nabla f(X),\dot{X}\right\rangle \\
	& =\dot{\beta}e^{\beta}D_{f}\left(x^{*},X\right)+e^{\alpha+\beta}\left\langle \nabla f(X),x^{*}-X\right\rangle +\dot{\beta}e^{\beta}\left(f(X)-f\left(x^{*}\right)\right)\\
	& =\left(e^{\alpha}-\dot{\beta}\right)e^{\beta}\left\langle \nabla f(X),x^{*}-X\right\rangle \\
	& \leq\left(e^{\alpha}-\dot{\beta}\right)e^{\beta}\left(f\left(x^{*}\right)-f(X)\right)\\
	& \leq0,
\end{align*}
where the last two inequalities follows from the ideal scaling condition \eqref{eq:ideal_scaling_b}, the convexity of $f$, and the fact that $x^*$ is a minimizer of $f$.

\subsection{Proof of Theorem~\ref{thm:time_dilation}}
\label{app:thm:time_dilation}

The derivatives of $X_2$ and $\nabla h(Z_{2})$ can be computed as
\begin{align*}
	\dot{X}_{2}(t) & =\dot{\bfT}(t)\dot{X}_{1}(\bfT(t))\\
	& =\dot{\bfT}(t)e^{\alpha_{1}(\bfT(t))}(Z_{1}(\bfT(t))-X_{1}(\bfT(t))\\
	& =\dot{\bfT}(t)e^{\alpha_{1}(\bfT(t))}(Z_{2}(t)-X_{2}(t))\\
	& =e^{\alpha_{2}(t)}(Z_{2}(t)-X_{2}(t))
\end{align*}
and
\begin{align*}
	\frac{d}{dt}\nabla h(Z_{2}(t)) & =\dot{\bfT}(t)\frac{d(\nabla h\circ Z_{1})}{dt}(\bfT(t))\\
	& =\dot{\bfT}(t)\Bigg(\frac{\mu\dot{\beta}_{1}(\bfT(t))e^{\beta_{1}(\bfT(t))}}{1+\mu e^{\beta_{1}(\bfT(t))}}\left(\nabla h(X_{1}(\bfT(t))-\nabla h(Z_{1}(\bfT(t)))\right)\\
	& \qquad-\frac{e^{\alpha_{1}(\bfT(t))+\beta_{1}(\bfT(t))}}{1+\mu e^{\beta_{1}(\bfT(t))}}\nabla f(X_{1}(\bfT(t)))\Bigg)\\
	& =\frac{\mu\dot{\beta}_{2}(t)e^{\beta_{2}(t)}}{1+\mu e^{\beta_{2}(t)}}\left(\nabla h(X_{2}(t))-\nabla h(Z_{2}(t))\right)-\frac{e^{\alpha_{2}(t)+\beta_{2}(t)}}{1+\mu e^{\beta_{2}(t)}}\nabla f(X_{2}(t)).
\end{align*}
Thus, we obtain the desired system of ODEs.

\subsection{Recovering Lyapunov analysis for the second Bregman Lagrangian flow}
\label{app:recover_lagrangian}

In this section, we recover the second Bregman Lagrangian flow \eqref{eq:second_family} with constant coefficients and its Lyapunov analysis from the unified Bregman Lagrangian flow \eqref{eq:unified_family} and its Lyapunov analysis (Theorem~\ref{thm:mainthm_lagrangian}). In particular, we recover \nagscsp ODE \eqref{eq:sc_ode} and its Lyapunov analysis from the unified NAG ODE \eqref{eq:u_ode_oneline} and its Lyapunov analysis (Theorem~\ref{thm:mainthm_continuous}).

For the parameter functions $\alpha,\beta:[0,\infty)\rightarrow\mathbb{R}$ of the unified Bregman Lagrangian flow \eqref{eq:unified_family}, assume that the limits $\alpha(\infty):=\lim_{t\rightarrow\infty}\alpha(t)$ and $\dot{\beta}(\infty):=\lim_{t\rightarrow\infty}\dot{\beta}(t)>0$ exist. We consider the following second Bregman Lagrangian flow \eqref{eq:second_family} with ${\alpha}_{\mathrm{2nd}}(t):\equiv\alpha(\infty)$ and ${\beta}_{\mathrm{2nd}}(t):=\dot{\beta}(\infty)t$:
\begin{equation}
	\label{eq:limiting_of_unified_lagrangian}
	\begin{aligned}
		\dot{X} & =e^{\alpha(\infty)}(Z-X)\\
		\frac{d}{dt}\nabla h(Z) & =\dot{\beta}(\infty)\left(\nabla h(X)-\nabla h(Z)\right)-\frac{e^{\alpha(\infty)}}{\mu}\nabla f(X).
	\end{aligned}
\end{equation}
Then, it follows from $\lim_{t\rightarrow\infty}e^{\alpha(t)}=e^{\alpha(\infty)}$, $\lim_{t\rightarrow\infty}\frac{\mu\dot{\beta}e^{\beta}}{1+\mu e^{\beta}}=\dot{\beta}(\infty)$, and $\lim_{t\rightarrow\infty}\frac{e^{\alpha+\beta}}{1+\mu e^{\beta}}=\frac{e^{\alpha(\infty)}}{\mu}$ that the coefficients in the unified Bregman Lagrangian flow \eqref{eq:unified_family} converge to those in the dynamics \eqref{eq:limiting_of_unified_lagrangian} as $t\to\infty$. Thus, roughly speaking, the dynamics \eqref{eq:limiting_of_unified_lagrangian} is the \emph{asymptotic version} of the unified Bregman Lagrangian flow in the sense that [the flow corresponding to \eqref{eq:unified_family}, starting at time $t_0$] converges to [the flow corresponding to \eqref{eq:limiting_of_unified_lagrangian}, starting at time $0$] as $t_0\to\infty$.

Note that the time derivative of the Lyapunov function \eqref{eq:lyapunov_lagrangian} for the unified Bregman Lagrangian flow can be written as
\begin{align*}
	\frac{d}{dt}\left\{ V(X(t),Z(t),t)\right\}  & =\frac{d}{dt}\left\{ 1+\mu e^{\beta}\right\} D_{h}\left(x^{*},Z\right)+\left(1+\mu e^{\beta}\right)\frac{d}{dt}\left\{ D_{h}\left(x^{*},Z\right)\right\} \\
	& \quad+\frac{d}{dt}\left\{ e^{\beta}\right\} \left(f(X)-f\left(x^{*}\right)\right)+e^{\beta}\frac{d}{dt}\left\{ f(X)-f\left(x^{*}\right)\right\} .
\end{align*}
Thus, we have
\begin{align*}
	0 & \geq e^{-\beta(t_{0}+t)}\frac{d}{dt}\left\{ V(X(t_{0}+t),Z(t_{0}+t),t_{0}+t)\right\} \\
	& =\mu\dot{\beta}(t_{0}+t)D_{h}\left(x^{*},Z(t_{0}+t)\right)+\frac{1+\mu e^{\beta(t_{0}+t)}}{e^{\beta(t_{0}+t)}}\frac{d}{dt}\left\{ D_{h}\left(x^{*},Z(t_{0}+t)\right)\right\} \\
	& \quad+\dot{\beta}(t_{0}+t)\left(f(X(t_{0}+t))-f\left(x^{*}\right)\right)+\frac{d}{dt}\left\{ f(X(t_{0}+t))-f\left(x^{*}\right)\right\} 
\end{align*}
for all $t>0$, where $t_0>0$ is the initial time of the flow. Fix $x_0=X(t_0)$ and $z_0=Z(t_0)$ in $\mathbb{R}^n$. Note that as $t_0\to\infty$, the flow $t\mapsto (X(t_0+t),Z(t_0+t))$ converges to the flow $t\mapsto (X_{\mathrm{2nd}}(t),Z_{\mathrm{2nd}}(t))$ corresponding to \eqref{eq:limiting_of_unified_lagrangian} with $X_{\mathrm{2nd}}(0)=x_0$ and $Z_{\mathrm{2nd}}(0)=z_0$. Now, taking the limit $t_0\to\infty$ in the inequality above yields
yields
\begin{align*}
	0 & \geq\mu\dot{\beta}(\infty)D_{h}\left(x^{*},Z(t)\right)+\mu\frac{d}{dt}\left\{ D_{h}\left(x^{*},Z(t)\right)\right\} \\
	& \quad+\dot{\beta}(\infty)\left(f(X(t))-f\left(x^{*}\right)\right)+\frac{d}{dt}\left\{ f(X(t))-f\left(x^{*}\right)\right\} \\
	& =e^{-{\beta}_{\mathrm{2nd}}(t)}\frac{d}{dt}\left\{ {V}_{\mathrm{2nd}}(X(t),Z(t),t)\right\} ,
\end{align*}
where ${V}_{\mathrm{2nd}}$ is the Lyapunov function \eqref{eq:second_lagrangian_lyapunov} for the second Bregman Lagrangian flow with the parameters ${\alpha}_{\mathrm{2nd}}$ and ${\beta}_{\mathrm{2nd}}$. Because $e^{-{\beta}_{\mathrm{2nd}}(t)}>0$, we recover the Lyapunov analysis for the second Bregman Lagrangian flow.

\paragraph{Recovering \nagscsp ODE from the unified ODE.}

Note that the unified Bregman Lagrangian flow \eqref{eq:unified_family} and its Lyapunov analysis (Theorem~\ref{thm:mainthm_lagrangian}) with $h(x)=\frac{1}{2}\left\Vert x\right\Vert ^{2}$, $\alpha(t)=\log\left(\frac{2}{t}\cothc\left(\frac{\sqrt{\mu}}{2}t\right)\right)$, and $\beta(t)=\log\left(\frac{t^{2}}{4}\sinhc^{2}\left(\frac{\sqrt{\mu}}{2}t\right)\right)$ recover the unified NAG system \eqref{eq:u_ode} and its Lyapunov analysis (Theorem~\ref{thm:mainthm_continuous}). Also, note that the second Bregman Lagrangian flow \eqref{eq:second_family} and the corresponding Lyapunov function \eqref{eq:second_lagrangian_lyapunov} with ${\alpha}_{\mathrm{2nd}}(t)=\log\left(\sqrt{\mu}\right)$ and ${\beta}_{\mathrm{2nd}}(t)=\sqrt{\mu}t$ recover \nagscsp system \eqref{eq:sc_ode_twoline} and the corresponding Lyapunov function \eqref{eq:sc_ode_lyapunov}. When $\mu>0$, because $\alpha(\infty)=\log\left(\sqrt{\mu}\right)$ and $\dot{\beta}(\infty)=\sqrt{\mu}$, the results above shows that \nagscsp ODE is the asymptotic version of the unified NAG ODE and that the Lyapunov analysis of \nagscsp ODE can be obtained by taking the limit $t\rightarrow\infty$ into the coeffiicients of the inequality (rigorously, taking the limit $t_0\rightarrow\infty$ of the initial time as in the preceding paragraph)
\[
\frac{4}{t^{2}}\cschc^{2}\left(\frac{\sqrt{\mu}}{2}t\right)\frac{d}{dt}\left\{ V(X(t),Z(t),t)\right\} \leq0,
\]}
where $V$ is the Lyapunov function \eqref{eq:lyapunov_continuous} for the unified NAG ODE.

\section{Unified NAG ODE}

\subsection{Choosing $\alpha$ and $\beta$}
\label{app:alphabeta}

We first note some properties of the functions $\alpha$ and $\beta$ that recover \nagcsp ODE (or \nagscsp ODE) from the first Bregman Lagrangian flow (or the second Bregman Lagrangian flow, respectively).

The first Bregman Lagrangian flow \eqref{eq:first_family} with $h(x)=\frac{1}{2}\left\Vert x\right\Vert ^{2}$ can be written as the following ODE:
\[
\ddot{X}+\left(-\dot{\alpha}+e^{\alpha}\right)\dot{X}+e^{2\alpha+\beta}\nabla f(X)=0.
\]
The choices $\alpha(t)=\log\frac{2}{t}$ and $\beta(t)=\log\frac{t^{2}}{4}$, which recover \nagcsp ODE, satisfy the ideal scaling condition \eqref{eq:ideal_scaling_b} with equality and make the coefficient of $\nabla f(X)$ equal to the coefficient of $\ddot{X}$. 

The second Bregman Lagrangian flow \eqref{eq:second_family} with $h(x)=\frac{1}{2}\left\Vert x\right\Vert ^{2}$ can be written as
\[
\ddot{X}+\left(-\dot{\alpha}+e^{\alpha}+\dot{\beta}\right)\dot{X}+\frac{e^{2\alpha}}{\mu}\nabla f(X)=0.
\]
The choices $\alpha(t)=\log\sqrt{\mu}$ and $\beta(t)=\log\left(\sqrt{\mu}t\right)$, which recover \nagscsp ODE, satisfy the ideal scaling condition \eqref{eq:ideal_scaling_b} with equality and make the coefficient of $\nabla f(X)$ equal to the coefficient of $\ddot{X}$. 

Inspired by these facts, for the unified Bregman Lagrangian, we
construct functions $\alpha(t)$ and $\beta(t)$ so that the ideal scaling condition \eqref{eq:ideal_scaling_b} holds with equality and that the coefficient
of $\nabla f(X)$ is equal to the coefficient of $\ddot{X}$. The unified Bregman Lagrangian flow \eqref{eq:unified_family} with $h(x)=\frac{1}{2}\left\Vert x\right\Vert ^{2}$ can be written as
\[
\ddot{X}+\left(-\dot{\alpha}+e^{\alpha}+\frac{\mu\dot{\beta}e^{\beta}}{1+\mu e^{\beta}}\right)\dot{X}+\frac{e^{2\alpha+\beta}}{1+\mu e^{\beta}}\nabla f(X)=0.
\]
Now, we solve the following system of ODEs:
\begin{align*}
	\dot{\beta} & =e^{\alpha}\\
	e^{2\alpha+\beta} & =1+\mu e^{\beta}.
\end{align*}
Let $A(t)=e^{\beta(t)}>0$. Then, we have $\dot{A}=\dot{\beta}e^{\beta}=e^{\alpha+\beta}>0$. Because $(\dot{A})^{2}=e^{2\alpha+\beta}e^{\beta}=A(1+\mu A)$, we have  $\dot{A}=\sqrt{A(1+\mu A)}$. Solving this differential equation with the initial condition $A(0)=0$ yields $A=\frac{t^{2}}{4}\sinhc^{2}(\frac{\sqrt{\mu}}{2}t)$. Thus, we have $\beta(t)=\log(\frac{t^{2}}{4}\sinhc^{2}(\frac{\sqrt{\mu}}{2}t))$ and $\alpha(t)=\log(\dot{\beta}(t))=\log(\frac{2}{t}\cothc(\frac{\sqrt{\mu}}{2}t))$.

\subsection{Equivalent forms of the unified NAG system and the unified NAG-G system}
\label{app:equivalent_of_ode}

When $\mu=0$, the unified NAG system is equivalent to \nagcsp system. Thus, we assume $\mu>0$ for the sake of simplicity. 

\paragraph{Second-order ODE form of the unified NAG system.}
When $\mu>0$, we can write the unified NAG system \eqref{eq:u_ode} as
\begin{align*}
	\dot{X} & =\sqrt{\mu}\coth\left(\frac{\sqrt{\mu}}{2}t\right)(Z-X)\\
	\dot{Z} & =\frac{1}{\sqrt{\mu}}\tanh\left(\frac{\sqrt{\mu}}{2}t\right)\left(\mu X-\mu Z-\nabla f(X)\right).
\end{align*}
Substituting $Z=X+\frac{1}{\sqrt{\mu}}\tanh(\frac{\sqrt{\mu}}{2}t)\dot{X}$ into $\dot{Z}=\frac{1}{\sqrt{\mu}}\tanh(\frac{\sqrt{\mu}}{2}t)(\mu X-\mu Z-\nabla f(X))$, we have
\begin{align*}
	& \frac{1}{\sqrt{\mu}}\tanh\left(\frac{\sqrt{\mu}}{2}t\right)\ddot{X}+\left(1+\frac{1}{2}\sech^{2}\left(\frac{\sqrt{\mu}}{2}t\right)\right) \\
	& =\frac{1}{\sqrt{\mu}}\tanh\left(\frac{\sqrt{\mu}}{2}t\right)\left(\mu X-\mu Z-\nabla f(X)\right)\\
	& =-\sqrt{\mu}\tanh\left(\frac{\sqrt{\mu}}{2}t\right)\left(Z-X\right)-\frac{1}{\sqrt{\mu}}\tanh\left(\frac{\sqrt{\mu}}{2}t\right)\nabla f(X)\\
	& =-\tanh^{2}\left(\frac{\sqrt{\mu}}{2}t\right)\dot{X}-\frac{1}{\sqrt{\mu}}\tanh\left(\frac{\sqrt{\mu}}{2}t\right)\nabla f(X).
\end{align*}
Multiplying by $\sqrt{\mu}\coth(\frac{\sqrt{\mu}}{2}t)$ and rearranging the terms, we have
\begin{multline*}
	\ddot{X}+\bigg(\sqrt{\mu}\tanh\left(\frac{\sqrt{\mu}}{2}t\right)+\sqrt{\mu}\coth\left(\frac{\sqrt{\mu}}{2}t\right)\\
	+\frac{\sqrt{\mu}}{2}\sech\left(\frac{\sqrt{\mu}}{2}t\right)\csch\left(\frac{\sqrt{\mu}}{2}t\right)\bigg)\dot{X}+\nabla f(X)=0.
\end{multline*}
Using the identity $\tanh(x)-\coth(x)+\sech(x)\csch(x)=0$, we can equivalently write this ODE as
\[
\ddot{X}+\left(\frac{\sqrt{\mu}}{2}\tanh\left(\frac{\sqrt{\mu}}{2}t\right)+\frac{3\sqrt{\mu}}{2}\coth\left(\frac{\sqrt{\mu}}{2}t\right)\right)\dot{X}+\nabla f(X)=0.
\]

\paragraph{Differential kernel for the unified NAG ODE.}

Substituting $b(t)=\frac{\sqrt{\mu}}{2}\tanh(\frac{\sqrt{\mu}}{2}t)+\frac{3\sqrt{\mu}}{2}\coth(\frac{\sqrt{\mu}}{2}t)$ and $c(t)=0$ into \eqref{eq:kernel_for_two}, we yield the following differential kernel corresponding to the unified NAG ODE:
\begin{align*}
	H(t,\tau) & =e^{-\int_{\tau}^{t}\left(\frac{\sqrt{\mu}}{2}\tanh\left(\frac{\sqrt{\mu}}{2}s\right)+\frac{3\sqrt{\mu}}{2}\coth\left(\frac{\sqrt{\mu}}{2}s\right)\right)\,ds}\\
	& =e^{-\left[3\log\left(\sinh\left(\frac{\sqrt{\mu}}{2}s\right)\right)+\log\left(\cosh\left(\frac{\sqrt{\mu}}{2}s\right)\right)\right]_{\tau}^{t}}\\
	& =\frac{\sinh^{3}\left(\frac{\sqrt{\mu}}{2}\tau\right)\cosh\left(\frac{\sqrt{\mu}}{2}\tau\right)}{\sinh^{3}\left(\frac{\sqrt{\mu}}{2}t\right)\cosh\left(\frac{\sqrt{\mu}}{2}t\right)}.
\end{align*}

\paragraph{Differential kernel for the unified NAG-G ODE.}

Substituting $b(t)=\frac{\sqrt{\mu}}{2}\tanh(\frac{\sqrt{\mu}}{2}(T-t))+\frac{3\sqrt{\mu}}{2}\coth(\frac{\sqrt{\mu}}{2}(T-t))$ and $c(t)=0$ into \eqref{eq:kernel_for_two}, we yield the following differential kernel corresponding to the unified NAG-G ODE:
\begin{align*}
	H(t,\tau) & =e^{-\int_{\tau}^{t}\left(\frac{\sqrt{\mu}}{2}\tanh\left(\frac{\sqrt{\mu}}{2}(T-s)\right)+\frac{3\sqrt{\mu}}{2}\coth\left(\frac{\sqrt{\mu}}{2}(T-s)\right)\right)\,ds}\\
	& =e^{\left[3\log\left(\sinh\left(\frac{\sqrt{\mu}}{2}(T-s)\right)\right)+\log\left(\cosh\left(\frac{\sqrt{\mu}}{2}(T-s)\right)\right)\right]_{\tau}^{t}}\\
	& =\frac{\sinh^{3}\left(\frac{\sqrt{\mu}}{2}(T-t)\right)\cosh\left(\frac{\sqrt{\mu}}{2}(T-t)\right)}{\sinh^{3}\left(\frac{\sqrt{\mu}}{2}(T-\tau)\right)\cosh\left(\frac{\sqrt{\mu}}{2}(T-\tau)\right)}.
\end{align*}

\section{Unified NAG Family}

\subsection{Proof of \cref{thm:mainthm_discrete}}
\label{app:thm:mainthm_discrete}
{
Note that when $\mu>0$, the inequality \eqref{eq:tk_condition1} can be written as
\begin{align*}
	0 & \geq\left(1-\sqrt{\mu s}\coth\left(\frac{\sqrt{\mu}}{2}\bft_{k+1}\right)\right)\frac{1}{\mu}\sinh^{2}\left(\frac{\sqrt{\mu}}{2}\bft_{k+1}\right)-\frac{1}{\mu}\sinh^{2}\left(\frac{\sqrt{\mu}}{2}\bft_{k}\right)\\
	& =\frac{1}{\mu}\sinh^{2}\left(\frac{\sqrt{\mu}}{2}\bft_{k+1}\right)-\sqrt{\frac{s}{\mu}}\sinh\left(\frac{\sqrt{\mu}}{2}\bft_{k+1}\right)\cosh\left(\frac{\sqrt{\mu}}{2}\bft_{k+1}\right)-\frac{1}{\mu}\sinh^{2}\left(\frac{\sqrt{\mu}}{2}\bft_{k}\right)\\
	& =\frac{1}{\mu}\cosh^{2}\left(\frac{\sqrt{\mu}}{2}\bft_{k+1}\right)-\sqrt{\frac{s}{\mu}}\sinh\left(\frac{\sqrt{\mu}}{2}\bft_{k+1}\right)\cosh\left(\frac{\sqrt{\mu}}{2}\bft_{k+1}\right)-\frac{1}{\mu}\cosh^{2}\left(\frac{\sqrt{\mu}}{2}\bft_{k}\right)\\
	& =\left(1-\sqrt{\mu s}\tanh\left(\frac{\sqrt{\mu}}{2}\bft_{k+1}\right)\right)\frac{1}{\mu}\cosh^{2}\left(\frac{\sqrt{\mu}}{2}\bft_{k+1}\right)-\frac{1}{\mu}\cosh^{2}\left(\frac{\sqrt{\mu}}{2}\bft_{k}\right).
\end{align*}
Thus, the following inequality holds for all $\mu\geq0$ (it clearly holds for $\mu=0$):
\begin{equation}
	\label{eq:cosh_ineq}
	\left(1-\frac{\mu\sqrt{s}\bft_{k+1}}{2}\tanhc\left(\frac{\sqrt{\mu}}{2}\bft_{k+1}\right)\right)\cosh^{2}\left(\frac{\sqrt{\mu}}{2}\bft_{k+1}\right)\leq\cosh^{2}\left(\frac{\sqrt{\mu}}{2}\bft_{k}\right).
\end{equation}
Using \eqref{eq:tk_condition1} and \eqref{eq:cosh_ineq}, we have
\begin{align*}
	& \mathcal{E}_{k+1} - \mathcal{E}_k\\
	& =\frac{1}{2}\cosh^{2}\left(\frac{\sqrt{\mu}}{2}\bft_{k+1}\right)\left\Vert z_{k+1}-x^{*}\right\Vert ^{2}-\frac{1}{2}\cosh^{2}\left(\frac{\sqrt{\mu}}{2}\bft_{k}\right)\left\Vert z_{k}-x^{*}\right\Vert ^{2}\\
	& \quad+\frac{\bft_{k+1}^{2}}{4}\sinhc^{2}\left(\frac{\sqrt{\mu}}{2}\bft_{k+1}\right)\left(f\left(x_{k+1}\right)-f\left(x^{*}\right)\right)-\frac{\bft_{k}^{2}}{4}\sinhc^{2}\left(\frac{\sqrt{\mu}}{2}\bft_{k}\right)\left(f\left(x_{k}\right)-f\left(x^{*}\right)\right)\\
	& \leq\frac{1}{2}\cosh^{2}\left(\frac{\sqrt{\mu}}{2}\bft_{k+1}\right)\left\Vert z_{k+1}-x^{*}\right\Vert ^{2}-\frac{1}{2}\left(1-\frac{\mu\sqrt{s}\bft_{k+1}}{2}\tanhc\left(\frac{\sqrt{\mu}}{2}\bft_{k+1}\right)\right)\cosh^{2}\left(\frac{\sqrt{\mu}}{2}\bft_{k+1}\right)\left\Vert z_{k}-x^{*}\right\Vert ^{2}\\
	& \quad+\frac{\bft_{k+1}^{2}}{4}\sinhc^{2}\left(\frac{\sqrt{\mu}}{2}\bft_{k+1}\right)\left(f\left(x_{k+1}\right)-f\left(x^{*}\right)\right)\\
	& \quad-\left(1-\frac{2\sqrt{s}}{\bft_{k+1}}\cothc\left(\frac{\sqrt{\mu}}{2}\bft_{k+1}\right)\right)\frac{\bft_{k+1}^{2}}{4}\sinhc^{2}\left(\frac{\sqrt{\mu}}{2}\bft_{k+1}\right)\left(f\left(x_{k}\right)-f\left(x^{*}\right)\right).
\end{align*}
Substituting
\[
z_{k+1}=y_{k}+\left(1-\frac{\mu\sqrt{s}\bft_{k+1}}{2}\tanhc\left(\frac{\sqrt{\mu}}{2}\bft_{k+1}\right)\right)\left(z_{k}-y_{k}\right)-\frac{\sqrt{s}\bft_{k+1}}{2}\tanhc\left(\frac{\sqrt{\mu}}{2}\bft_{k+1}\right)\nabla f\left(y_{k}\right)
\]
into the inequality above, we have
\begin{align*}
	& \mathcal{E}_{k+1}-\mathcal{E}_{k}\\
	& \leq\frac{1}{2}\cosh^{2}\left(\frac{\sqrt{\mu}}{2}\bft_{k+1}\right)\\
	& \quad\times\Bigg\Vert\left(1-\frac{\mu\sqrt{s}\bft_{k+1}}{2}\tanhc\left(\frac{\sqrt{\mu}}{2}\bft_{k+1}\right)\right)\left(z_{k}-y_{k}\right)\\
	& \qquad\quad-\frac{\sqrt{s}\bft_{k+1}}{2}\tanhc\left(\frac{\sqrt{\mu}}{2}\bft_{k+1}\right)\nabla f\left(y_{k}\right)-\left(x^{*}-y_{k}\right)\Bigg\Vert^{2}\\
	& \quad-\frac{1}{2}\left(1-\frac{\mu\sqrt{s}\bft_{k+1}}{2}\tanhc\left(\frac{\sqrt{\mu}}{2}\bft_{k+1}\right)\right)\cosh^{2}\left(\frac{\sqrt{\mu}}{2}\bft_{k+1}\right)\left\Vert \left(z_{k}-y_{k}\right)-\left(x^{*}-y_{k}\right)\right\Vert ^{2}\\
	& \quad+\frac{\bft_{k+1}^{2}}{4}\sinhc^{2}\left(\frac{\sqrt{\mu}}{2}\bft_{k+1}\right)\left(f\left(x_{k+1}\right)-f\left(x^{*}\right)\right)\\
	& \quad-\left(1-\frac{2\sqrt{s}}{\bft_{k+1}}\cothc\left(\frac{\sqrt{\mu}}{2}\bft_{k+1}\right)\right)\frac{\bft_{k+1}^{2}}{4}\sinhc^{2}\left(\frac{\sqrt{\mu}}{2}\bft_{k+1}\right)\left(f\left(x_{k}\right)-f\left(x^{*}\right)\right)\\
	& =\frac{1}{2}\cosh^{2}\left(\frac{\sqrt{\mu}}{2}\bft_{k+1}\right)\Bigg(\left(1-\frac{\mu\sqrt{s}\bft_{k+1}}{2}\tanhc\left(\frac{\sqrt{\mu}}{2}\bft_{k+1}\right)\right)^{2}\\
	& \qquad\quad-\left(1-\frac{\mu\sqrt{s}\bft_{k+1}}{2}\tanhc\left(\frac{\sqrt{\mu}}{2}\bft_{k+1}\right)\right)\Bigg)\left\Vert z_{k}-y_{k}\right\Vert ^{2}\\
	& \quad+\frac{\sqrt{s}\bft_{k+1}}{2}\sinhc\left(\frac{\sqrt{\mu}}{2}\bft_{k+1}\right)\cosh\left(\frac{\sqrt{\mu}}{2}\bft_{k+1}\right)\left\langle \nabla f\left(y_{k}\right),x^{*}-y_{k}\right\rangle \\
	& \quad+\frac{\mu\sqrt{s}\bft_{k+1}}{4}\sinhc\left(\frac{\sqrt{\mu}}{2}\bft_{k+1}\right)\cosh\left(\frac{\sqrt{\mu}}{2}\bft_{k+1}\right)\left\Vert x^{*}-y_{k}\right\Vert ^{2}\\
	& \quad-\frac{\sqrt{s}\bft_{k+1}}{2}\sinhc\left(\frac{\sqrt{\mu}}{2}\bft_{k+1}\right)\cosh\left(\frac{\sqrt{\mu}}{2}\bft_{k+1}\right)\\
	& \qquad\quad\times\left(1-\frac{\mu\sqrt{s}\bft_{k+1}}{2}\tanhc\left(\frac{\sqrt{\mu}}{2}\bft_{k+1}\right)\right)\left\langle \nabla f(y_{k}),z_{k}-y_{k}\right\rangle \\
	& \quad+\frac{s\bft_{k+1}^{2}}{8}\sinhc^{2}\left(\frac{\sqrt{\mu}}{2}\bft_{k+1}\right)\left\Vert \nabla f\left(y_{k}\right)\right\Vert ^{2}+\frac{\bft_{k+1}^{2}}{4}\sinhc^{2}\left(\frac{\sqrt{\mu}}{2}\bft_{k+1}\right)\left(f\left(x_{k+1}\right)-f\left(x^{*}\right)\right)\\
	& \quad-\left(1-\frac{2\sqrt{s}}{\bft_{k+1}}\cothc\left(\frac{\sqrt{\mu}}{2}\bft_{k+1}\right)\right)\frac{\bft_{k+1}^{2}}{4}\sinhc^{2}\left(\frac{\sqrt{\mu}}{2}\bft_{k+1}\right)\left(f\left(x_{k}\right)-f\left(x^{*}\right)\right).
\end{align*}
Since
\[
0\leq1-\sqrt{\mu s}\leq1-\frac{\mu\sqrt{s}\bft_{k+1}}{2}\tanhc\left(\frac{\sqrt{\mu}}{2}\bft_{k+1}\right)\leq1,
\]
we have
\begin{multline*}
	\frac{1}{2}\cosh^{2}\left(\frac{\sqrt{\mu}}{2}\bft_{k+1}\right)\Bigg(\left(1-\frac{\mu\sqrt{s}\bft_{k+1}}{2}\tanhc\left(\frac{\sqrt{\mu}}{2}\bft_{k+1}\right)\right)^{2}\\
	-\left(1-\frac{\mu\sqrt{s}\bft_{k+1}}{2}\tanhc\left(\frac{\sqrt{\mu}}{2}\bft_{k+1}\right)\right)\Bigg) \left\Vert z_{k}-y_{k}\right\Vert ^{2}\leq 0.
\end{multline*}
Therefore, we deduce that
\begin{align*}
	& \mathcal{E}_{k+1}-\mathcal{E}_k\\
	& \leq\frac{\sqrt{s}\bft_{k+1}}{2}\sinhc\left(\frac{\sqrt{\mu}}{2}\bft_{k+1}\right)\cosh\left(\frac{\sqrt{\mu}}{2}\bft_{k+1}\right)\left\langle \nabla f\left(y_{k}\right),x^{*}-y_{k}\right\rangle\\
	& \quad  +\frac{\mu\sqrt{s}\bft_{k+1}}{4}\sinhc\left(\frac{\sqrt{\mu}}{2}\bft_{k+1}\right)\cosh\left(\frac{\sqrt{\mu}}{2}\bft_{k+1}\right)\left\Vert x^{*}-y_{k}\right\Vert ^{2}\\
	& \quad-\frac{\sqrt{s}\bft_{k+1}}{2}\sinhc\left(\frac{\sqrt{\mu}}{2}\bft_{k+1}\right)\cosh\left(\frac{\sqrt{\mu}}{2}\bft_{k+1}\right) \\
	& \qquad\quad \times\left(1-\frac{\mu\sqrt{s}\bft_{k+1}}{2}\tanhc\left(\frac{\sqrt{\mu}}{2}\bft_{k+1}\right)\right)\left\langle \nabla f(y_{k}),z_{k}-y_{k}\right\rangle \\
	& \quad +\frac{s\bft_{k+1}^{2}}{8}\sinhc^{2}\left(\frac{\sqrt{\mu}}{2}\bft_{k+1}\right)\left\Vert \nabla f\left(y_{k}\right)\right\Vert ^{2}\\
	& \quad+\frac{\bft_{k+1}^{2}}{4}\sinhc^{2}\left(\frac{\sqrt{\mu}}{2}\bft_{k+1}\right)\left(f\left(x_{k+1}\right)-f\left(x^{*}\right)\right)\\
	& \quad -\left(1-\frac{2\sqrt{s}}{\bft_{k+1}}\cothc\left(\frac{\sqrt{\mu}}{2}\bft_{k+1}\right)\right)\frac{\bft_{k+1}^{2}}{4}\sinhc^{2}\left(\frac{\sqrt{\mu}}{2}\bft_{k+1}\right)\left(f\left(x_{k}\right)-f\left(x^{*}\right)\right).
\end{align*}
Now, it suffices to show that the right-hand side (RHS) of the inequality above is non-positive. By the $\mu$-strong convexity of $f$, we have
\[
0\geq f\left(y_{k}\right)-f\left(x^{*}\right)+\left\langle \nabla f\left(y_{k}\right),x^{*}-y_{k}\right\rangle +\frac{\mu}{2}\left\Vert x^{*}-y_{k}\right\Vert ^{2}.
\]
Moreover, it follows from the convexity and the $\frac{1}{s}$-smoothness of $f$ that
\[
0\geq f\left(y_{k}\right)-f\left(x_{k}\right)+\left\langle \nabla f(y_{k}),x_{k}-y_{k}\right\rangle
\]
and
\[
0\geq f(x_{k+1})-f(y_{k})+\frac{s}{2}\left\Vert \nabla f(y_{k})\right\Vert ^{2},
\]
respectively. 
Note that $$x_{k}-y_{k}=-\frac{\tau_{k}}{1-\tau_{k}}\left(z_{k}-y_{k}\right)=-\frac{\frac{2\sqrt{s}}{\bft_{k+1}}\cothc\left(\frac{\sqrt{\mu}}{2}\bft_{k+1}\right)-\mu s}{1-\frac{2\sqrt{s}}{\bft_{k+1}}\cothc\left(\frac{\sqrt{\mu}}{2}\bft_{k+1}\right)}\left(z_{k}-y_{k}\right).$$ Taking a weighted sum of the inequalities above yields (the assumption \eqref{eq:tk_condition2} ensures that these weights are non-negative for $k\geq1$, and the case $k=0$ is trivial because $y_0=x_0$)
\begin{align*}
	0 & \geq\frac{2\sqrt{s}}{\bft_{k+1}}\cothc\left(\frac{\sqrt{\mu}}{2}\bft_{k+1}\right)\frac{\bft_{k+1}^{2}}{4}\sinhc^{2}\left(\frac{\sqrt{\mu}}{2}\bft_{k+1}\right)\\
	& \qquad\quad\times\left[f\left(y_{k}\right)-f\left(x^{*}\right)+\left\langle \nabla f\left(y_{k}\right),x^{*}-y_{k}\right\rangle +\frac{\mu}{2}\left\Vert x^{*}-y_{k}\right\Vert ^{2}\right]\\
	& \quad+\left(1-\frac{2\sqrt{s}}{\bft_{k+1}}\cothc\left(\frac{\sqrt{\mu}}{2}\bft_{k+1}\right)\right)\frac{\bft_{k+1}^{2}}{4}\sinhc^{2}\left(\frac{\sqrt{\mu}}{2}\bft_{k+1}\right)\\
	& \qquad\quad\times\left[f\left(y_{k}\right)-f\left(x_{k}\right)+\left\langle \nabla f(y_{k}),x_{k}-y_{k}\right\rangle \right]\\
	& \quad+\frac{\bft_{k+1}^{2}}{4}\sinhc^{2}\left(\frac{\sqrt{\mu}}{2}\bft_{k+1}\right)\left[f(x_{k+1})-f(y_{k})+\frac{s}{2}\left\Vert \nabla f(y_{k})\right\Vert ^{2}\right]\\
	& =\frac{\sqrt{s}\bft_{k+1}}{2}\sinhc\left(\frac{\sqrt{\mu}}{2}\bft_{k+1}\right)\cosh\left(\frac{\sqrt{\mu}}{2}\bft_{k+1}\right)\left\langle \nabla f\left(y_{k}\right),x^{*}-y_{k}\right\rangle \\
	& \quad+\frac{\mu\sqrt{s}\bft_{k+1}}{4}\sinhc\left(\frac{\sqrt{\mu}}{2}\bft_{k+1}\right)\cosh\left(\frac{\sqrt{\mu}}{2}\bft_{k+1}\right)\left\Vert x^{*}-y_{k}\right\Vert ^{2}\\
	& \quad-\left(\frac{2\sqrt{s}}{\bft_{k+1}}\cothc\left(\frac{\sqrt{\mu}}{2}\bft_{k+1}\right)-\mu s\right)\frac{\bft_{k+1}^{2}}{4}\sinhc^{2}\left(\frac{\sqrt{\mu}}{2}\bft_{k+1}\right)\left\langle \nabla f(y_{k}),z_{k}-y_{k}\right\rangle \\
	& \quad+\frac{s\bft_{k+1}^{2}}{8}\sinhc^{2}\left(\frac{\sqrt{\mu}}{2}\bft_{k+1}\right)\left\Vert \nabla f\left(y_{k}\right)\right\Vert ^{2}\\
	& \quad+\frac{2\sqrt{s}}{\bft_{k+1}}\cothc\left(\frac{\sqrt{\mu}}{2}\bft_{k+1}\right)\frac{\bft_{k+1}^{2}}{4}\sinhc^{2}\left(\frac{\sqrt{\mu}}{2}\bft_{k+1}\right)\left(f\left(y_{k}\right)-f\left(x^{*}\right)\right)\\
	& \quad+\left(1-\frac{2\sqrt{s}}{\bft_{k+1}}\cothc\left(\frac{\sqrt{\mu}}{2}\bft_{k+1}\right)\right)\frac{\bft_{k+1}^{2}}{4}\sinhc^{2}\left(\frac{\sqrt{\mu}}{2}\bft_{k+1}\right)\left(f\left(y_{k}\right)-f\left(x_{k}\right)\right)\\
	& \quad+\frac{\bft_{k+1}^{2}}{4}\sinhc^{2}\left(\frac{\sqrt{\mu}}{2}\bft_{k+1}\right)\left(f(x_{k+1})-f(y_{k})\right)\\
	& =\frac{\sqrt{s}\bft_{k+1}}{2}\sinhc\left(\frac{\sqrt{\mu}}{2}\bft_{k+1}\right)\cosh\left(\frac{\sqrt{\mu}}{2}\bft_{k+1}\right)\left\langle \nabla f\left(y_{k}\right),x^{*}-y_{k}\right\rangle \\
	& \quad+\frac{\mu\sqrt{s}\bft_{k+1}}{4}\sinhc\left(\frac{\sqrt{\mu}}{2}\bft_{k+1}\right)\cosh\left(\frac{\sqrt{\mu}}{2}\bft_{k+1}\right)\left\Vert x^{*}-y_{k}\right\Vert ^{2}\\
	& \quad-\frac{\sqrt{s}\bft_{k+1}}{2}\sinhc\left(\frac{\sqrt{\mu}}{2}\bft_{k+1}\right)\cosh\left(\frac{\sqrt{\mu}}{2}\bft_{k+1}\right)\\
	& \qquad\quad\times\left(1-\frac{\mu\sqrt{s}\bft_{k+1}}{2}\tanhc\left(\frac{\sqrt{\mu}}{2}\bft_{k+1}\right)\right)\left\langle \nabla f(y_{k}),z_{k}-y_{k}\right\rangle \\
	& \quad+\frac{s\bft_{k+1}^{2}}{8}\sinhc^{2}\left(\frac{\sqrt{\mu}}{2}\bft_{k+1}\right)\left\Vert \nabla f\left(y_{k}\right)\right\Vert ^{2}\\
	& \quad+\frac{\bft_{k+1}^{2}}{4}\sinhc^{2}\left(\frac{\sqrt{\mu}}{2}\bft_{k+1}\right)\left(f\left(x_{k+1}\right)-f\left(x^{*}\right)\right)\\
	& \quad-\left(1-\frac{2\sqrt{s}}{\bft_{k+1}}\cothc\left(\frac{\sqrt{\mu}}{2}\bft_{k+1}\right)\right)\frac{\bft_{k+1}^{2}}{4}\sinhc^{2}\left(\frac{\sqrt{\mu}}{2}\bft_{k+1}\right)\left(f\left(x_{k}\right)-f\left(x^{*}\right)\right).
\end{align*}
This completes the proof.}

\subsection{Constant timestep scheme}
\label{app:constant}

In this section, we show that the sequence $(\bft_k)$ defined in \eqref{eq:constant_tk} satisfies the conditions \eqref{eq:tk_condition2} and \eqref{eq:tk_condition1}. For convenience, we assume $\mu>0$ (the case $\mu=0$ can be handled easily). The condition \eqref{eq:tk_condition2} follows from
\begin{align*}
	\frac{2\sqrt{s}}{\bft_{k}}\cothc\left(\frac{\sqrt{\mu}}{2}\bft_{k}\right) & =\sqrt{\mu s}\coth\left(\frac{\sqrt{\mu}}{2}\bft_{k}\right)\\
	& \leq\sqrt{\mu s}\coth\left(\frac{\sqrt{\mu}}{2}\bft_{2}\right)\\
	& =\sqrt{\mu s}\coth\left(-\log\left(1-\sqrt{\mu s}\right)\right)\\
	& =\sqrt{\mu s}\frac{1+e^{2\log\left(1-\sqrt{\mu s}\right)}}{1-e^{2\log\left(1-\sqrt{\mu s}\right)}}\\
	& =\sqrt{\mu s}\frac{1+\left(1-\sqrt{\mu s}\right)^{2}}{1-\left(1-\sqrt{\mu s}\right)^{2}}\\
	& \leq1,
\end{align*}
where the last inequality holds because $\sqrt{\mu s}\in (0,1)$. To prove \eqref{eq:tk_condition1}, it suffices to show that the inequality
\[
\sinh^{2}\left(\frac{\sqrt{\mu}}{2}t\right)-\sqrt{\mu s}\sinh\left(\frac{\sqrt{\mu}}{2}t\right)\cosh\left(\frac{\sqrt{\mu}}{2}t\right) - \sinh^{2}\left(\frac{\sqrt{\mu}}{2}t+\frac{1}{2}\log\left(1-\sqrt{\mu s}\right)\right)\leq0
\]
holds for all $t\in\mathbb{R}$. Letting $r=e^{\frac{\sqrt{\mu}}{2}t}$, this inequality can be expressed as
\[
\frac{r^{2}+r^{-2}-2}{4}-\sqrt{\mu s}\frac{r^{2}-r^{-2}}{4}- \frac{\left(1-\sqrt{\mu s}\right)r^{2}+\left(1-\sqrt{\mu s}\right)^{-1}r^{-2}-2}{4}\leq0.
\]
Letting $q=r^2$ and multiplying both sides by $4q$,  the inequality can be rewritten as
\begin{align*}
	0 & \geq q^{2}+1-2q-\sqrt{\mu s}\left(q^{2}-1\right)-\left(1-\sqrt{\mu s}\right)q^{2}-\left(1-\sqrt{\mu s}\right)^{-1}+2q\\
	& =1+\sqrt{\mu s}-\frac{1}{1-\sqrt{\mu s}}\\
	& =\frac{-\mu s}{1-\sqrt{\mu s}},
\end{align*}
which clearly holds.

\subsection{Adaptive timestep scheme}
\label{app:variable}

In this section, we show that for the sequence $\left(\bft_{k}\right)$ defined by \eqref{eq:adaptive_timestep},
\begin{itemize}
	\item the sequence $\left(\bft_{k}\right)$ is well-defined, and
	\item the conditions \eqref{eq:tk_condition0} and \eqref{eq:tk_condition} hold when $\lim_{s\to0}\bft_0=0$.
\end{itemize}

\paragraph{The sequence $\left(\bft_{k}\right)$ is well-defined.}

Because $$\frac{4}{\bft_{k+1}^{2}}\cschc^{2}\left(\frac{\sqrt{\mu}}{2}\bft_{k+1}\right)+\mu=\frac{4}{\bft_{k+1}^{2}}\cothc^{2}\left(\frac{\sqrt{\mu}}{2}\bft_{k+1}\right),$$ the updating rule \eqref{eq:adaptive_timestep} is equivalent to
\begin{multline}
		\label{eq:update_in_cothc}
		\frac{4}{\bft_{k+1}^{2}}\cothc^{2}\left(\frac{\sqrt{\mu}}{2}\bft_{k+1}\right)=\left(1-\frac{2\sqrt{s}}{\bft_{k+1}}\cothc\left(\frac{\sqrt{\mu}}{2}\bft_{k+1}\right)\right)\frac{4}{\bft_{k}^{2}}\cothc^{2}\left(\frac{\sqrt{\mu}}{2}\bft_{k}\right)\\
		+\frac{2\mu\sqrt{s}}{\bft_{k+1}}\cothc\left(\frac{\sqrt{\mu}}{2}\bft_{k+1}\right),\;\bft_{k+1}>0.
\end{multline}
Introduce a sequence $\left(\alpha_{k}\right)_{k=-1}^{\infty}$ such that  $\alpha_{k}=\frac{2\sqrt{s}}{\bft_{k+1}}\cothc\left(\frac{\sqrt{\mu}}{2}\bft_{k+1}\right)$. As $t\mapsto\frac{2\sqrt{s}}{t}\cothc\left(\frac{\sqrt{\mu}}{2}t\right)$ is a bijective map from $(0,\infty)$ to $\left(\sqrt{\mu s},\infty\right)$, the sequences $\left(\bft_{k}\right)$ and $\left(\alpha_{k}\right)$ have a one-to-one relationship. Thus, the updating rule \eqref{eq:update_in_cothc} is equivalent to
\begin{equation}
	\label{eq:update_in_alpha}
	\alpha_{k}^{2}=\left(1-\alpha_{k}\right)\alpha_{k-1}^{2}+\mu s\alpha_{k},\quad \alpha_{k}>\sqrt{\mu s},
\end{equation}
which admits a unique solution in $(\sqrt{\mu s},\infty)$ when $\alpha_{k-1}>\sqrt{\mu s}$. Thus, the sequence $\left(\bft_{k}\right)$ is well-defined.

\paragraph{The sequence $\left(t_{k}\right)$ satisfies the conditions \eqref{eq:tk_condition0} and \eqref{eq:tk_condition}.}

Define a function $A(t)$ as
\begin{equation}
	\label{eq:At}
	A(t):=\frac{ t^{2}}{4}\sinhc^{2}\left(\frac{\sqrt{\mu}}{2}t\right).
\end{equation}
For $t\in\left(0,\infty\right)$, it follows from \eqref{eq:adaptive_timestep} that 
\[
\dot{A}\left(\bft_{\bfk(t)+1}\right)=\frac{A\left(\bft_{\bfk(t)+1}\right)-A\left(t\right)}{\sqrt{s}}=\frac{A\left(\bft_{\bfk(t)+1}\right)-A\left(t\right)}{\bft_{\bfk(t)+1}-t}\frac{\bft_{\bfk(t)+1}-t}{\sqrt{s}}.
\]
Because $\bft_{\bfk(t)+1}\rightarrow t$ as $s\rightarrow0$, taking the limit $s\rightarrow0$ in the equation above yields
\[
1=\lim_{s\rightarrow0}\frac{\bft_{\bfk(t)+1}-t}{\sqrt{s}}.
\]
Thus, the condition \eqref{eq:tk_condition1} holds. 

\subsection{Equivalence between the adaptive timestep scheme and the original NAG}
\label{app:prop:equivalence}

In this section, we show that the adaptive timestep scheme (\cref{sec:variable}) with $\bft_0>0$ is equivalent to the original NAG \eqref{eq:original_nag} with $\gamma_{0}=\frac{4}{\bft_{0}^{2}}\cothc^{2}\left(\frac{\sqrt{\mu}}{2}\bft_{0}\right)>\mu$.

We first show that the sequences $\left(\alpha_{k}\right)_{k=0}^{\infty}$ and $\left(\gamma_{k}\right)_{k=0}^{\infty}$ generated in the original NAG \eqref{eq:original_nag} with $\gamma_{0}=\frac{4}{\bft_{0}^{2}}\cothc^{2}\left(\frac{\sqrt{\mu}}{2}\bft_{0}\right)>\mu$ can be written as $\alpha_{k}=\frac{2\sqrt{s}}{\bft_{k+1}}\cothc\left(\frac{\sqrt{\mu}}{2}\bft_{k+1}\right)$ and $\gamma_{k}=\frac{4}{\bft_{k}^{2}}\cothc^{2}\left(\frac{\sqrt{\mu}}{2}\bft_{k}\right)$, where the sequence $\left(\bft_{k}\right)_{k=0}^{\infty}$ is defined as \eqref{eq:adaptive_timestep}. Note that the equality \eqref{eq:alpha_k_update} implies
\[
\gamma_{k+1}=\left(1-\alpha_{k}\right)\gamma_{k}+\mu\alpha_{k}=\frac{\alpha_{k}^{2}}{s}.
\]
Thus, the updating rule for $\alpha_k$ \eqref{eq:alpha_k_update} can be written as
\[
\frac{1}{s}\alpha_{k}^{2}=\left(1-\alpha_{k}\right)\frac{\alpha_{k-1}^{2}}{s}+\mu\alpha_{k},
\]
where we define $\alpha_{-1}:=\sqrt{s\gamma_{0}}=\frac{2\sqrt{s}}{\bft_{0}}\cothc\left(\frac{\sqrt{\mu}}{2}\bft_{0}\right)>\sqrt{\mu s}$. This implies that  the sequence $\left(\alpha_{k}\right)_{k=-1}^{\infty}$ in the original NAG and the sequence $\left(\alpha_{k}\right)_{k=-1}^{\infty}$ defined in Section~\ref{app:variable} are identical. Thus, we have $\alpha_{k}=\frac{2\sqrt{s}}{\bft_{k+1}}\cothc\left(\frac{\sqrt{\mu}}{2}\bft_{k+1}\right)$ and $\gamma_{k}=\frac{\alpha_{k-1}^2}{s}=\frac{4}{\bft_{k}^{2}}\cothc^{2}\left(\frac{\sqrt{\mu}}{2}\bft_{k}\right)$.

Now, we  show that the parameters $\tau_k$ and $\delta_k$ for the original NAG are equal to those for our adaptive timestep scheme.
In the original NAG, we have
\begin{align*}
	\left(\alpha_{k}-\mu s\right)\left(\gamma_{k}+\mu\alpha_{k}\right) & =\alpha_{k}\gamma_{k}+\mu\alpha_{k}^{2}-\mu s\gamma_{k}-\mu^{2}s\alpha_{k}\\
	& =\mu s\gamma_{k+1}+\alpha_{k}\gamma_{k}-\mu s\gamma_{k}-\mu^{2}s\alpha_{k}\\
	& =\mu s\left(\left(1-\alpha_{k}\right)\gamma_{k}+\mu\alpha_{k}\right)+\alpha_{k}\gamma_{k}-\mu s\gamma_{k}-\mu^{2}s\alpha_{k}\\
	& =(1-\mu s)\alpha_{k}\gamma_{k}.
\end{align*}
Therefore, we have
$$\tau_{k}=\frac{\alpha_{k}\gamma_{k}}{\gamma_{k}+\mu\alpha_{k}}=\frac{\alpha_{k}-\mu s}{1-\mu s}=\frac{\frac{2\sqrt{s}}{\bft_{k+1}}\cothc\left(\frac{\sqrt{\mu}}{2}\bft_{k+1}\right)-\mu s}{1-\mu s}$$
and
$$\delta_{k}=\frac{\alpha_{k}}{\gamma_{k+1}}=\frac{s}{\alpha_{k}}=\frac{\sqrt{s}\bft_{k+1}}{2}\tanhc\left(\frac{\sqrt{\mu}}{2}\bft_{k+1}\right).$$
Thus, the ogirinal Nesterov's method with $\gamma_{0}=\frac{4}{\bft_{0}^{2}}\cothc^{2}\left(\frac{\sqrt{\mu}}{2}\bft_{0}\right)>\mu$ is equivalent to the adaptive timestep scheme.

\section{Higher-Order Extension}
\label{app:higher-order}

\subsection{Limiting ODE}
\label{app:limiting_ode_higher}

\paragraph{Limiting ODE of the unified accelerated tensor method family.}
We show that if the sequence $(\bft_k)$ satisfies the conditions \eqref{eq:tk_condition0_higher} and \eqref{eq:tk_condition_higher}, then the unified accelerated tensor method family \eqref{eq:hyperbolic_scheme} converges to the unified accelerated tensor flow \eqref{eq:hyperbolic_flow} under the identifications $x_k=X(\bft_k)$ and $z_k=Z(\bft_k)$.

For convenience, we assume that $\mu>0$ (the case $\mu=0$ can be handled easily). Define a function $A:[0,\infty)\rightarrow\mathbb{R}$ as
\begin{equation}
	\label{eq:At_higher}
	A(t)=Ct^{p}\sinhc_{p}^{p}\left(C^{1/p}\mu^{1/p}t\right)=\frac{1}{\mu}\sinh_{p}^{p}\left(C^{1/p}\mu^{1/p}t\right)
\end{equation}
so that $A_k=A(\bft_k)$. It follows from the step \eqref{eq:hyperbolic_scheme2} that
\begin{align*}
	\dot{X}(t) & =\lim_{s\rightarrow0}\frac{x_{\bfk(t)+1}-x_{\bfk(t)}}{\bft_{\bfk(t)+1}-t}\\
	& =\lim_{s\rightarrow0}\frac{x_{\bfk(t)+1}-x_{\bfk(t)}}{s^{1/p}}\\
	& =\lim_{s\rightarrow0}\frac{y_{\bfk(t)}-x_{\bfk(t)}}{s^{1/p}}\\
	& =\lim_{s\rightarrow0}\frac{A_{\bfk(t)+1}-A_{\bfk(t)}}{s^{1/p}A_{\bfk(t)+1}}\left(z_{\bfk(t)}-x_{\bfk(t)}\right)\\
	& =\lim_{s\rightarrow0}\frac{A\left(\bft_{\bfk(t)+1}\right)-A(t)}{s^{1/p}A\left(\bft_{\bfk(t)+1}\right)}\left(Z(t)-X(t)\right)\\
	& =\frac{\dot{A}(t)}{A(t)}\left(Z(t)-X(t)\right)\\
	& =pC^{1/p}\mu^{1/p}\coth_{p}\left(C^{1/p}\mu^{1/p}t\right)\left(Z(t)-X(t)\right),
\end{align*}
where we used $\Vert x_{k+1}-y_{k}\Vert=o(s^{1/p})$ \citep[see][Lemma~2.2]{wibisono2016} for the third equality. Using the step \eqref{eq:our_zk}, we have
\begin{align*}
	\frac{d}{dt}\nabla h(Z(t)) & =\lim_{s\rightarrow0}\frac{\nabla h\left(z_{\bfk(t)+1}\right)-\nabla h\left(z_{\bfk(t)}\right)}{\bft_{\bfk(t)+1}-t}\\
	& =\lim_{s\rightarrow0}\frac{\nabla h\left(z_{\bfk(t)+1}\right)-\nabla h\left(z_{\bfk(t)}\right)}{s^{1/p}}\\
	& =\lim_{s\rightarrow0}\frac{A_{\bfk(t)+1}-A_{\bfk(t)}}{s^{1/p}\left(1+\mu A_{\bfk(t)}\right)}\left(\mu\nabla h\left(x_{\bfk(t)+1}\right)-\mu\nabla h\left(z_{\bfk(t)+1}\right)-\nabla f\left(x_{\bfk(t)+1}\right)\right)\\
	& =\lim_{s\rightarrow0}\frac{A\left(\bft_{\bfk(t)+1}\right)-A(t)}{s^{1/p}\left(1+\mu A(t)\right)}\left(\mu\nabla h\left(X(t)\right)-\mu\nabla h\left(X(t)\right)-\nabla f\left(X(t)\right)\right)\\
	& =\frac{\dot{A}(t)}{1+\mu A(t)}\left(\mu\nabla h\left(X(t)\right)-\mu\nabla h\left(X(t)\right)-\nabla f\left(X(t)\right)\right)\\
	& =\frac{C^{1/p}p}{\mu^{(p-1)/p}}\tanh_{p}^{p-1}\left(C^{1/p}\mu^{1/p}t\right)\left(\mu\nabla h\left(X(t)\right)-\mu\nabla h\left(X(t)\right)-\nabla f\left(X(t)\right)\right).
\end{align*}
Thus, we obtain the system of ODEs \eqref{eq:hyperbolic_flow}.

\paragraph{Limiting ODE of the unified accelerated tensor method.}
We check that the sequence $(\bft_k)$ defined in \eqref{eq:our_specific_Ak} satisfies the condition \eqref{eq:tk_condition_higher}. It is easy to check that the function $A(t)$ defined in \eqref{eq:At_higher} satisfies
\[
\dot{A}(t)=C^{1/p}p\mu^{\frac{1-p}{p}}\sinh_{p}^{p-1}\left(C^{1/p}\mu^{1/p}t\right)\cosh_{p}\left(C^{1/p}\mu^{1/p}t\right)=C^{1/p}pA(t)^{\frac{p-1}{p}}(1+\mu A(t))^{\frac{1}{p}}
\]
and that the sequence $(\bft_k)$ defined in \eqref{eq:our_specific_Ak} satisfies
\[
\frac{A\left(\bft_{k+1}\right)-A\left(\bft_{k}\right)}{s^{1/p}}-C^{1/p}pA\left(\bft_{k+1}\right)^{\frac{p-1}{p}}\left(1+\mu A\left(\bft_{k}\right)\right)^{\frac{1}{p}}=0.
\]
Now, substituting $k=\bfk(t)$ into the above equality and taking the limit $s\rightarrow0$, we have $\lim_{s\rightarrow0}\frac{\bft_{\bfk(t)+1}-t}{s^{1/p}}=1$.

\subsection{Proof of Theorem~\ref{thm:mainthm_discrete_higher}}
\label{app:thm:mainthm_discrete_higher}

By the Bregman three-point identity \eqref{eq:three-point} with $x=x^{*}$, $y=z_{k+1}$, $z=x_{k+1}$ and the non-negativity of Bregman divergence, we have
\begin{align*}
	D_{h}\left(x^{*},z_{k+1}\right) & =D_{h}\left(x^{*},x_{k+1}\right)-\left\langle \nabla h\left(z_{k+1}\right)-\nabla h\left(x_{k+1}\right),x^{*}-z_{k+1}\right\rangle -D_{h}\left(z_{k+1},x_{k+1}\right)\\
	& \leq D_{h}\left(x^{*},x_{k+1}\right)-\left\langle \nabla h\left(z_{k+1}\right)-\nabla h\left(x_{k+1}\right),x^{*}-z_{k+1}\right\rangle .
\end{align*}
Thus, we can bound the difference of the discrete-time energy function \eqref{eq:discrete_energy_higher} as follows:
\begin{align*}
	& \mathcal{E}_{k+1}-\mathcal{E}_{k} \\
	& =\left(1+\mu A_{k+1}\right)D_{h}\left(x^{*},z_{k+1}\right)-\left(1+\mu A_{k}\right)D_{h}\left(x^{*},z_{k}\right)\\
	& \quad+A_{k+1}\left(f\left(x_{k+1}\right)-f\left(x^{*}\right)\right)-A_{k}\left(f\left(x_{k}\right)-f\left(x^{*}\right)\right)\\
	& =\mu\left(A_{k+1}-A_{k}\right)D_{h}\left(x^{*},z_{k+1}\right)\\
	& \quad+\left(A_{k+1}-A_{k}\right)\left(f\left(x_{k+1}\right)-f\left(x^{*}\right)\right)+A_{k}\left(f\left(x_{k+1}\right)-f\left(x_{k}\right)\right)\\
	& \quad+\left(1+\mu A_{k}\right)\left(-h\left(z_{k+1}\right)-\left\langle \nabla h\left(z_{k+1}\right),x^{*}-z_{k+1}\right\rangle +h\left(z_{k}\right)+\left\langle \nabla h\left(z_{k}\right),x^{*}-z_{k}\right\rangle \right)\\
	& \leq\mu\left(A_{k+1}-A_{k}\right)D_{h}\left(x^{*},x_{k+1}\right)-\mu\left(A_{k+1}-A_{k}\right)\left\langle \nabla h\left(z_{k+1}\right)-\nabla h\left(x_{k+1}\right),x^{*}-z_{k+1}\right\rangle \\
	& \quad+\left(A_{k+1}-A_{k}\right)\left(f\left(x_{k+1}\right)-f\left(x^{*}\right)\right)+A_{k}\left(f\left(x_{k+1}\right)-f\left(x_{k}\right)\right)\\
	& \quad+\left(1+\mu A_{k}\right)\left(-h\left(z_{k+1}\right)-\left\langle \nabla h\left(z_{k+1}\right),x^{*}-z_{k+1}\right\rangle +h\left(z_{k}\right)+\left\langle \nabla h\left(z_{k}\right),x^{*}-z_{k}\right\rangle \right).
\end{align*}

By the ($\mu$-uniform) convexity of $f$ with respect to $h$, the $p$-th order $1$-uniform convexity of $h$, and the property \eqref{eq:M_ineq} of the higher-order gradient update operator $G_{p,M}$, the following inequalities hold:
\begin{align*}
	0 & \geq f\left(x_{k+1}\right)-f\left(x^{*}\right)+\left\langle \nabla f\left(x_{k+1}\right),x^{*}-x_{k+1}\right\rangle +\mu D_{h}\left(x^{*},x_{k+1}\right)\\
	0 & \geq f\left(x_{k+1}\right)-f\left(x_{k}\right)+\left\langle \nabla f\left(x_{k+1}\right),x_{k}-x_{k+1}\right\rangle \\
	0 & \geq Ms^{\frac{1}{p-1}}\left\Vert \nabla f\left(x_{k+1}\right)\right\Vert ^{\frac{p}{p-1}}-\left\langle \nabla f\left(x_{k+1}\right),y_{k}-x_{k+1}\right\rangle \\
	0 & \geq h\left(z_{k}\right)-h\left(z_{k+1}\right)+\left\langle \nabla h\left(z_{k}\right),z_{k+1}-z_{k}\right\rangle +\frac{1}{p}\left\Vert z_{k+1}-z_{k}\right\Vert ^{p}.
\end{align*}
Taking a weighted sum of these inequalities yields
\begin{align*}
	0 & \geq\left(A_{k+1}-A_{k}\right)\left[f\left(x_{k+1}\right)-f\left(x^{*}\right)+\left\langle \nabla f\left(x_{k+1}\right),x^{*}-x_{k+1}\right\rangle +\mu D_{h}\left(x^{*},x_{k+1}\right)\right]\\
	& \quad+A_{k}\left[f\left(x_{k+1}\right)-f\left(x_{k}\right)+\left\langle \nabla f\left(x_{k+1}\right),x_{k}-x_{k+1}\right\rangle \right]\\
	& \quad+A_{k+1}\left[Ms^{\frac{1}{p-1}}\left\Vert \nabla f\left(x_{k+1}\right)\right\Vert ^{\frac{p}{p-1}}-\left\langle \nabla f\left(x_{k+1}\right),y_{k}-x_{k+1}\right\rangle \right]\\
	& \quad+\left(1+\mu A_{k}\right)\left[h\left(z_{k}\right)-h\left(z_{k+1}\right)+\left\langle \nabla h\left(z_{k}\right),z_{k+1}-z_{k}\right\rangle +\frac{1}{p}\left\Vert z_{k+1}-z_{k}\right\Vert ^{p}\right]\\
	& \geq\mathcal{E}_{k+1}-\mathcal{E}_{k}\\
	& \quad-\mu\left(A_{k+1}-A_{k}\right)D_{h}\left(x^{*},x_{k+1}\right)+\mu\left(A_{k+1}-A_{k}\right)\left\langle \nabla h\left(z_{k+1}\right)-\nabla h\left(x_{k+1}\right),x^{*}-z_{k+1}\right\rangle \\
	& \quad-\left(A_{k+1}-A_{k}\right)\left(f\left(x_{k+1}\right)-f\left(x^{*}\right)\right)-A_{k}\left(f\left(x_{k+1}\right)-f\left(x_{k}\right)\right)\\
	& \quad-\left(1+\mu A_{k+1}\right)\left(-h\left(z_{k+1}\right)-\left\langle \nabla h\left(z_{k+1}\right),x^{*}-z_{k+1}\right\rangle +h\left(z_{k}\right)+\left\langle \nabla h\left(z_{k}\right),x^{*}-z_{k}\right\rangle \right)\\
	& \quad+\left(A_{k+1}-A_{k}\right)\left[f\left(x_{k+1}\right)-f\left(x^{*}\right)+\left\langle \nabla f\left(x_{k+1}\right),x^{*}-x_{k+1}\right\rangle +\mu D_{h}\left(x^{*},x_{k+1}\right)\right]\\
	& \quad+A_{k}\left[f\left(x_{k+1}\right)-f\left(x_{k}\right)+\left\langle \nabla f\left(x_{k+1}\right),x_{k}-x_{k+1}\right\rangle \right]\\
	& \quad+A_{k+1}\left[Ms^{\frac{1}{p-1}}\left\Vert \nabla f\left(x_{k+1}\right)\right\Vert ^{\frac{p}{p-1}}-\left\langle \nabla f\left(x_{k+1}\right),y_{k}-x_{k+1}\right\rangle \right]\\
	& \quad+\left(1+\mu A_{k}\right)\left[h\left(z_{k}\right)-h\left(z_{k+1}\right)+\left\langle \nabla h\left(z_{k}\right),z_{k+1}-z_{k}\right\rangle +\frac{1}{p}\left\Vert z_{k+1}-z_{k}\right\Vert ^{p}\right]\\
	& =\mathcal{E}_{k+1}-\mathcal{E}_{k}\\
	& \quad+\left\langle \nabla f\left(x_{k+1}\right),\left(A_{k+1}-A_{k}\right)\left(x^{*}-x_{k+1}\right)+A_{k}\left(x_{k}-x_{k+1}\right)+A_{k+1}\left(x_{k+1}-y_{k}\right)\right\rangle \\
	& \quad+\left(1+\mu A_{k}\right)\left\langle \nabla h\left(z_{k+1}\right)-\nabla h\left(z_{k}\right),x^{*}-z_{k+1}\right\rangle +\frac{1+\mu A_{k}}{p}\left\Vert z_{k+1}-z_{k}\right\Vert ^{p}\\
	& \quad+\mu\left(A_{k+1}-A_{k}\right)\left\langle \nabla h\left(z_{k+1}\right)-\nabla h\left(x_{k+1}\right),x^{*}-z_{k+1}\right\rangle +MA_{k+1}s^{\frac{1}{p-1}}\left\Vert \nabla f\left(x_{k+1}\right)\right\Vert ^{\frac{p}{p-1}}.
\end{align*}
Substituting \eqref{eq:our_zk} with the term $\left(1+\mu A_{k}\right)\left\langle \nabla h\left(z_{k+1}\right)-\nabla h\left(z_{k}\right),x^{*}-z_{k+1}\right\rangle $, we have
\begin{align*}
	0& \geq\mathcal{E}_{k+1}-\mathcal{E}_{k}\\
	& \quad+\left\langle \nabla f\left(x_{k+1}\right),\left(A_{k+1}-A_{k}\right)\left(x^{*}-x_{k+1}\right)+A_{k}\left(x_{k}-x_{k+1}\right)+A_{k+1}\left(x_{k+1}-y_{k}\right)\right\rangle \\
	& \quad+\left(A_{k+1}-A_{k}\right)\left\langle \mu\nabla h\left(x_{k+1}\right)-\mu\nabla h\left(z_{k+1}\right)-\nabla f\left(x_{k+1}\right),x^{*}-z_{k+1}\right\rangle \\
	& \quad +\frac{1+\mu A_{k}}{p}\left\Vert z_{k+1}-z_{k}\right\Vert ^{p}\\
	& \quad+\mu\left(A_{k+1}-A_{k}\right)\left\langle \nabla h\left(z_{k+1}\right)-\nabla h\left(x_{k+1}\right),x^{*}-z_{k+1}\right\rangle +MA_{k+1}s^{\frac{1}{p-1}}\left\Vert \nabla f\left(x_{k+1}\right)\right\Vert ^{\frac{p}{p-1}}\\
	& =\mathcal{E}_{k+1}-\mathcal{E}_{k}\\
	& \quad+\left\langle \nabla f\left(x_{k+1}\right),\left(A_{k+1}-A_{k}\right)\left(z_{k+1}-x_{k+1}\right)+A_{k}\left(x_{k}-x_{k+1}\right)+A_{k+1}\left(x_{k+1}-y_{k}\right)\right\rangle \\
	& \quad+\frac{1+\mu A_{k}}{p}\left\Vert z_{k+1}-z_{k}\right\Vert ^{p}+MA_{k+1}s^{\frac{1}{p-1}}\left\Vert \nabla f\left(x_{k+1}\right)\right\Vert ^{\frac{p}{p-1}}.
\end{align*}
We also notice that
\begin{align*}
	& \left(A_{k+1}-A_{k}\right)\left(z_{k+1}-x_{k+1}\right)+A_{k}\left(x_{k}-x_{k+1}\right)+A_{k+1}\left(x_{k+1}-y_{k}\right)\\
	& =\left(A_{k+1}-A_{k}\right)z_{k+1}+A_{k}x_{k}-A_{k+1}y_{k}\\
	& =\left(A_{k+1}-A_{k}\right)\left(z_{k+1}-z_{k}\right)+\left(A_{k+1}-A_{k}\right)z_{k}+A_{k}x_{k}-A_{k+1}y_{k}\\
	& =\left(A_{k+1}-A_{k}\right)\left(z_{k+1}-z_{k}\right),
\end{align*}
where the last equality follows from $y_{k}=x_{k}+\frac{A_{k+1}-A_{k}}{A_{k+1}}\left(z_{k}-x_{k}\right)$. Therefore, 
\begin{align*}
	0 & \geq\mathcal{E}_{k+1}-\mathcal{E}_{k}\\
	& \quad+\left(A_{k+1}-A_{k}\right)\left\langle \nabla f\left(x_{k+1}\right),z_{k+1}-z_{k}\right\rangle \\
	& \quad+\frac{1+\mu A_{k}}{p}\left\Vert z_{k+1}-z_{k}\right\Vert ^{p}+MA_{k+1}s^{\frac{1}{p-1}}\left\Vert \nabla f\left(x_{k+1}\right)\right\Vert ^{\frac{p}{p-1}}.
\end{align*}
Now, we use the Fenchel-Young inequality $\left\langle s,u\right\rangle +\frac{1}{p}\left\Vert u\right\Vert ^{p}\geq-\frac{p-1}{p}\left\Vert s\right\Vert ^{\frac{p}{p-1}}$ with $u=\left(1+\mu A_{k}\right)^{\frac{1}{p}}\left(z_{k+1}-z_{k}\right)$
and $\ensuremath{s=\left(A_{k+1}-A_{k}\right)\left(1+\mu A_{k}\right)^{-\frac{1}{p}}\nabla f\left(x_{k+1}\right)}$ to obtain that
\begin{multline*}
	\left(A_{k+1}-A_{k}\right)\left\langle \nabla f\left(x_{k+1}\right),z_{k+1}-z_{k}\right\rangle +\frac{1+\mu A_{k}}{p}\left\Vert z_{k+1}-z_{k}\right\Vert ^{p}\\
	\geq-\frac{p-1}{p}\left(A_{k+1}-A_{k}\right)^{\frac{p}{p-1}}\left(1+\mu A_{k}\right)^{-\frac{1}{p-1}}\left\Vert \nabla f\left(x_{k+1}\right)\right\Vert ^{\frac{p}{p-1}}.
\end{multline*}
Hence, we have
\begin{align*}
	0 & \geq\mathcal{E}_{k+1}-\mathcal{E}_{k}\\
	& \quad+\left(MA_{k+1}s^{\frac{1}{p-1}}-\frac{p-1}{p}\left(A_{k+1}-A_{k}\right)^{\frac{p}{p-1}}\left(1+\mu A_{k}\right)^{-\frac{1}{p-1}}\right)\left\Vert \nabla f\left(x_{k+1}\right)\right\Vert ^{\frac{p}{p-1}}\\
	& =\mathcal{E}_{k+1}-\mathcal{E}_{k}\\
	& \quad+\left((p-1)p^{\frac{1}{p-1}}C^{\frac{1}{p-1}}A_{k+1}s^{\frac{1}{p-1}}-\frac{p-1}{p}\left(A_{k+1}-A_{k}\right)^{\frac{p}{p-1}}\left(1+\mu A_{k}\right)^{-\frac{1}{p-1}}\right)\left\Vert \nabla f\left(x_{k+1}\right)\right\Vert ^{\frac{p}{p-1}},
\end{align*}
where $C=\frac{1}{p}(\frac{M}{p-1})^{p-1}$. It is easy to see that the condition \eqref{eq:our_Ak} implies that the term
\[
\left((p-1)p^{\frac{1}{p-1}}C^{\frac{1}{p-1}}A_{k+1}s^{\frac{1}{p-1}}-\frac{p-1}{p}\left(A_{k+1}-A_{k}\right)^{\frac{p}{p-1}}\left(1+\mu A_{k}\right)^{-\frac{1}{p-1}}\right)\left\Vert \nabla f\left(x_{k+1}\right)\right\Vert ^{\frac{p}{p-1}}
\]
is non-negative. Thus, we conclude that
\[
0\geq\mathcal{E}_{k+1}-\mathcal{E}_{k}
\]
as desired.

\subsection{Lower bounds for the sequence $(A_k)$}
\label{app:lower_bound_Ak}

Let $(A_{k}^{\mathrm{best}})$ denote
 the sequence $(A_k)$ determined by \eqref{eq:our_specific_Ak}. In this section, we prove that that the following inequality holds:
\[
A_{k}^{\mathrm{best}}\geq \max\left\{ O\left(k^{p}\right),O\left(\left(1+C^{1/p}p\mu^{1/p}s^{1/p}\right)^{k}\right)\right\}. 
\]
 We use the following lemma.
\begin{lemma}
	\label{lem:Ak_lemma}
	For any sequence $(A_{k})$ satisfying $A_{0}=0$ and the condition \eqref{eq:our_Ak}, we have
	\begin{equation}
		\label{eq:A_best_is_best}
		A_{k}\leq A_{k}^{\mathrm{best}}\quad \forall k\geq0.
	\end{equation}
\end{lemma}
Its proof can be found in the following subsection. Now, we claim that the following two sequences satisfy the condition \eqref{eq:our_Ak}:
\[
A_{k}=Csk(k+1)\cdots(k+p-1)
\]
and
\[
A_{k}=\begin{cases}
	0, & k=0\\
	Cp^{p}s\left(1+C^{1/p}p\mu^{1/p}s^{1/p}\right)^{k-1} & k=1.
\end{cases}
\]
For the first sequence, we have
\begin{align*}
	& \left(A_{k+1}-A_{k}\right)^{p}-Cp^{p}sA_{k+1}^{p-1}\left(1+\mu A_{k}\right)\\
	& \leq\left(A_{k+1}-A_{k}\right)^{p}-Cp^{p}sA_{k+1}^{p-1}\\
	& =\left(Cps(k+1)\cdots(k+p-1)\right)^{p}-Cp^{p}s\left(Cs(k+1)\cdots(k+p)\right)^{p-1}\\
	& =C^{p}p^{p}s^{p}\left(\left((k+1)\cdots(k+p-1)\right)^{p}-\left((k+1)\cdots(k+p)\right)^{p-1}\right)\\
	& \leq0,
\end{align*}
which implies that \eqref{eq:our_Ak} holds.

For the second sequence, \eqref{eq:our_Ak} holds because
\begin{align*}
	& \left(A_{k+1}-A_{k}\right)^{p}-Cp^{p}sA_{k+1}^{p-1}\left(1+\mu A_{k}\right)\\
	& \leq\left(A_{k+1}-A_{k}\right)^{p}-C\mu p^{p}sA_{k+1}^{p-1}A_{k}\\
	& \leq\left(A_{k+1}-A_{k}\right)^{p}-C\mu p^{p}sA_{k}^{p}\\
	& =\left(\left(\frac{A_{k+1}}{A_{k}}-1\right)^{p}-C\mu p^{p}s\right)A_{k}^{p}\\
	& =\left(\left(C^{1/p}p\mu^{1/p}s^{1/p}\right)^{p}-C\mu p^{p}s\right)A_{k}^{p}\\
	& =0
\end{align*}
for all $k\geq1$ (the case $k=0$ is trivial). Thus, it follows from Lemma~\ref{lem:Ak_lemma} that
\begin{align*}
	A_{k}^{\mathrm{best}} & \geq\max\left\{ Csk(k+1)\cdots(k+p-1),Cp^{p}s\left(1+C^{1/p}p\mu^{1/p}s^{1/p}\right)^{k-1}\right\} \\
	& =\max\left\{ O\left(k^{p}\right),O\left(\left(1+C^{1/p}p\mu^{1/p}s^{1/p}\right)^{k}\right)\right\},
\end{align*}
as desired.

\subsubsection{Proof of Lemma~\ref{lem:Ak_lemma}}

For $r\geq0$, we define
\begin{align*}
	S(r) & :=\left\{ x:(x-r)^{p}-Cp^{p}sx^{p-1}(1+\mu r)\leq0\right\} \\
	U(r) & :=\max S_{r}.
\end{align*}
Then, it is straightforward to see the following:
\begin{itemize}
	\item The set $S(r)$ is nonempty. In particular, $r\in S(r)$ (which implies	$U(r)\geq r$).
	\item For any sequence $(A_{k})$ satisfying the condition \eqref{eq:our_Ak}, we have $A_{k+1}\in S(A_{k})$ for all $k\geq0$. 
	\item For the sequence $(A_{k})$ defined in \eqref{eq:our_specific_Ak}, we have $A_{k+1}=U(A_{k})$ for all $k\geq0$.
\end{itemize}
If we have
\begin{equation}
	\label{eq:Ur_monotonicity}
	U\left(r_{1}\right)\leq U\left(r_{2}\right)\textrm{ whenever }r_{1}\leq r_{2},
\end{equation}
then we can prove \eqref{eq:A_best_is_best} using mathematical induction on $k$. It clearly holds when $k=0$. If \eqref{eq:A_best_is_best} holds for $k$, then it holds for $k+1$ because
\[
A_{k+1}\leq U\left(A_{k}\right)\leq U(A_{k}^{\mathrm{best}})=A_{k+1}^{\mathrm{best}}.
\]

It remains to prove \eqref{eq:Ur_monotonicity}. Let $r_{1}$ and $r_{2}$ be positive real numbers with $r_{1}\leq r_{2}$. Then, it is easy to check that $r_{2}+U(r_{1})-r_{1}\in S(r_{2})$. Thus, we have
\[
U\left(r_{1}\right)\leq U\left(r_{1}\right)+\left(r_{2}-r_{1}\right)\leq U\left(r_{2}\right).
\]
This completes the proof.

\section{Existence and Uniqueness Theorems}
\label{app:existence_uniqueness}

\subsection{Proof of Theorem~\ref{thm:existence_unique_unified}}
\label{app:existence_uniqueness_unified_nag}

We prove a stronger result, that the unified Bregman Lagrangian flow \eqref{eq:unified_family} with $\alpha(t)=\log(\frac{2}{t}\cothc(\frac{\sqrt{\mu}}{2}t))$, $\beta(t)=\log(\frac{t^{2}}{4}\sinhc^{2}(\frac{\sqrt{\mu}}{2}t))$:
\begin{equation}
	\label{eq:mirror_u_ode}
	\begin{aligned}
		\dot{X} & =\frac{2}{t}\cothc\left(\frac{\sqrt{\mu}}{2}t\right)(Z-X)\\
		\frac{d}{dt}\nabla h(Z) & =\frac{t}{2}\tanhc\left(\frac{\sqrt{\mu}}{2}t\right)\left(\mu \nabla h(X)-\mu \nabla h(Z)-\nabla f(X)\right)
	\end{aligned}
\end{equation}
with the initial conditions $X(0)=Z(0)=x_0$ has a unique global solution $(X,Z)$ in $C^1([0,\infty),\mathbb{R}^n\times\mathbb{R}^n)$. Following \citep{krichene2015accelerated}, we assume that $\nabla f$ is $L_f$-Lipschitz continuous and $\nabla h$ is $L_h$-Lipschitz continuous. The strong convexity of $h$ implies a $L_{h^*}$-Lipschitz continuity of $\nabla h^*$ for some $L_{h^*}>0$ \citep[see][Proposition~12.60]{rockafellar2009variational}. 

\subsubsection{Proof of existence}

Fix $t_1>0$. We show the existence of solution to the system \eqref{eq:mirror_u_ode} on $[0,t_1]$. To remove the singularity of the system \eqref{eq:mirror_u_ode} at $t=0$, fix $\delta>0$, and consider the following system of ODEs:
\begin{equation}
	\label{eq:mirror_u_ode2}
	\begin{aligned}
		\dot{X} & =\frac{2}{\max\{\delta,t\}}\cothc\left(\frac{\sqrt{\mu}}{2}\max\{\delta,t\}\right)(Z-X)\\
		\frac{d}{dt}\nabla h(Z) & =\frac{t}{2}\tanhc\left(\frac{\sqrt{\mu}}{2}t\right)\left(\mu\nabla h(X)-\mu\nabla h(Z)-\nabla f(X)\right)
	\end{aligned}
\end{equation}
with $X(0)=Z(0)=x_0$, which does not have singularities. Denote the image of $Z$ under the mirror map as $W(t)=\nabla h(Z(t))$. Denote the convex conjugate of $h$ by $h^*:\mathbb{R}^n \rightarrow \mathbb{R}$. Then, $\nabla h$ and $\nabla h^*$ are inverses of each other \citep[see][Section~11]{rockafellar2009variational}. Now, we can equivalently write the system \eqref{eq:mirror_u_ode2} as
\begin{subequations}
	\label{eq:smoothed_amf_u}
	\begin{align}
		\label{eq:smoothed_amf_u_a}
		\dot{X} & =\frac{2}{\max\{\delta,t\}}\cothc\left(\frac{\sqrt{\mu}}{2}\max\{\delta,t\}\right)\left(\nabla h^{*}(W)-X\right)\\
		\label{eq:smoothed_amf_u_b}
		\dot{W} & =\frac{t}{2}\tanhc\left(\frac{\sqrt{\mu}}{2}t\right)\left(\mu\nabla h(X)-\mu W-\nabla f(X)\right)
	\end{align}
\end{subequations}
with $X(0)=x_0$ and $W(0)=w_0 := \nabla h\left(x_0\right)$. By the Cauchy-Lipschitz theorem, the system of ODEs \eqref{eq:smoothed_amf_u} has a unique solution $\left( X_\delta , W_\delta \right)$ in $C^1 ([0,t_1],\mathbb{R}^n\times\mathbb{R}^n)$. If we prove the following lemma, then one can prove the existence of solution to the ODE system \eqref{eq:mirror_u_ode2} following the argument in \citep[Section~3.2]{krichene2015accelerated}.
\begin{lemma}
	\label{lem:existence_lemma}
	Define a constant $T$ as
	\[
	T=\min\left\{ \sqrt{\frac{2}{\mu}},\frac{1}{2}\sqrt{\frac{1}{K_{2}K_{3}}}\right\} ,
	\]
	where $K_2$ and $K_3$ are constants defined in \eqref{eq:K_constants}. Then, the family of solutions $((X_{\delta},Z_{\delta})|_{[0,T]})_{\delta\in (0,T]}$ is equi-Lipschitz-continuous and uniformly bounded.
\end{lemma}
We now prove this lemma. We follow the argument of \citet{krichene2015accelerated} and omit the detailed calculations that can be found in \citep[Appendix~2]{krichene2015accelerated}. Fix $\delta$. For $t>0$, define
\begin{align*}
	A_{\delta}(t) & :=\sup_{u\in[0,t]}\frac{\left\Vert \dot{W}_{\delta}(u)\right\Vert }{u}\\
	B_{\delta}(t) & :=\sup_{u\in[0,t]}\frac{\left\Vert X_{\delta}(u)-x_{0}\right\Vert }{u}\\
	C_{\delta}(t) & :=\sup_{u\in[0,t]}\left\Vert \dot{X}_{\delta}(u)\right\Vert .
\end{align*}
Then, these quantities are finite. We first prove the following inequalities, which correspond to \citep[Lemma~3]{krichene2015accelerated}.
\begin{subequations}
	\label{eq:auxiliary_lemma_subeqs}
	\begin{align}
		\label{eq:auxiliary_lemma_subeq1}
		A_{\delta}(t) & \leq\mu\left\Vert w_{0}\right\Vert +\mu\left\Vert \nabla h\left(x_{0}\right)\right\Vert +\left\Vert \nabla f\left(x_{0}\right)\right\Vert +\left(\mu L_{h}+L_{f}\right)tB_{\delta}(t)\\
		\label{eq:auxiliary_lemma_subeq2}
		B_{\delta}(t) & \leq\frac{L_{h^{*}}t}{3}\cothc\left(\frac{\sqrt{\mu}}{2}T\right)A_{\delta}(t)\\
		\label{eq:auxiliary_lemma_subeq3}
		C_{\delta}(t) & \leq\cothc\left(\frac{\sqrt{\mu}}{2}T\right)\left(L_{h^{*}}TA_{\delta}(t)+2B_{\delta}(t)\right).
	\end{align}
\end{subequations}
\paragraph{Proof of \eqref{eq:auxiliary_lemma_subeq1}.}
Using $A_\delta$ and $B_\delta$, we can bound $\Vert W_{\delta}(t)-w_{0}\Vert$ and $\Vert X_{\delta}(t)-x_{0}\Vert$ as
\begin{align*}
	\left\Vert W_{\delta}(t)-w_{0}\right\Vert  & \leq\frac{t^{2}}{2}A_{\delta}(t)\\
	\left\Vert X_{\delta}(t)-x_{0}\right\Vert  & \leq tB_{\delta}(t).
\end{align*}
From \eqref{eq:smoothed_amf_u_b}, we have
\begin{align*}
	2\frac{\left\Vert \dot{W}_{\delta}(t)\right\Vert }{t} & =\tanhc\left(\frac{\sqrt{\mu}}{2}t\right)\left\Vert \mu\nabla h(X_{\delta})-\mu W_{\delta}-\nabla f(X_{\delta})\right\Vert \\
	& \leq\left\Vert \mu\nabla h(X_{\delta})-\mu W_{\delta}-\nabla f(X_{\delta})\right\Vert \\
	& \leq\mu\left\Vert W_{\delta}\right\Vert +\mu\left\Vert \nabla h(X_{\delta})\right\Vert +\left\Vert \nabla f(X_{\delta})\right\Vert \\
	& \leq\mu\left\Vert w_{0}\right\Vert +\frac{\mu t^{2}}{2}A_{\delta}(t)+\mu\left\Vert \nabla h\left(x_{0}\right)\right\Vert +\mu L_{h}tB_{\delta}(t)+\left\Vert \nabla f\left(x_{0}\right)\right\Vert +L_{f}tB_{\delta}(t).
\end{align*}
Thus,
\begin{align*}
	2A_{\delta}(t) & \leq\mu\left\Vert w_{0}\right\Vert +\mu\left\Vert \nabla h\left(x_{0}\right)\right\Vert +\left\Vert \nabla f\left(x_{0}\right)\right\Vert \\
	& \quad+\frac{\mu t^{2}}{2}A_{\delta}(t)+\left(\mu L_{h}+L_{f}\right)tB_{\delta}(t).
\end{align*}
Because $ T\leq \sqrt{{2}/{\mu}}$, we obtain the inequality \eqref{eq:auxiliary_lemma_subeq1}.

\paragraph{Proof of \eqref{eq:auxiliary_lemma_subeq2}.}

To bound the function $B_{\delta}(t) =\sup_{u\in[0,t]}\frac{\Vert X_{\delta}(u)-x_{0}\Vert }{u}$, we first compute an upper bound of $\left\Vert X_{\delta}(t)-x_{0}\right\Vert $ in the case $0\leq t\leq\delta$ and the case $t\geq\delta$ separately. First, consider the case $t\in [0,\delta]$. By \eqref{eq:smoothed_amf_u_a}, we have
\[
\dot{X}_{\delta}+\frac{2}{\delta}\cothc\left(\frac{\sqrt{\mu}}{2}\delta\right)\left(X_{\delta}-x_{0}\right)=\frac{2}{\delta}\cothc\left(\frac{\sqrt{\mu}}{2}\delta\right)\left(\nabla h^{*}(W_{\delta}-\nabla h^{*}\left(w_{0}\right)\right).
\]
Multiplying $e^{\frac{2}{\delta}\cothc\left(\frac{\sqrt{\mu}}{2}\delta\right)t}$, we obtain
\begin{multline*}
	e^{\frac{2}{\delta}\cothc\left(\frac{\sqrt{\mu}}{2}\delta\right)t}\left[\dot{X}_{\delta}+\frac{2}{\delta}\cothc\left(\frac{\sqrt{\mu}}{2}\delta\right)\left(X_{\delta}-x_{0}\right)\right]\\
	=\frac{2}{\delta}\cothc\left(\frac{\sqrt{\mu}}{2}\delta\right)e^{\frac{2}{\delta}\cothc\left(\frac{\sqrt{\mu}}{2}\delta\right)t}\left(\nabla h^{*}(W_{\delta})-\nabla h^{*}\left(w_{0}\right)\right).
\end{multline*}
This equality can be written as
\[
\frac{d}{dt}\left(\left(X_{\delta}(t)-x_{0}\right)e^{\frac{2}{\delta}\cothc\left(\frac{\sqrt{\mu}}{2}\delta\right)t}\right)=\frac{2}{\delta}\cothc\left(\frac{\sqrt{\mu}}{2}\delta\right)e^{\frac{2}{\delta}\cothc\left(\frac{\sqrt{\mu}}{2}\delta\right)t}\left(\nabla h^{*}\left(W_{\delta}(t)\right)-\nabla h^{*}\left(w_{0}\right)\right).
\]
Integrating both sides yields
\begin{multline*}
	\left(X_{\delta}(t)-x_{0}\right)e^{\frac{2}{\delta}\cothc\left(\frac{\sqrt{\mu}}{2}\delta\right)t}\\
	=\frac{2}{\delta}\cothc\left(\frac{\sqrt{\mu}}{2}\delta\right)\int_{0}^{t}\left[e^{\frac{2}{\delta}\cothc\left(\frac{\sqrt{\mu}}{2}\delta\right)s}\left(\nabla h^{*}\left(W_{\delta}(s)\right)-\nabla h^{*}\left(w_{0}\right)\right)\right]\,ds.
\end{multline*}
Taking norms, we have
\begin{align*}
	\left\Vert X_{\delta}(t)-x_{0}\right\Vert  & \leq\frac{2}{\delta}\cothc\left(\frac{\sqrt{\mu}}{2}\delta\right)\int_{0}^{t}\left\Vert \nabla h^{*}\left(W_{\delta}(s)\right)-\nabla h^{*}\left(w_{0}\right)\right\Vert \,ds\\
	& \leq\frac{2L_{h^{*}}}{\delta}\cothc\left(\frac{\sqrt{\mu}}{2}\delta\right)\int_{0}^{t}\left\Vert W_{\delta}(s)-w_{0}\right\Vert \,ds\\
	& \leq\frac{2L_{h^{*}}}{\delta}\cothc\left(\frac{\sqrt{\mu}}{2}\delta\right)\int_{0}^{t}\frac{s^{2}}{2}A_{\delta}(t)\,ds\\
	& =\frac{2L_{h^{*}}}{\delta}\cothc\left(\frac{\sqrt{\mu}}{2}\delta\right)A_{\delta}(t)\frac{t^{3}}{6}\\
	& \leq\frac{2L_{h^{*}}}{t}\cothc\left(\frac{\sqrt{\mu}}{2}\delta\right)A_{\delta}(t)\frac{t^{3}}{6}\\
	& =\frac{L_{h^{*}}t^{2}}{3}\cothc\left(\frac{\sqrt{\mu}}{2}\delta\right)A_{\delta}(t).
\end{align*}
So far, we provide an upper bound of $\left\Vert X_{\delta}(t)-x_{0}\right\Vert $ in the case $0\leq t\leq\delta$. We now consider the case $t\geq \delta$. By \eqref{eq:smoothed_amf_u_a}, we have
\[
\dot{X}_{\delta}+\frac{2}{t}\cothc\left(\frac{\sqrt{\mu}}{2}t\right)\left(X_{\delta}-x_{0}\right)=\frac{2}{t}\cothc\left(\frac{\sqrt{\mu}}{2}t\right)\left(\nabla h^{*}(W_{\delta})-\nabla h^{*}\left(w_{0}\right)\right).
\]
Multiplying $\frac{t^{2}}{4}\sinhc^{2}(\frac{\sqrt{\mu}}{2}t)$ to both sides, we obtain
\begin{multline*}
	\frac{t^{2}}{4}\sinhc^{2}\left(\frac{\sqrt{\mu}}{2}t\right)\dot{X}_{\delta}+\frac{t}{2}\sinhc\left(\frac{\sqrt{\mu}}{2}t\right)\cosh\left(\frac{\sqrt{\mu}}{2}t\right)\left(X_{\delta}-x_{0}\right)\\
	=\frac{t}{2}\sinhc\left(\frac{\sqrt{\mu}}{2}t\right)\cosh\left(\frac{\sqrt{\mu}}{2}t\right)\left(\nabla h^{*}(W_{\delta})-\nabla h^{*}\left(w_{0}\right)\right).
\end{multline*}
This equality can be written as
\begin{multline*}
	\frac{d}{dt}\left(\frac{t^{2}}{4}\sinhc^{2}\left(\frac{\sqrt{\mu}}{2}t\right)\left(X_{\delta}(t)-x_{0}\right)\right)\\
	=\frac{t}{2}\sinhc\left(\frac{\sqrt{\mu}}{2}t\right)\cosh\left(\frac{\sqrt{\mu}}{2}t\right)\left(\nabla h^{*}\left(W_{\delta}(t)\right)-\nabla h^{*}\left(w_{0}\right)\right).
\end{multline*}
Integrating both sides, we obtain
\begin{multline*}
	\frac{t^{2}}{4}\sinhc^{2}\left(\frac{\sqrt{\mu}}{2}t\right)\left(X_{\delta}(t)-x_{0}\right)\\
	=\int_{0}^{t}\left(\frac{s}{2}\sinhc\left(\frac{\sqrt{\mu}}{2}s\right)\cosh\left(\frac{\sqrt{\mu}}{2}s\right)\left(\nabla h^{*}\left(W_{\delta}(s)\right)-\nabla h^{*}\left(w_{0}\right)\right)\right)\,ds.
\end{multline*}
Taking norms, we have the following upper bound on $\left\Vert X_{\delta}(t)-x_{0}\right\Vert $:
\begin{align*}
	\left\Vert X_{\delta}(t)-x_{0}\right\Vert  & \leq\frac{2}{t}\cothc\left(\frac{\sqrt{\mu}}{2}t\right)\int_{0}^{t}\left\Vert \nabla h^{*}\left(W_{\delta}(s)\right)-\nabla h^{*}\left(w_{0}\right)\right\Vert \,ds.\\
	& \leq\frac{2L_{h^{*}}}{t}\cothc\left(\frac{\sqrt{\mu}}{2}t\right)\int_{0}^{t}\left\Vert W_{\delta}(s)-w_{0}\right\Vert \,ds.\\
	& \leq\frac{2L_{h^{*}}}{t}\cothc\left(\frac{\sqrt{\mu}}{2}t\right)\int_{0}^{t}\frac{s^{2}}{2}A_{\delta}(t)\,ds\\
	& =\frac{2L_{h^{*}}}{t}\cothc\left(\frac{\sqrt{\mu}}{2}t\right)A_{\delta}(t)\frac{t^{3}}{6}\\
	& =\frac{L_{h^{*}}t^{2}}{3}\cothc\left(\frac{\sqrt{\mu}}{2}t\right)A_{\delta}(t).
\end{align*}
Combining both cases $0\leq t\leq\delta$ and $t\geq\delta$, we have
\[
\left\Vert X_{\delta}(t)-x_{0}\right\Vert  \leq \frac{L_{h^{*}}t^{2}}{3}\cothc\left(\frac{\sqrt{\mu}}{2}T\right)A_{\delta}(t)
\]
for all $t\geq0$. Dividing by $t$ and taking the supremum, we obtain
\[
B_{\delta}(t)\leq\frac{L_{h^{*}}t}{3}\cothc\left(\frac{\sqrt{\mu}}{2}T\right)A_{\delta}(t).
\]
\paragraph{Proof of \eqref{eq:auxiliary_lemma_subeq3}.}
By \eqref{eq:smoothed_amf_u_a}, we have
\begin{align*}
	\left\Vert \dot{X}\right\Vert  & =\frac{2}{\max\{\delta,t\}}\cothc\left(\frac{\sqrt{\mu}}{2}\max\{\delta,t\}\right)\left\Vert \nabla h^{*}\left(W_{\delta}(t)\right)-X_{\delta}(t)\right\Vert \\
	& \leq\frac{2}{\max\{\delta,t\}}\cothc\left(\frac{\sqrt{\mu}}{2}\max\{\delta,t\}\right)\left(\left\Vert \nabla h^{*}\left(W_{\delta}(t)\right)-\nabla h^{*}\left(z_{0}\right)\right\Vert +\left\Vert X_{\delta}(t)-x_{0}\right\Vert \right)\\
	& \leq\frac{2}{\max\{\delta,t\}}\cothc\left(\frac{\sqrt{\mu}}{2}\max\{\delta,t\}\right)\left(\frac{t^{2}}{2}L_{h^{*}}A_{\delta}(t)+tB_{\delta}(t)\right)\\
	& \leq\cothc\left(\frac{\sqrt{\mu}}{2}T\right)\frac{2}{t}\left(\frac{t^{2}}{2}L_{h^{*}}A_{\delta}(t)+tB_{\delta}(t)\right)\\
	& \leq\cothc\left(\frac{\sqrt{\mu}}{2}T\right)\left(L_{h^{*}}TA_{\delta}(t)+2B_{\delta}(t)\right).
\end{align*}

\paragraph{Complete the proof of Lemma~\ref{lem:existence_lemma}.}

Define five positive constants $K_1$, $\ldots$, $K_5$ as
\begin{equation}
	\label{eq:K_constants}
	\begin{aligned}
		K_{1} & :=\mu\left\Vert w_{0}\right\Vert +\mu\left\Vert \nabla h\left(x_{0}\right)\right\Vert +\left\Vert \nabla f\left(x_{0}\right)\right\Vert \\
		K_{2} & :=\mu L_{h}+L_{f}\\
		K_{3} & :=\frac{2L_{h^{*}}}{3}\\
		K_{4} & :=2L_{h^{*}}\\
		K_{5} & :=4.
	\end{aligned}
\end{equation}
Because $T\leq\frac{2}{\sqrt{\mu}}$, we have $\cothc(\frac{\sqrt{\mu}}{2}T)\leq\cothc(1)\leq2$.
Thus, the inequalities \eqref{eq:auxiliary_lemma_subeqs} imply
\begin{subequations}
	\label{eq:a_a_lemma_subeqs}
	\begin{align}
		\label{eq:a_a_lemma_subeq1}
		A_{\delta}(t) & \leq K_{1}+K_{2}TB_{\delta}(t)\\
		\label{eq:a_a_lemma_subeq2}
		B_{\delta}(t) & \leq K_{3}TA_{\delta}(t)\\
		\label{eq:a_a_lemma_subeq3}
		C_{\delta}(t) & \leq K_{4}TA_{\delta}(t)+K_{5}B_{\delta}(t).
	\end{align}
\end{subequations}
Combining \eqref{eq:a_a_lemma_subeq1} and \eqref{eq:a_a_lemma_subeq2}, we have
\[
\left(\frac{1}{K_{3}T}-K_{2}T\right)B_{\delta}(t)\leq K_{1}.
\]
Because $T\mapsto\frac{1}{K_{3}T}-K_{2}T$
is a positive decreasing funtion on $[0,\frac{1}{2}\sqrt{\frac{1}{K_{2}K_{3}}}]$ and $T\leq\frac{1}{2}\sqrt{\frac{1}{K_{2}K_{3}}}$, we have
\begin{equation}
	\label{eq:a_a_result1}
	B_{\delta}(T)\leq\left(\frac{1}{K_{3}\cdot \frac{1}{2}\sqrt{\frac{1}{K_{2}K_{3}}}}-K_{2}\cdot \frac{1}{2}\sqrt{\frac{1}{K_{2}K_{3}}} \right)^{-1}K_{1}=\frac{2}{3}K_{1}\sqrt{\frac{K_{3}}{K_{2}}}.
\end{equation}
The inequalities \eqref{eq:a_a_lemma_subeq1}, \eqref{eq:a_a_result1}, and $T\leq\frac{1}{2}\sqrt{\frac{1}{K_{2}K_{3}}}$ imply
\begin{equation}
	\label{eq:a_a_result2}
	A_{\delta}(T)\leq K_{1}+K_{2}TB_{\delta}(T)\leq K_{1}+K_{2}\left(\frac{1}{2}\sqrt{\frac{1}{K_{2}K_{3}}}\right)\left(\frac{2}{3}K_{1}\sqrt{\frac{K_{3}}{K_{2}}}\right).
\end{equation}
The inequalities \eqref{eq:a_a_lemma_subeq1}, \eqref{eq:a_a_result1}, \eqref{eq:a_a_result2}, and $T\leq\frac{1}{2}\sqrt{\frac{1}{K_{2}K_{3}}}$ imply
\begin{equation}
	\label{eq:a_a_result3}
	\begin{aligned}
		C_{\delta}(T) & \leq K_{4}TA_{\delta}(T)+K_{5}B_{\delta}(T)\\
		& \leq K_{4}\left(\frac{1}{2}\sqrt{\frac{1}{K_{2}K_{3}}}\right)\left(K_{1}+K_{2}\left(\frac{1}{2}\sqrt{\frac{1}{K_{2}K_{3}}}\right)\left(\frac{2}{3}K_{1}\sqrt{\frac{K_{3}}{K_{2}}}\right)\right)+K_{5}\left(\frac{2}{3}K_{1}\sqrt{\frac{K_{3}}{K_{2}}}\right).
	\end{aligned}
\end{equation}
Therefore, $\|\dot{W}\|$ and $\|\dot{X}\|$ are bounded uniformly in $\delta$ because
\begin{align*}
	\left\Vert \dot{W}_{\delta}(t)\right\Vert  & \leq TA_{\delta}(T)\\
	\left\Vert \dot{X}_{\delta}(t)\right\Vert  & \leq C_{\delta}(T)
\end{align*}
for all $t\in[0,T]$. This implies that the family of solutions $((X_{\delta},Z_{\delta})|_{[0,T]})_{\delta\in (0,T]}$ is equi-Lipschitz-continuous and uniformly bounded.

\subsubsection{Proof of uniqueness}

We follow the argument in \citep[Appendix~3]{krichene2015accelerated} and omit the detailed calculations that can be found in \citep{krichene2015accelerated}.
Because we only need to prove the uniqueness of solution near $t=0$, we assume $t<T$ for some $T>0$. Let $(X,W)$ and $\left(\bar{X},\bar{W}\right)$ be solutions to the following system of ODEs, which is equivalent to \eqref{eq:mirror_u_ode}:
\begin{align*}
	\dot{X} & =\frac{2}{t}\cothc\left(\frac{\sqrt{\mu}}{2}t\right)\left(\nabla h^{*}(W)-X\right)\\
	\dot{W} & =\frac{t}{2}\tanhc\left(\frac{\sqrt{\mu}}{2}t\right)\left(\mu\nabla h(X)-\mu W-\nabla f(X)\right).
\end{align*}
Let $\Delta_{W}=W-\bar{W}$ and $\Delta_{X}=X-\bar{X}$. Then, we have
\begin{align*}
	\dot{\Delta}_{W} & =\frac{t}{2}\tanhc\left(\frac{\sqrt{\mu}}{2}t\right)\left(\mu\nabla h(X)-\mu W-\nabla f(X)-\mu\nabla h\left(\bar{X}\right)+\mu\bar{W}+\nabla f\left(\bar{X}\right)\right)\\
	\dot{\Delta}_{X} & =\frac{2}{t}\cothc\left(\frac{\sqrt{\mu}}{2}t\right)\left(\nabla h^{*}(W)-\nabla h^{*}\left(\bar{W}\right)-\Delta_{X}\right)
\end{align*}
with $\Delta_{X}(0)=\Delta_{W}(0)=0$. Define
\begin{align*}
	A(t) & :=\sup_{[0,t]}\frac{\left\Vert \dot{\Delta}_{W}(u)\right\Vert }{u}\\
	B(t) & :=\sup_{[0,t]}\left\Vert \Delta_{X}\right\Vert .
\end{align*}
Then, $B(t)$ and $C(t)$ are finite because $\Delta_{X}$ and $\Delta_W$ are continuous. First, we compute an upper bound of $A(t)$. We have
\begin{equation}
	\label{eq:unique1}
	\begin{aligned}
		\left\Vert \dot{\Delta}_{W}(t)\right\Vert  & =\frac{t}{2}\tanhc\left(\frac{\sqrt{\mu}}{2}t\right)\left\Vert \mu\nabla h(X)-\mu W-\nabla f(X)-\mu\nabla h\left(\bar{X}\right)+\mu\bar{W}+\nabla f\left(\bar{X}\right)\right\Vert \\
		& \leq\frac{t}{2}\tanhc\left(\frac{\sqrt{\mu}}{2}t\right)\left(\mu\left\Vert \nabla h(X)-\nabla h\left(\bar{X}\right)\right\Vert +\mu\left\Vert W-\bar{W}\right\Vert +\left\Vert \nabla f(X)-\nabla f\left(\bar{X}\right)\right\Vert \right)\\
		& \leq\frac{t}{2}\tanhc\left(\frac{\sqrt{\mu}}{2}t\right)\left(\left(\mu L_{h}+L_{f}\right)\left\Vert \Delta_{X}\right\Vert +\mu\left\Vert \Delta_{W}\right\Vert \right)\\
		& \leq\frac{t}{2}\tanhc\left(\frac{\sqrt{\mu}}{2}t\right)\left(\left(\mu L_{h}+L_{f}\right)B(t)+\frac{\mu t^{2}}{2}A(t)\right),
	\end{aligned}
\end{equation}
where we used $\Vert \Delta_{W}(t)\Vert \leq\Vert \int_{0}^{t}\dot{\Delta}_{W}(s)\,ds\Vert \leq\int_{0}^{t}sA(s)\,ds\leq\int_{0}^{t}sA(t)\,ds=\frac{t^{2}}{2}A(t)$ for the last inequality. Dividing both sides of \eqref{eq:unique1} by $t$ and then taking the supremum, we obtain
\begin{equation}
	\label{eq:unique2}
	A(t)\leq\frac{1}{2}\tanhc\left(\frac{\sqrt{\mu}}{2}t\right)\left(\left(\mu L_{h}+L_{f}\right)B(t)+\frac{\mu t^{2}}{2}A(t)\right).
\end{equation}
Nest, we compute an upper boudn of $B(t)$. We have
\[
\dot{\Delta}_{X}+\frac{2}{t}\cothc\left(\frac{\sqrt{\mu}}{2}t\right)\Delta_{X}=\frac{2}{t}\cothc\left(\frac{\sqrt{\mu}}{2}t\right)\left(\nabla h^{*}(W)-\nabla h^{*}\left(\bar{W}\right)\right).
\]
Multiplying both sides by $\frac{t^{2}}{4}\sinhc^{2}\left(\frac{\sqrt{\mu}}{2}t\right)$, we have
\begin{multline*}
	\frac{t^{2}}{4}\sinhc^{2}\left(\frac{\sqrt{\mu}}{2}t\right)\dot{\Delta}_{X}+\frac{t}{2}\sinhc\left(\frac{\sqrt{\mu}}{2}t\right)\cosh\left(\frac{\sqrt{\mu}}{2}t\right)\Delta_{X}\\
	=\frac{t}{2}\sinhc\left(\frac{\sqrt{\mu}}{2}t\right)\cosh\left(\frac{\sqrt{\mu}}{2}t\right)\left(\nabla h^{*}(W)-\nabla h^{*}\left(\bar{W}\right)\right).
\end{multline*}
This equality can be written as
\[
\frac{d}{dt}\left(\frac{t^{2}}{4}\sinhc^{2}\left(\frac{\sqrt{\mu}}{2}t\right)\Delta_{X}\right)=\frac{t}{2}\sinhc\left(\frac{\sqrt{\mu}}{2}t\right)\cosh\left(\frac{\sqrt{\mu}}{2}t\right)\left(\nabla h^{*}(W)-\nabla h^{*}\left(\bar{W}\right)\right).
\]
Integrating both sides, we obtain
\[
\frac{t^{2}}{4}\sinhc^{2}\left(\frac{\sqrt{\mu}}{2}t\right)\Delta_{X}=\int_{0}^{t}\left[\frac{s}{2}\sinhc\left(\frac{\sqrt{\mu}}{2}s\right)\cosh\left(\frac{\sqrt{\mu}}{2}s\right)\left(\nabla h^{*}(W(s))-\nabla h^{*}\left(\bar{W}(s)\right)\right)\right]\,ds.
\]
Taking norms, we have
\begin{align*}
	\left\Vert \Delta_{X}(t)\right\Vert  & \leq\frac{2}{t}\cothc\left(\frac{\sqrt{\mu}}{2}t\right)\int_{0}^{t}\left\Vert \nabla h^{*}(W(s))-\nabla h^{*}\left(\bar{W}(s)\right)\right\Vert \,ds.\\
	& \leq\frac{2L_{h^{*}}}{t}\cothc\left(\frac{\sqrt{\mu}}{2}t\right)\int_{0}^{t}\left\Vert \Delta_{W}(s)\right\Vert \,ds.\\
	& \leq\frac{2L_{h^{*}}}{t}\cothc\left(\frac{\sqrt{\mu}}{2}t\right)\int_{0}^{t}\frac{s^{2}}{2}A(t)\,ds\\
	& =\frac{L_{h^{*}}2}{t}\cothc\left(\frac{\sqrt{\mu}}{2}t\right)A(t)\frac{t^{3}}{6}\\
	& =\frac{L_{h^{*}}t^{2}}{3}\cothc\left(\frac{\sqrt{\mu}}{2}t\right)A(t).
\end{align*}
Taking the supremum yields
\begin{equation}
	\label{eq:unique3}
	B(t)\leq\frac{L_{h^{*}}t^{2}}{3}\cothc\left(\frac{\sqrt{\mu}}{2}t\right)A(t).
\end{equation}
Now, combining the inequalities \eqref{eq:unique2} and \eqref{eq:unique3}, we have
\begin{align*}
	A(t) & \leq\frac{1}{2}\tanhc\left(\frac{\sqrt{\mu}}{2}t\right)\left(\left(\mu L_{h}+L_{f}\right)B(t)+\frac{\mu t^{2}}{2}A(t)\right)\\
	& \leq\frac{1}{2}\tanhc\left(\frac{\sqrt{\mu}}{2}t\right)\left(\left(\mu L_{h}+L_{f}\right)\frac{L_{h^{*}}t^{2}}{3}\cothc\left(\frac{\sqrt{\mu}}{2}t\right)+\frac{\mu t^{2}}{2}\right)A(t).
\end{align*}
Using continuity, it is easy to see that there is $T_{\mathrm{small}}>0$ such that the following inequality
holds whenever $t\in(0,T_{\mathrm{small}})$: 
\[
\frac{1}{2}\tanhc\left(\frac{\sqrt{\mu}}{2}t\right)\left(\left(\mu L_{h}+L_{f}\right)\frac{L_{h^{*}}t^{2}}{3}\cothc\left(\frac{\sqrt{\mu}}{2}t\right)+\frac{\mu t^{2}}{2}\right)<1.
\]
Thus, for $t\in(0,T_{\mathrm{small}})$, we have $A(t)\leq 1\cdot A(t)$, which implies $A(t)=0$ because $A(t)$ is nonnegative by its definition. Finally, $B(t)=0$ follows from \eqref{eq:unique3}. This completes the proof.

\subsection{Existence and uniqueness of solution to the unified accelerated tensor flow}
\label{app:existence_uniqueness_hyperbolic}

We first note that
\begin{itemize}
	\item The system of ODEs \eqref{eq:mirror_u_ode} is the unified Bregman Lagrangian flow \eqref{eq:unified_family} with $\beta_1=\log\left(\frac{t^{2}}{4}\sinhc^{2}\left(\frac{\sqrt{\mu}}{2}t\right)\right)$ and $\alpha_1=\log \dot{\beta}_1$.
	\item The unified accelerated tensor flow \eqref{eq:hyperbolic_flow} is the unified Bregman Lagrangian flow \eqref{eq:unified_family} with $\beta_2 = p\log t+\log C+p\log\left(\sinhc_{p}\left(C^{1/p}\mu^{1/p}t\right)\right)$ and $\alpha_2 =\log \dot{\beta}_2$.
\end{itemize}
Define a function $\bfT:[0,\infty)\to[0,\infty)$ as $\bfT = \beta_1^{-1} \circ \beta_2$. Then, we have
\begin{align*}
	{\alpha_2}(t) & =\alpha_1(\bfT(t))+\log\dot{\bfT}(t)\\
	{\beta_2}(t) & =\beta_1(\bfT(t)).
\end{align*}
Thus, by Theorem~\ref{thm:time_dilation}, if $(X_1,Z_1)$ is a solution to the unified NAG system, then $X_2(t)=X_1(\bfT(t))$ and $Z_2(t)=Z_1(\bfT(t))$ is a solution to the unified accelerated tensor system. Thus, the existence of solution to the unified NAG system implies the existence of solution to the unified accelerated tensor system.

A similar argument shows that if $(X_2,Z_2)$ is a solution to the unified accelerated tensor system, then $X_1(t)=X_2(\bfT^{-1}(t))$ and $Z_1(t)=Z_2(\bfT^{-1}(t))$ is a solution to the unified NAG system. It is easy to show that this correspondence is one-to-one. Thus, the uniqueness of solution to the unified NAG system implies the uniqueness of solution to the unified accelerated tensor system.

\section{Further Exploration: ODE Model for Minimizing Gradient Norms of Strongly Convex Functions}

\subsection{Limiting ODE of OGM}
\label{app:limiting_ogm}

For the sequence $\theta_k$ defined in \eqref{eq:theta_sequence}, \citet{su2016differential} showed that the algorithm 
\begin{equation}
	\label{eq:theta_nag-c}
	\begin{aligned}
		y_{k} & =\left(1-\frac{1}{\theta_{k}}\right)x_{k}+\frac{1}{\theta_{k}}z_{k}\\
		x_{k+1} & =y_{k}-s\nabla f\left(y_{k}\right)\\
		z_{k+1} & =z_{k}-s\theta_{k}\nabla f\left(y_{k}\right)
	\end{aligned}
\end{equation}
converges to \nagcsp ODE as $s\to0$ \citep[see][Section~2]{su2016differential} (in fact, this algorithm is equivalent to the original NAG \eqref{eq:original_nag} with $\mu=0$ and $\gamma_0=\infty$). Because $\Vert x_{k+1}-y_{k}\Vert=o(\sqrt{s})$, we can ignore the gradient descent step $x_{k+1}=y_{k}-s\nabla f\left(y_{k}\right)$ in both \eqref{eq:theta_nag-c} and OGM. Then, applying OGM to the objective function $f$ is equivalent to applying the algorithm \eqref{eq:theta_nag-c} to the objective function $2f$. Thus, the limiting ODE of OGM is given by
\[
\ddot{X}+\frac{3}{t}\dot{X}+2\nabla f(X)=0.
\]

\subsection{Proof of Theorem~\ref{thm:mainthm_unified_ogm-g}}
\label{app:unified_ogm-g}

For convenience, we assume $\mu>0$ (the case $\mu=0$ can be handled easily). We denote $X:=X(t)$ and $x^{T}:=X(T)$. We also omit the input $\frac{\sqrt{\mu}}{2}(T-t)$ of each hyperbolic function. For example, we  write the unified NAG-G ODE \eqref{eq:unified_ogm-g} as
\[
\ddot{X}+\left(\frac{\sqrt{\mu}}{2}\tanh+\frac{3\sqrt{\mu}}{2}\coth\right)\dot{X}+\nabla f(X)=0
\]
and the continuous-time energy function \eqref{eq:energy_ogm-g} as
\[
\mathcal{E}(t)=\mu^2\csch^{4}\left(\frac{\sinh^{2}}{\mu}\left(f(X)-f\left(x^{T}\right)\right)-\frac{1}{2}\left\Vert X-x^{T}\right\Vert ^{2}+\frac{\cosh^{2}}{2}\left\Vert X+\frac{\tanh}{\sqrt{\mu}}\dot{X}-x^{T}\right\Vert ^{2}\right).
\]
Then, we have
\begin{align*}
	& \frac{\sinh^{4}}{\mu^2}\dot{\mathcal{E}}(t) \\
	& =\sinh^{4}\frac{d}{dt}\left\{ \csch^{4}\right\} \left(\frac{\sinh^{2}}{\mu}\left(f(X)-f\left(x^{T}\right)\right)-\frac{1}{2}\left\Vert X-x^{T}\right\Vert ^{2}+\frac{\cosh^{2}}{2}\left\Vert X+\frac{\tanh}{\sqrt{\mu}}\dot{X}-x^{T}\right\Vert ^{2}\right)\\
	& \quad+\frac{d}{dt}\left\{ \frac{\sinh^{2}}{\mu}\left(f(X)-f\left(x^{T}\right)\right)-\frac{1}{2}\left\Vert X-x^{T}\right\Vert ^{2}+\frac{\cosh^{2}}{2}\left\Vert X+\frac{\tanh}{\sqrt{\mu}}\dot{X}-x^{T}\right\Vert ^{2}\right\} \\
	& =2\sqrt{\mu}\coth\left(\frac{\sinh^{2}}{\mu}\left(f(X)-f\left(x^{T}\right)\right)-\frac{1}{2}\left\Vert X-x^{T}\right\Vert ^{2}+\frac{\cosh^{2}}{2}\left\Vert X+\frac{\tanh}{\sqrt{\mu}}\dot{X}-x^{T}\right\Vert ^{2}\right)\\
	& \quad-\frac{\sinh\cosh}{\sqrt{\mu}}\left(f(X)-f\left(x^{T}\right)\right)+\frac{\sinh^{2}}{\mu}\left\langle \nabla f(X),\dot{X}\right\rangle -\left\langle X-x^{T},\dot{X}\right\rangle \\
	& \quad-\frac{\sqrt{\mu}\sinh\cosh}{2}\left\Vert X+\frac{\tanh}{\sqrt{\mu}}\dot{X}-x^{T}\right\Vert ^{2}+\cosh^{2}\left\langle X+\frac{\tanh}{\sqrt{\mu}}\dot{X}-x^{T},-\dot{X}-\frac{\tanh}{\sqrt{\mu}}\nabla f(X)\right\rangle ,
\end{align*}
where we used 
\begin{align*}
	\frac{d}{dt}\left\{ X+\frac{\tanh}{\sqrt{\mu}}\dot{X}-x^{T}\right\}  & =\frac{\tanh}{\sqrt{\mu}}\ddot{X}+\left(1-\frac{1}{2}\sech^{2}\right)\dot{X}\\
	& =\left(-\frac{1}{2}\tanh^{2}-\frac{1}{2}-\frac{1}{2}\sech^{2}\right)\dot{X}-\frac{\tanh}{\sqrt{\mu}}\nabla f(X)\\
	& =-\dot{X}-\frac{\tanh}{\sqrt{\mu}}\nabla f(X)
\end{align*}
for the last equality. We further simplify as
\begin{align*}
	\frac{\sinh^{4}}{\mu^2}\dot{\mathcal{E}}(t) & =\frac{2\sinh\cosh}{\sqrt{\mu}}\left(f(X)-f\left(x^{T}\right)\right)-\sqrt{\mu}\coth\left\Vert X-x^{T}\right\Vert ^{2}\\
	& \quad+\sqrt{\mu}\coth\cosh^{2}\left(\left\Vert X-x^{T}\right\Vert ^{2}+\frac{\tanh^{2}}{\mu}\left\Vert \dot{X}\right\Vert ^{2}+\frac{2\tanh}{\sqrt{\mu}}\left\langle X-x^{T},\dot{X}\right\rangle \right)\\
	& \quad-\frac{\sinh\cosh}{\sqrt{\mu}}\left(f(X)-f\left(x^{T}\right)\right)+\frac{\sinh^{2}}{\mu}\left\langle \nabla f(X),\dot{X}\right\rangle -\left\langle X-x^{T},\dot{X}\right\rangle \\
	& \quad-\frac{\sqrt{\mu}\sinh\cosh}{2}\left(\left\Vert X-x^{T}\right\Vert ^{2}+\frac{\tanh^{2}}{\mu}\left\Vert \dot{X}\right\Vert ^{2}+\frac{2\tanh}{\sqrt{\mu}}\left\langle X-x^{T},\dot{X}\right\rangle \right)\\
	& \quad-\cosh^{2}\bigg(\left\langle X-x^{T},\dot{X}\right\rangle +\frac{\tanh}{\sqrt{\mu}}\left\Vert \dot{X}\right\Vert ^{2}\\
	& \qquad \quad +\frac{\tanh}{\sqrt{\mu}}\left\langle X-x^{T},\nabla f(X)\right\rangle +\frac{\tanh^{2}}{\mu}\left\langle \dot{X},\nabla f(X)\right\rangle \bigg)\\
	& =\left(\frac{2\sinh\cosh}{\sqrt{\mu}}-\frac{\sinh\cosh}{\sqrt{\mu}}\right)\left(f(X)-f\left(x^{T}\right)\right)\\
	& \quad+\left(-\sqrt{\mu}\coth+\sqrt{\mu}\coth\cosh^{2}-\frac{\sqrt{\mu}\sinh\cosh}{2}\right)\left\Vert X-x^{T}\right\Vert ^{2}\\
	& \quad+\left(\frac{\sinh\cosh}{\sqrt{\mu}}-\frac{\sinh^{2}\tanh}{2\sqrt{\mu}}-\frac{\sinh\cosh}{\sqrt{\mu}}\right)\left\Vert \dot{X}\right\Vert ^{2}\\
	& \quad+\left(2\cosh^{2}-1-\sinh^{2}-\cosh^{2}\right)\left\langle X-x^{T},\dot{X}\right\rangle \\
	& \quad+\left(\frac{\sinh^{2}}{\mu}-\frac{\sinh^{2}}{\mu}\right)\left\langle \nabla f(X),\dot{X}\right\rangle \\
	& \quad-\frac{\sinh\cosh}{\sqrt{\mu}}\left\langle X-x^{T},\nabla f(X)\right\rangle \\
	& =\frac{\sinh\cosh}{\sqrt{\mu}}\left(f(X)-f\left(x^{T}\right)\right)\\
	& \quad+\frac{\sqrt{\mu}\sinh\cosh}{2}\left\Vert X-x^{T}\right\Vert ^{2}-\frac{\sinh^{2}\tanh}{2\sqrt{\mu}}\left\Vert \dot{X}\right\Vert ^{2}-\frac{\sinh\cosh}{\sqrt{\mu}}\left\langle X-x^{T},\nabla f(X)\right\rangle .
\end{align*}
It follows from the $\mu$-strong convexity of $f$ that $f(X)-f\left(x^{T}\right)\leq\left\langle X-x^{T},\nabla f(X)\right\rangle -\frac{\mu}{2}\left\Vert X-x^{T}\right\Vert ^{2}$.
Thus, we have
\begin{align*}
	\frac{\sinh^{4}}{\mu^2}\dot{\mathcal{E}}(t) & \leq\frac{\sinh\cosh}{\sqrt{\mu}}\left(\left\langle X-x^{T},\nabla f(X)\right\rangle -\frac{\mu}{2}\left\Vert X-x^{T}\right\Vert ^{2}\right)\\
	& \quad+\frac{\sqrt{\mu}\sinh\cosh}{2}\left\Vert X-x^{T}\right\Vert ^{2}-\frac{\sinh^{2}\tanh}{2\sqrt{\mu}}\left\Vert \dot{X}\right\Vert ^{2}-\frac{\sinh\cosh}{\sqrt{\mu}}\left\langle X-x^{T},\nabla f(X)\right\rangle \\
	& =-\frac{\sinh^{2}\tanh}{2\sqrt{\mu}}\left\Vert \dot{X}\right\Vert ^{2}\\
	& \leq0.
\end{align*}

\subsection{Computing $\dot{X}(T)$ and $\ddot{X}(T)$}
\label{app:dotX_and_ddotX}

For simplicity, we assume that the limits $\lim_{t\rightarrow T^{-}}\dot{X}(T)$ and $\lim_{t\rightarrow T^{-}}\ddot{X}(T)$ exist.\footnote{The proof to prove the existence of these limits is similar to that in \citep[Appendix~D.3]{suh2022continuous}, so we omit it.} Consider the energy function
\begin{multline}
	\mathcal{E}(t)=\frac{1}{2}\left\Vert \dot{X}(t)\right\Vert ^{2}+\left(f(X(t))-f\left(x^{*}\right)\right)\\
	+\int_{0}^{t}\left[\frac{\sqrt{\mu}}{2}\tanh\left(\frac{\sqrt{\mu}}{2}(T-s)\right)+\frac{3}{T-s}\cothc\left(\frac{\sqrt{\mu}}{2}(T-s)\right)\right]\left\Vert \dot{X}(s)\right\Vert ^{2}ds.
\end{multline}
Then, it is easy to show that $\mathcal{E}(t)=\mathcal{E}(0)$ for all $t\in[0,T)$. Because the terms $\frac{1}{2}\left\Vert \dot{X}(t)\right\Vert ^{2}$ and $f(X(t))-f\left(x^{*}\right)$ are non-negtive, we have
\[
\int_{0}^{T}\left[\frac{\sqrt{\mu}}{2}\tanh\left(\frac{\sqrt{\mu}}{2}(T-s)\right)+\frac{3}{T-s}\cothc\left(\frac{\sqrt{\mu}}{2}(T-s)\right)\right]\left\Vert \dot{X}(s)\right\Vert ^{2}ds<\infty.
\]
This implies $\lim_{t\rightarrow T^{-}} \dot{X}(t)=0$. By L'H\^{o}pital's rule, we obtain that
\[
\lim_{t\rightarrow T^{-}}\left[\frac{\sqrt{\mu}}{2}\tanh\left(\frac{\sqrt{\mu}}{2}(T-t)\right)+\frac{3}{T-t}\cothc\left(\frac{\sqrt{\mu}}{2}(T-t)\right)\right]\dot{X}(t)=-3\ddot{X}(T).
\]
Now, we have
\begin{align*}
	0 & =\lim_{t\rightarrow T^{-}}\bigg\{\ddot{X}(t)+\left[\frac{\sqrt{\mu}}{2}\tanh\left(\frac{\sqrt{\mu}}{2}(T-t)\right)+\frac{3}{T-t}\cothc\left(\frac{\sqrt{\mu}}{2}(T-t)\right)\right]\dot{X}(t)+\nabla f(X(t))\bigg\}\\
	& =-2\ddot{X}(T)+\nabla f(X(T)).
\end{align*}
Thus, $\ddot{X}(T)=\frac{1}{2}\nabla f(X(T))$.

\vskip 0.2in
\bibliography{unag_references}

\end{document}